\documentclass[11pt,a4paper,leqno]{article}
\usepackage{a4wide} 
\usepackage{latexsym}
\usepackage{amsmath}
\usepackage{amssymb} 
\usepackage{stackrel}  
\usepackage{color} 
\usepackage{graphicx}
\usepackage[linesnumbered,ruled,vlined]{algorithm2e}

\usepackage{authblk}

\usepackage{hyperref}

\newtheorem{defin}{Definition}
\newtheorem{lemma}{Lemma}
\newtheorem{prop}{Proposition}
\newtheorem{theo}{Theorem}
\newtheorem{corol}{Corollary}

\newenvironment{proof}{\medskip\par\noindent{\bf Proof}}{\hfill $\Box$
\medskip\par}

\newcommand{\C}{\mathbb{C}}
\newcommand{\N}{\mathbb{N}}
\newcommand{\R}{\mathbb{R}}

\begin{document}
\title{Parametric formal Gevrey asymptotic expansions in two complex time variable problems}

\author[1]{Guoting Chen}
\author[2]{Alberto Lastra}
\author[3]{St\'ephane Malek}
\affil[1]{Great Bay University, Dongguan, Guangdong, China. {\tt guoting.chen@gbu.edu.cn}}
\affil[2]{Universidad de Alcal\'a, Dpto. F\'isica y Matem\'aticas, Alcal\'a de Henares, Madrid, Spain. {\tt alberto.lastra@uah.es}}
\affil[3]{University of Lille, Laboratoire Paul Painlev\'e, Villeneuve d'Ascq cedex, France. {\tt stephane.malek@univ-lille.fr}}


\date{}

\maketitle
\thispagestyle{empty}
{ \small \begin{center}
{\bf Abstract}
\end{center}

The analytic and formal solutions to a family of singularly perturbed partial differential equations in the complex domain involving two complex time variables are considered. The analytic continuation properties of the solution of an auxiliary problem in the Borel plane overcomes the absence of adequate domains which would guarantee summability of the formal solution. 

Moreover, several exponential decay rates of the difference of analytic solutions with respect to the perturbation parameter at the origin are observed, leading to several asymptotic levels relating the analytic and the formal solution.

\smallskip


}
\bigskip \bigskip

\section{Introduction}

The main aim of the present work is to describe the asymptotic relation existing between the analytic and the formal solutions of a family of singularly perturbed nonlinear partial differential equations in two complex time variables of the form
\begin{multline}
Q(\partial_z)u(t_1,t_2,z,\epsilon)\\
= \epsilon^{\Delta_0}(t_1^{k_1+1}\partial_{t_1})^{\delta_1}(t_2^{k_2+1}\partial_{t_2})^{\delta_2}R(\partial_z) u(t_1,t_2,z,\epsilon)+ P(t_1,t_2,\partial_{t_1},\partial_{t_2},z,\partial_z,\epsilon)u(t_1,t_2,z,\epsilon)\label{epralintro}\\ 
+\left( P_1(\epsilon,\partial_z)u(t_1,t_2,z,\epsilon)\right)\left( P_2(\epsilon,\partial_z)u(t_1,t_2,z,\epsilon)\right)+f(t_1,t_2,z,\epsilon),
\end{multline}
under initial data $u(t_1,0,z,\epsilon)\equiv u(0,t_2,z,\epsilon)\equiv 0$. In the previous equation $\epsilon$ acts as a small complex perturbation parameter. In addition to this, $Q(X),R(X)\in\C[X]$ and $P_1,P_2$ are polynomials in their second variable, with coefficients in the space of holomorphic functions on some neighborhood of the origin. $\Delta_0,k_1,k_2,\delta_1,\delta_2$ are nonnegative integers. The function $P(T_1,T_2,S_1,S_2,Z,S_3,\epsilon)$ turns out to be a polynomial in $T_1,T_2,S_1,S_2,S_3$ with coefficients being holomorphic and bounded functions on a horizontal strip, say $H$, w.r.t. the space variable $Z$, and some neighborhood of the origin in the perturbation parameter, say $D$. The forcing term $f$ is a polynomial in its two first variables, with coefficients being holomorphic and bounded functions on $H\times D$. The precise shape and assumptions on the elements involved in the equation is described in detail in Section~\ref{secmainproblem}.

\vspace{0.3cm}

This work puts a step forward in the theory of analytic and asymptotic solutions to singularly perturbed partial differential equations in the complex domain, in several complex time variables. Different advances in this direction have been achieved in the last years. In~\cite{family1}, the two last authors observed a multilevel Gevrey asymptotic phenomenon with respect to the perturbation parameter relating the formal and the analytic solutions to a family of the form (\ref{epralintro}) where the leading operator 
\begin{equation}\label{e00}
Q(\partial_z)-\epsilon^{\Delta_0}(t_1^{k_1+1}\partial_{t_1})^{\delta_1}(t_2^{k_2+1}\partial_{t_2})^{\delta_2}R(\partial_z)
\end{equation}
is replaced by a product of operators in the form 
$$Q_j(\partial_z)-\epsilon^{\Delta_j}(t_j^{k_j+1}\partial_{t_j})^{\delta_j}R_j(\partial_z),$$
$j=1,2$. This leads to symmetric asymptotic properties of the solutions. Afterwards, the  situation in which the role of the Borel time variables is asymmetric was initially considered in~\cite{family2}, observing a small division phenomena. In the present work, the configuration of the main problem does not allow us to follow the previous procedures developed in~\cite{family1,family2} due to the shape of the operator (\ref{e00}).

It is worth mentioning the phenomena observed in the previous works from a geometric point of view. Indeed, the analytic solution of the main problem is constructed as a Laplace-like transform of a function defined in the Borel space. In the work~\cite{family1}, this function can be defined with respect to the Borel time variables on sets of the form $(D_1\cup S_1)\times (D_2\cup S_2)$. Here, $S_j$ denotes an infinite sector with vertex at the origin, and $D_j$ stands for a disc centered at the origin, $j=1,2$. In~\cite{family2}, the asymmetric behavior of the time variables causes that the function in the Borel space is only defined on sets of the form $(D_1\cup S_1)\times S_2$, but not on sets of the form $S_1\times (D_2\cup S_2)$. In the present work, the auxiliary function in the Borel domain can not be defined on sets of the previous form: $(D_1\cup S_1)\times (D_2\cup S_2)$, nor $S_1\times (D_2\cup S_2)$, nor $(D_1\cup S_1)\times S_2$, but only on sets of the form $S_1\times S_2$. As a result, the deformation of the integration path defining the solution as a double Laplace-like operator is not available. This is essential in order to estimate the difference of two solutions which share a common domain in the perturbation parameter. Therefore, a strategy of a direct application of a Ramis-Sibuya theorem is no longer valid in the present situation. Alternatively, the analytic continuation of the auxiliary functions in the Borel space and an adequate splitting of the integration paths involved in the definition of the solutions will play the key point to achieve an asymptotic meaning of the solution. In the end, a fine structure related to the analytic solution is observed, involving two Gevrey orders, $k_1$ and $k_2$ which remain linked to the singular operators involved in the leading term of equation (\ref{epralintro}). This fine structure appears in the form of a decomposition of the formal and the analytic solution as the sum of terms involving different Gevrey asymptotic orders, appearing in  Balser's decomposition approach to multisummability (see Section 7.5,~\cite{loday}). 

The appearance of a scheme involving several levels relating the analytic and the formal solutions to ordinary differential equations in the complex domain (under the action of a small perturbation complex parameter or not) has been a field of interest in the scientific community during the last three decades. See the classical references~\cite{ba97,babrrasi,br92,lo94,malra92,rasi94,takei}. Regarding partial differential equations in the complex domain, the advances on multilevel solutions are much more limited. We refer to the works~\cite{ba4,mi10,michalik12,ou,ta,taya} and the references therein, among many others. Multilevel results have also been recently extended to other more general functional equations in the last years, such as~\cite{jikalasa,lamisu4,michalik12}, in the context of moment partial differential equations.

The study of solutions to singularly perturbed partial differential equations involving several complex time variables has also been considered regarding boundary layer expansions, distinguishing outer and inner solutions, together with their asymptotic representation. This is the case of~\cite{lamaboundarylayer} in which a failure of the application of a Borel-Laplace procedure is also observed. 

The main motivation for the present study consists on considering the nonlinear situation related to two previous research, namely~\cite{lamaboundarylayer2} and~\cite{chenlastramalek}. In~\cite{lamaboundarylayer2}, the two last authors consider a linearized problem with respect to (\ref{epralintro}). In both works, the leading term coincides with that in (\ref{epralintro}), that we are about to study. However, in our new setting there is more freedom in the choice of the other terms in the linear part which not only allow powers of irregular operators in the time variables of the same nature as those in the distinguished term, i.e. of the form $(t_1^{k_1+1}\partial_{t_1})^{\delta_{\ell_1}}(t_2^{k_2+1}\partial_{t_2})^{\delta_{\ell_2}}$, but a wider variety of irregular operators $t_1^{\ell_1}\partial_{t_1}^{\ell_2}t_2^{\ell_3}\partial_{t_2}^{\ell_4}$, regarding the hypothesis (\ref{e167}) of the present work.

In that previous work, a small divisor phenomenon occurs so that the Borel-Laplace classical procedure in two time variables does not apply, and one has to search for solutions in the form
$$\frac{1}{(2\pi)^{1/2}}\int_{-\infty}^{\infty}\int_{L_d}\omega(u,m,\epsilon)\exp\left(-\left(\frac{u}{\epsilon t_1}\right)^{k_1}-\left(\frac{u}{\epsilon t_2}\right)^{k_2}\right)\frac{du}{u}e^{izm}dm,$$
where $L_d$ is an infinite ray with direction $d\in\R$, for some function $\omega$. This approach leads to the construction of analytic solutions which adopt inner and outer asymptotic solutions as asymptotic representations.

On the other hand, in the work~\cite{chenlastramalek} one searches for analytic solutions of a linearized version of the equation under study in the form of a Fourier, truncated Laplace and Laplace transform of certain function
$$\frac{1}{(2\pi)^{1/2}}\int_{-\infty}^{\infty}\int_{L_{d_1,\epsilon}}\int_{L_{d_2}}\omega(u_1,u_2,m,\epsilon)\exp\left(-\left(\frac{u_1}{\epsilon t_1}\right)^{k_1}-\left(\frac{u_2}{\epsilon t_2}\right)^{k_2}\right)\frac{du_2}{u_2}\frac{du_1}{u_1}e^{izm}dm,$$
for some $d_1,d_2\in\R$ and where $L_{d_2}$ is an infinite ray with direction $d_2$, with $L_{d_1,\epsilon}$ is the segment $[0,h_1(\epsilon)e^{id_1}]$, for some holomorphic function $\epsilon\mapsto h_1(\epsilon)$. In that previous work, all the terms except from the leading term are of the form $(t_1^{k_1+1}\partial_{t_1})^{\delta_{\ell_1}}t_2^{d_{\ell_2}}\partial_{t_2}^{\delta_{\ell_2}}$, in contrast to the freedom acquired in the present work. One can additionally observe that the freedom acquired on the possible values of the parameters involved in the second time variable coincide in both works (see the second condition in (\ref{e111}).

The strategy of the present work is as follows. First, we search for solutions to the main problem in the form of a double Laplace and inverse Fourier transform of some auxiliary function. This allows us to substitute the main equation (\ref{epralintro}) by an auxiliary convolution problem in the Borel domain (see (\ref{e309})). The solution to this auxiliary equation, say $(\tau_1,\tau_2,m,\epsilon)\mapsto \omega(\tau_1,\tau_2,m,\epsilon)$ is in principle defined on $S_1\times S_2\times \R\times D$, where  $S_1\times S_2$ stands for a product of unbounded sectors, and $D$ is a punctured disc at the origin. However, it is proved (see Proposition~\ref{propauxpral}) that such function can be analytically extended as follows:
\begin{itemize}
\item[-] for all $\tau_1\in S_1$ near the origin, $m\in\R$ and $\epsilon\in D$, the map $\tau_2\mapsto \omega(\tau_1,\tau_2,m,\epsilon)$, defined on $S_2$, can be analytically continued to some neighborhood of the origin.
\item[-] for all $\tau_2\in S_2$ near the origin, $m\in\R$ and $\epsilon\in D$, the map $\tau_1\mapsto \omega(\tau_1,\tau_2,m,\epsilon)$, defined on $S_1$, can be analytically continued to some neighborhood of the origin.
\end{itemize}

The construction of the analytic solution $u(t_1,t_2,z,\epsilon)$ is achieved by means of a double Laplace transform with respect to $(\tau_1,\tau_2)$ (see Theorem~\ref{teo1}) following a classical procedure. It turns out to be holomorphic and bounded on $\mathcal{T}_1\times\mathcal{T}_2\times H_{\beta'}\times \mathcal{E}$, with $\mathcal{T}_j,\mathcal{E}$ being bounded sectors with vertex at the origin, and $H_{\beta'}$ a horizontal strip. Although the mentioned construction of the analytic solution $u(t_1,t_2,z,\epsilon)$ needs no analytic continuation of the auxiliary function $\omega(\tau_1,\tau_2,m,\epsilon)$ to be defined, such analytic continuation is indeed needed for setting the existence of an asymptotic expansion.

In order to attain such asymptotic expansion, we split the solution as the sum of three terms, $J_1$ (see Section~\ref{secj1}), $J_2$ and $J_3$ (see Section~\ref{secj23}). The term $J_2$ (resp. $J_3$) has null Gevrey asymptotic expansion of order $1/k_2$ (resp. $1/k_1$), see Proposition~\ref{prop8} and Proposition~\ref{prop9}. On the other hand, a parametric Gevrey series expansion associated to $J_1$ is attained by completing $J_1$ into a set $(J_{1,p})_{0\le p\le \varsigma-1}$, each of them holomorphic on a finite sector $\mathcal{E}_p$, where $(\mathcal{E}_p)_{0\le p\le \varsigma}$ describes a good covering (see Definition~\ref{defi662}) and by describing the exponential decay of the difference of two consecutive maps (i.e. with nonempty intersection of their domains of definition) in the perturbation parameter (see Proposition~\ref{prop689}). 

The application of a multilevel Ramis-Sibuya Theorem (RS) to $J_1$ allows to conclude the main result of the present work (Theorem~\ref{teo2}). More precisely, this result states the existence of a decomposition of the analytic solution of the main problem in the form
$$u(t_1,t_2,z,\epsilon)=a(t_1,t_2,z,\epsilon)+u_{1}(t_1,t_2,z,\epsilon)+u_{2}(t_1,t_2,z,\epsilon),$$
where $a(t_1,t_2,z,\epsilon)$ turns out to be a holomorphic function on some neighborhood of the origin with respect to the perturbation parameter $\epsilon$, and with coefficients belonging to the Banach space of holomorphic and bounded functions on $\mathcal{T}_1\times \mathcal{T}_2\times H_{\beta'}$. Moreover, $u_{j}$ are holomorphic and bounded functions on $\mathcal{T}_1\times\mathcal{T}_2\times H_{\beta'}\times \mathcal{E}$, for $j=1,2$. On the other hand, there exist formal power series in $\epsilon$ with coefficients in the Banach space of holomorphic and bounded functions on $\mathcal{T}_1\times\mathcal{T}_2\times H_{\beta'}$, say $\hat{u}_1(t_1,t_2,z,\epsilon)$ and $\hat{u}_2(t_1,t_2,z,\epsilon)$, which satisfy that $u_{j}$ admits $\hat{u}_{j}$ as its common Gevrey asymptotic expansion of order $1/k_j$ with respect to $\epsilon$ on $\mathcal{E}$, $j=1,2$.

\vspace{0.3cm}

\textbf{Notation:}

We write $\N$ for the set of positive integers, and $\N_0=\N\cup\{0\}$.

For all $z_0\in\C$ and $r>0$, we denote the open disc centered at $z_0$ and radius $r$ by $D(z_0,r)$. 

Given two open and bounded sectors in the complex domain $\mathcal{E}$ and $\mathcal{T}$, with vertex at the origin, we say that $T$ is a proper subsector of $\mathcal{E}$, and denote it by $T\prec \mathcal{E}$ whenever $\overline{T}\setminus\{0\}\subseteq \mathcal{E}$. 

Given a complex Banach space $\mathbb{E}$ and a nonempty open set $U\subseteq\C$, we write $\mathcal{O}_b(U,\mathbb{E})$ for the set of holomorphic functions defined on $U$ with values in $\mathbb{E}$. We simply write $\mathcal{O}_b(U)$ in the case that $\mathbb{E}=\C$. We write $\mathbb{E}\{\epsilon\}$ for the set of holomorphic functions with values on $\mathbb{E}$, defined on $\epsilon$, which are convergent on some neighborhood of the origin.

\section{Statement of the main problem}\label{secmainproblem}

In this section, we establish the main Cauchy problem under study, providing the details on the elements involved in the problem.

Let $k_1,k_2$ be positive integer numbers. Let us assume that $k_1>k_2\ge1$ and fix $\epsilon_0>0$.  We also consider positive integers $\delta_{1},\delta_{2}$ and a finite set $I\subseteq \N^4$. 
We assume 
$$\delta_1k_1=\delta_2k_2,$$
and define
$$\Delta_{0}:=k_1\delta_{1}+k_2\delta_{2}.$$

We assume that for every $\underline{\ell}=(\ell_1,\ell_2,\ell_3,\ell_4)\in I$, one has $\ell_2$ and $\ell_4$ are positive integers such that
\begin{equation}\label{e114}
\ell_1=\ell_2(k_1+1)+d_{k_1,\ell_1,\ell_2} \hbox{ and } \ell_3=\ell_4(k_2+1)+d_{k_2,\ell_3,\ell_4},
\end{equation}
for some positive integers $d_{k_1,\ell_1,\ell_2}, d_{k_2,\ell_3,\ell_4}$. Moreover, for every $\underline{\ell}=(\ell_1,\ell_2,\ell_3,\ell_4)\in I$ we fix an integer number $\Delta_{\underline{\ell}}$ such that
\begin{equation}\label{e111}
\Delta_{\underline{\ell}}\ge \ell_1-\ell_2+\ell_3-\ell_4+1.
\end{equation}

We also assume that
\begin{equation}\label{e167}
k_1\le \ell_1-\ell_2=\ell_3-\ell_4\le \delta_1k_1,
\end{equation}
for every $(\ell_1,\ell_2,\ell_3,\ell_4)\in I$.

Let us fix polynomials $Q(X),R(X),R_{\underline{\ell}}(X)\in\C[X]$, for $\underline{\ell}\in I$. We assume that 
\begin{equation}\label{e120}
\hbox{deg}(Q)\ge \hbox{deg}(R)\ge\hbox{deg}(R_{\underline{\ell}}),\quad \underline{\ell}\in I,
\end{equation}
and
\begin{equation}\label{e121}
R(im)\neq 0,\quad Q(im)\neq 0,\quad R_{\underline{\ell}}(im)\neq 0
\end{equation}
for every $m\in\R$ and all $\underline{\ell}\in I$. In addition to this, we assume the existence of an infinite sector $S_{Q,R}$, centered at the origin such that 
\begin{equation}\label{e133}
\left\{\frac{Q(im)}{R(im)}:m\in\R\right\}\subseteq S_{Q,R}.
\end{equation} 

\vspace{0.3cm}

\textbf{Remark:} An example of a situation in which the previous assumption holds is that in which one considers polynomials $Q$ and $R$ with positive coefficients with only powers which are divisible by 4, and where $S_{Q,R}$ is a sector containing the set of positive real numbers.

\vspace{0.3cm}

We also fix polynomials $P_1(\epsilon,X),P_2(\epsilon,X)\in\mathcal{O}_b(D(0,\epsilon_0))[X]$. We assume that 
\begin{equation}\label{e136}
\hbox{deg}(R)\ge\max\{\hbox{deg}(P_1),\hbox{deg}(P_2)\}.
\end{equation}

In view of the assumptions made on $P_1$ and $P_2$, there exist $C_{P_1},C_{P_2}>0$ such that
\begin{equation}\label{e706}
|P_{1}(\epsilon,im)|\le C_{P_1}(1+|m|)^{\hbox{deg}(P_1)},\quad |P_{2}(\epsilon,im)|\le C_{P_2}(1+|m|)^{\hbox{deg}(P_2)},
\end{equation}
for every $m\in\R$ and $\epsilon\in D(0,\epsilon_0)$. We also assume that $C_{P_1}$ and $C_{P_2}$ are small enough (see Proposition~\ref{prop2}).

We choose $\beta>0$ and fix $\mu>0$ such that
\begin{equation}\label{e137}
\mu>\max\{\hbox{deg}(P_1),\hbox{deg}(P_2),\max_{\underline{\ell}\in I}\{\hbox{deg}(R_{\underline{\ell}})\}\}+1.
\end{equation}
We consider the following nonlinear initial value problem
\begin{multline}\label{epral}
Q(\partial_z)u(t_1,t_2,z,\epsilon)=\epsilon^{\Delta_{0}}(t_1^{k_1+1}\partial_{t_1})^{\delta_{1}}(t_2^{k_2+1}\partial_{t_2})^{\delta_{2}}R(\partial_z)u(t_1,t_2,z,\epsilon)\\
+\sum_{\underline{\ell}=(\ell_1,\ell_2,\ell_3,\ell_4)\in I}\epsilon^{\Delta_{\underline{\ell}}}c_{\underline{\ell}}(z,\epsilon)t_1^{\ell_1}\partial_{t_1}^{\ell_2}t_2^{\ell_3}\partial_{t_2}^{\ell_4}R_{\underline{\ell}}(\partial_z)u(t_1,t_2,z,\epsilon)\\
+\left( P_1(\epsilon,\partial_z)u(t_1,t_2,z,\epsilon)\right)\left( P_2(\epsilon,\partial_z)u(t_1,t_2,z,\epsilon)\right)+f(t_1,t_2,z,\epsilon),
\end{multline}
under vanishing initial data $u(0,t_2,z,\epsilon)=u(t_1,0,z,\epsilon)=0$.

\vspace{0.3cm}

\textbf{Remark:} Observe all the previous assumptions made on the parameters involved in the problem are compatible. An example of equation satisfying the previous assumptions is the following:
\begin{multline*}
(\partial_z^4+1)u(t_1,t_2,z,\epsilon)=\epsilon^{12}(t_1^{4}\partial_{t_1})^{2}(t_2^{3}\partial_{t_2})^3+\epsilon^{11}c_{(6,1,7,2)}(z,\epsilon)t_1^6\partial_{t_1}t_2^7\partial_{t_2}^2(-\partial_z^2+2)u(t_1,t_2,z,\epsilon)\\
+\left( P_1(\epsilon,\partial_z)u(t_1,t_2,z,\epsilon)\right)\left( P_2(\epsilon,\partial_z)u(t_1,t_2,z,\epsilon)\right)+f(t_1,t_2,z,\epsilon),
\end{multline*}
with $\max\{\hbox{deg}(P_1),\hbox{deg}(P_2)\}\le 2$. The functions $f$ and $c_{(6,1,7,2)}$ are attained to the construction described below.

\vspace{0.3cm}

Let $0<\beta'<\beta$ be fixed. We denote $H_{\beta'}$ the horizontal strip
$$H_{\beta'}=\left\{z\in\C:|\hbox{Im}(z)|<\beta'\right\}.$$
For every $\underline{\ell}\in I$, the function $c_{\underline{\ell}}$ belongs to $\mathcal{O}_b(H_{\beta'}\times D(0,\epsilon_0))$. The function $f$ turns out to be a polynomial in its two first variables, and a holomorphic and bounded function on $H_{\beta'}\times D(0,\epsilon_0)$ with respect to $(z,\epsilon)$. These functions are constructed as follows. 

Given $\underline{\ell}\in I$, we choose a function $\R\times D(0,\epsilon_0)\ni(m,\epsilon)\mapsto C_{\underline{\ell}}(m,\epsilon)$ under the following assumptions:
\begin{itemize}
\item For every $\epsilon\in D(0,\epsilon_0)$, the function $\R\ni m\mapsto C_{\underline{\ell}}(m,\epsilon)$ is continuous on $\R$, and it satisfies there exists $K_{\underline{\ell}}(\epsilon)>0$ with
$$|C_{\underline{\ell}}(m,\epsilon)|\le K_{\underline{\ell}}(\epsilon)\frac{1}{(1+|m|)^{\mu}}e^{-\beta |m|},\qquad m\in\R.$$
Moreover, there exists $K>0$ such that 
$$\sup_{\underline{\ell}\in I}\sup_{\epsilon\in D(0,\epsilon_0)}K_{\underline{\ell}}(\epsilon)\le K.$$
\item For all $m\in\R$, the mapping $D(0,\epsilon_0)\ni\epsilon\mapsto C_{\underline{\ell}}(m,\epsilon)$ is a holomorphic function on $D(0,\epsilon_0)$. 
\end{itemize}

In view of the previous assumptions, and the properties of inverse Fourier transform (see Appendix~\ref{secapendice1}), one defines
$$c_{\underline{\ell}}(z,\epsilon):=\mathcal{F}^{-1}\left(m\mapsto C_{\underline{\ell}}(m,\epsilon)\right)(z)$$
which represents a holomorphic and bounded function on $H_{\beta'}\times D(0,\epsilon_0)$.

\vspace{0.3cm}

\noindent\textbf{Remark:} Observe from Definition~\ref{defi937} in Section~\ref{secapendice}, that for every  $\epsilon\in D(0,\epsilon_0)$, the function $m\mapsto C_{\underline{\ell}}(m,\epsilon)$ belongs to the Banach space $E_{(\beta,\mu)}$, with 
\begin{equation}\label{e166}
\sup_{\underline{\ell}\in I}\sup_{\epsilon\in D(0,\epsilon_0)}\left\|C_{\underline{\ell}}(m,\epsilon)\right\|_{(\beta,\mu)}\le K.
\end{equation}

\vspace{0.3cm}

For the construction of the forcing term, we make use of the following well-known property of Laplace transform.

\begin{lemma}\label{lema1}
Let $k,n$ be positive integers. We also fix $d\in\R$. For very $T\in\C^{\star}$ with $\cos((d-\hbox{arg}(T))k)>0$, it holds that
$$k\int_{L_d}u^{n-1}\exp\left(-\left(\frac{u}{T}\right)^{k}\right)du=T^n\Gamma\left(\frac{n}{k}\right),$$
with $L_d$ being the integration path $[0,\infty)\ni r\mapsto re^{id}$. Here, $\Gamma(\cdot)$ stands for Gamma function.
\end{lemma}

The forcing term $f$ is constructed as follows. Let $N_1,N_2\subseteq\N$ be two nonempty finite subsets of positive integers. For every $(n_1,n_2)\in N_1\times N_2$ let $\R\times D(0,\epsilon_0)\ni (m,\epsilon)\mapsto F_{n_1,n_2}(m,\epsilon)$ be a function under the following assumptions:
\begin{itemize}
\item For every $\epsilon\in D(0,\epsilon_0)$, the function $m\mapsto F_{n_1,n_2}(m,\epsilon)$ is continuous on $\R$, and it satisfies there exists $K_{n_1,n_2}(\epsilon)>0$ such that
$$|F_{n_1,n_2}(m,\epsilon)|\le K_{n_1,n_2}(\epsilon)\frac{1}{(1+|m|)^{\mu}}e^{-\beta |m|},\qquad m\in\R.$$
\item For all $m\in\R$, the mapping $D(0,\epsilon_0)\ni\epsilon\mapsto F_{n_1,n_2}(m,\epsilon)$ is a holomorphic function. Moreover, there exists $\tilde{K}>0$  such that
$$\sup_{(n_1,n_2)\in N_1\times N_2}\sup_{\epsilon\in D(0,\epsilon_0)}K_{n_1,n_2}(\epsilon)\le \tilde{K}.$$
\end{itemize}
Let us define the function $\Psi$ on $\C^2\times\R\times D(0,\epsilon_0)$ by
$$\Psi(\tau_1,\tau_2,m,\epsilon):=\sum_{(n_1,n_2)\in N_1\times N_2} F_{n_1,n_2}(m,\epsilon)\frac{\tau_1^{n_1}}{\Gamma\left(\frac{n_1}{k_1}\right)}\frac{\tau_2^{n_2}}{\Gamma\left(\frac{n_2}{k_2}\right)},$$
and consider its inverse Fourier and double Laplace transform, giving rise to the function
\begin{multline*}
f(t_1,t_2,z,\epsilon):=\frac{k_1k_2}{(2\pi)^{1/2}}\int_{-\infty}^{\infty}\int_{L_{d_1}}\int_{L_{d_2}}\Psi(\tau_1,\tau_2,m,\epsilon)\\
\times \exp\left(-\left(\frac{\tau_1}{\epsilon t_1}\right)^{k_1}-\left(\frac{\tau_2}{\epsilon t_2}\right)^{k_2}\right)e^{izm}\frac{d\tau_2}{\tau_2}\frac{d\tau_1}{\tau_1}dm.
\end{multline*}
Here, $d_1,d_2\in\R$ are arbitrary numbers, and $L_{d_j}$ denotes the integration path $[0,\infty)e^{id_j}$, for $j=1,2$. Observe that, due to Cauchy Theorem and the fact that $\Psi$ is a polynomial in its two first variables, the directions $d_1$ and $d_2$ can be arbitrarily chosen. In addition to that, in view of Lemma~\ref{lema1}, it is straight to check that $f$ determines a bounded holomorphic function with respect to $(z,\epsilon)$ on $H_{\beta'}\times D(0,\epsilon_0)$, and a polynomial with respect to its first two variables. In addition to this, 
\begin{align*}
f(t_1,t_2,z,\epsilon)&=\sum_{(n_1,n_2)\in N_1\times N_2}\left[\frac{1}{(2\pi)^{1/2}}\int_{-\infty}^{\infty}F_{n_1,n_2}(m,\epsilon)e^{izm}dm\right](\epsilon t_1)^{n_1}(\epsilon t_2)^{n_2}\\
&=:\sum_{(n_1,n_2)\in N_1\times N_2}\mathfrak{F}_{n_1,n_2}(z,\epsilon)(\epsilon t_1)^{n_1}(\epsilon t_2)^{n_2},
\end{align*}
for $(t_1,t_2,z,\epsilon)\in \C^2\times H_{\beta'}\times D(0,\epsilon_0)$.

\vspace{0.3cm}

\noindent\textbf{Remark:}\label{rem1} Observe from Definition~\ref{defi937} in Section~\ref{secapendice} that the function $m\mapsto F_{n_1,n_2}(m,\epsilon)$ belongs to $E_{(\beta,\mu)}$ for every $\epsilon\in D(0,\epsilon_0)$, with 
\begin{equation}\label{e280}
\sup_{\epsilon\in D(0,\epsilon_0), (n_1,n_2)\in N_1\times N_2}\left\|F_{n_1,n_2}(m,\epsilon)\right\|_{(\beta,\mu)}\le \tilde{K}.
\end{equation}
It is straightforward to check that for every $\epsilon\in D(0,\epsilon_0)$ and $\rho_1,\rho_2>0$, one has that $\Psi\in B_{(\beta,\mu,\rho_1,\rho_2)}$ (see Definition~\ref{def963} from Subsection~\ref{sec452}). Indeed,
\begin{equation}\label{e215}
\left\|\Psi(\tau_1,\tau_2,m,\epsilon)\right\|_{(\beta,\mu,\rho_1,\rho_2)}\le \tilde{K}\sum_{(n_1,n_2)\in N_1\times N_2}\frac{\rho_1^{n_1-1}}{\Gamma\left(\frac{n_1}{k_1}\right)}\frac{\rho_2^{n_2-1}}{\Gamma\left(\frac{n_2}{k_2}\right)}=:C_{\Psi}<\infty.
\end{equation}
Observe this quantity can be diminish as desired as long as $\tilde{K}$ is small enough.

\section{Analytic solutions to the main problem}\label{secanal}

This section is devoted to the construction of analytic solutions to the main problem under study (\ref{epral}). First, we sketch the strategy to follow, and later we describe the solution by means of auxiliary problems solved by means of fixed point arguments involving certain operators defined on appropriate Banach spaces. 

\subsection{Strategy for the construction of the analytic solutions to (\ref{epral})}

We consider that the assumptions and constructions associated to the main problem (\ref{epral}), made in Section~\ref{secmainproblem} hold, and we search for solutions of (\ref{epral}) in the form of an inverse Fourier and double Laplace transform, for every fixed value of the perturbation parameter $\epsilon$. More precisely, we will search for solutions of (\ref{epral}) in the form
\begin{multline}\label{eanal}
u(t_1,t_2,z,\epsilon)=\frac{k_1k_2}{(2\pi)^{1/2}}\int_{-\infty}^{\infty}\int_{L_{d_1}}\int_{L_{d_2}}\omega(u_1,u_2,m,\epsilon)\\
\times \exp\left(-\left(\frac{u_1}{\epsilon t_1}\right)^{k_1}-\left(\frac{u_2}{\epsilon t_2}\right)^{k_2}\right)e^{izm}\frac{du_2}{u_2}\frac{du_1}{u_1}dm
\end{multline}
for some $d_1,d_2\in\R$ to be determined, and where $L_{d_j}$ stands for the integration path $[0,\infty)e^{id_j}$, for $j=1,2$.

We proceed by following several steps in the search of a solucion of (\ref{epral}) in the form (\ref{eanal}). First, let us search for a function $u(t_1,t_2,z,\epsilon)$ showing a monomial behavior with respect to its first two variables, i.e. assume that  
$$u(t_1,t_2,z,\epsilon)=U(\epsilon t_1,\epsilon t_2,z,\epsilon),$$
for some function $U$ defined on appropriate domains, to be specified. A direct inspection of the previous elements yields the next result.

\begin{lemma}
From the formal point of view, $u(t_1,t_2,z,\epsilon)$ solves (\ref{epral}) whenever $U(T_1,T_2,z,\epsilon)$ provides a formal solution of
\begin{multline}\label{e214}
Q(\partial_z)U(T_1,T_2,z,\epsilon)=(T_1^{k_1+1}\partial_{T_1})^{\delta_1}(T_2^{k_2+1}\partial_{T_2})^{\delta_2}R(\partial_z)U(T_1,T_2,z,\epsilon)\\
+\sum_{\underline{\ell}=(\ell_1,\ell_2,\ell_3,\ell_4)\in I}\epsilon^{\Delta_{\underline{\ell}}-\ell_1+\ell_2-\ell_3+\ell_4}c_{\underline{\ell}}(z,\epsilon)T_1^{\ell_1}\partial_{T_1}^{\ell_2}T_2^{\ell_3}\partial_{T_2}^{\ell_4}R_{\underline{\ell}}(\partial_z)U(T_1,T_2,z,\epsilon)\\
+\left(P_1(\epsilon,\partial_z)U(T_1,T_2,z,\epsilon)\right)\left( P_2(\epsilon,\partial_z)U(T_1,T_2,z,\epsilon)\right)+F(T_1,T_2,z,\epsilon),
\end{multline}
where 
\begin{multline}\label{eF}
F(T_1,T_2,z,\epsilon)=\frac{k_1k_2}{(2\pi)^{1/2}}\int_{-\infty}^{\infty}\int_{L_{d_1}}\int_{L_{d_2}}\Psi(\tau_1,\tau_2,m,\epsilon)\\
\times \exp\left(-\left(\frac{\tau_1}{T_1}\right)^{k_1}-\left(\frac{\tau_2}{T_2}\right)^{k_2}\right)e^{izm}\frac{d\tau_2}{\tau_2}\frac{d\tau_1}{\tau_1}dm.
\end{multline}
\end{lemma}

As a second step, we rewrite the terms involved in the equation (\ref{e214}) by means of the next result.

\begin{lemma}[Formula (8.7),~\cite{taya}]\label{lema3}
Let $m,k$ be positive integers. Then, for all $1\le \ell\le m-1$, there exists a constant $A_{m,\ell}\in\R$ such that
$$T^{m(k+1)}\partial_T^{m}=(T^{k+1}\partial_T)^m+\sum_{\ell=1}^{m-1}A_{m,\ell}T^{k(m-\ell)}(T^{k+1}\partial_T)^{\ell}.$$
\end{lemma}

Taking into account the previous Lemma, together with the assumption (\ref{e114}) made on the elements involved in the equation, one can write (\ref{e214}) in the form
\begin{multline}\label{e234}
Q(\partial_z)U(T_1,T_2,z,\epsilon)=(T_1^{k_1+1}\partial_{T_1})^{\delta_1}(T_2^{k_2+1}\partial_{T_2})^{\delta_2}R(\partial_z)U(T_1,T_2,z,\epsilon)\\
+\sum_{\underline{\ell}=(\ell_1,\ell_2,\ell_3,\ell_4)\in I}\epsilon^{\Delta_{\underline{\ell}}-\ell_1+\ell_2-\ell_3+\ell_4}c_{\underline{\ell}}(z,\epsilon)T_1^{d_{k_1,\ell_1,\ell_2}}\left[(T_1^{k_1+1}\partial_{T_1})^{\ell_2}+\sum_{h=1}^{\ell_2-1}A_{\ell_2,h}T_1^{k_1(\ell_2-h)}(T_1^{k_1+1}\partial_{T_1})^{h}\right]\\
\times T_2^{d_{k_2,\ell_3,\ell_4}}\left[(T_2^{k_2+1}\partial_{T_2})^{\ell_4}+\sum_{j=1}^{\ell_4-1}A_{\ell_4,j}T_2^{k_2(\ell_4-j)}(T_2^{k_2+1}\partial_{T_2})^{j}\right]R_{\underline{\ell}}(\partial_z)U(T_1,T_2,z,\epsilon)\\+\left( P_1(\epsilon,\partial_z)U(T_1,T_2,z,\epsilon)\right)\left( P_2(\epsilon,\partial_z)U(T_1,T_2,z,\epsilon)\right)+F(T_1,T_2,z,\epsilon),
\end{multline}
for $A_{\ell_2,h},A_{\ell_4,j}$ for $1\le h\le \ell_2-1$ and $1\le j\le \ell_4-1$ determined in Lemma~\ref{lema3}.

Our strategy of finding solutions of (\ref{epral}) in the form (\ref{eanal}) lead us to search for solutions of (\ref{e234}) in the form 
\begin{multline}\label{eanal2}
U(T_1,T_2,z,\epsilon)=\frac{k_1k_2}{(2\pi)^{1/2}}\int_{-\infty}^{\infty}\int_{L_{d_1}}\int_{L_{d_2}}\omega(u_1,u_2,m,\epsilon)\\
\times \exp\left(-\left(\frac{u_1}{T_1}\right)^{k_1}-\left(\frac{u_2}{T_2}\right)^{k_2}\right)e^{izm}\frac{du_2}{u_2}\frac{du_1}{u_1}dm,
\end{multline}
defined on certain domains to be clarified.

The following result describes the action of some operators on an expression of the form (\ref{eanal2}).

\begin{lemma}\label{lema2}
Assume the expression of $U(T_1,T_2,z,\epsilon)$ is formally defined by (\ref{eanal2}). The following statements hold:
\begin{itemize}
\item[(1)] For $j=1,2$, 
\begin{multline*}
T_j^{k_j+1}\partial_{T_j}U(T_1,T_2,z,\epsilon)=\frac{k_1k_2}{(2\pi)^{1/2}}\int_{-\infty}^{\infty}\int_{L_{d_1}}\int_{L_{d_2}}k_ju_j^{k_j}\omega(u_1,u_2,m,\epsilon)\\
\times \exp\left(-\left(\frac{u_1}{T_1}\right)^{k_1}-\left(\frac{u_2}{T_2}\right)^{k_2}\right)e^{izm}\frac{du_2}{u_2}\frac{du_1}{u_1}dm.
\end{multline*}
\item[(2)] For every positive integer $m_1$, the expression $T_1^{m_1}U(T_1,T_2,z,\epsilon)$ equals
\begin{multline*}
\frac{k_1k_2}{(2\pi)^{1/2}}\int_{-\infty}^{\infty}\int_{L_{d_1}}\int_{L_{d_2}}\left[\frac{u_1^{k_1}}{\Gamma\left(\frac{m_1}{k_1}\right)}\int_0^{u_1^{k_1}}(u_1^{k_1}-s_1)^{\frac{m_1}{k_1}-1}\omega(s_1^{\frac{1}{k_1}},u_2,m,\epsilon)\frac{ds_1}{s_1}\right]\\
\times \exp\left(-\left(\frac{u_1}{T_1}\right)^{k_1}-\left(\frac{u_2}{T_2}\right)^{k_2}\right)e^{izm}\frac{du_2}{u_2}\frac{du_1}{u_1}dm.
\end{multline*}
\item[(3)] For every positive integer $m_2$, the expression $T_2^{m_2}U(T_1,T_2,z,\epsilon)$ equals
\begin{multline*}
\frac{k_1k_2}{(2\pi)^{1/2}}\int_{-\infty}^{\infty}\int_{L_{d_1}}\int_{L_{d_2}}\left[\frac{u_2^{k_2}}{\Gamma\left(\frac{m_2}{k_2}\right)}\int_0^{u_2^{k_2}}(u_2^{k_2}-s_2)^{\frac{m_2}{k_2}-1}\omega(u_1,s_2^{\frac{1}{k_2}},m,\epsilon)\frac{ds_2}{s_2}\right]\\
\times \exp\left(-\left(\frac{u_1}{T_1}\right)^{k_1}-\left(\frac{u_2}{T_2}\right)^{k_2}\right)e^{izm}\frac{du_2}{u_2}\frac{du_1}{u_1}dm.
\end{multline*}
\item[(4)] 
\begin{multline*}
\left(P_1(\epsilon,\partial_z)U(T_1,T_2,z,\epsilon)\right)\left(P_2(\epsilon,\partial_z)U(T_1,T_2,z,\epsilon)\right)\\
=\frac{k_1k_2}{(2\pi)^{1/2}}\int_{-\infty}^{\infty}\int_{L_{d_1}}\int_{L_{d_2}}\left[u_1^{k_1}u_2^{k_2}\int_0^{u_1^{k_1}}\int_0^{u_2^{k_2}}\frac{1}{(2\pi)^{1/2}}\int_{-\infty}^{\infty}P_1(\epsilon,i(m-m_1))\right.\\
\times\omega\left((u_1^{k_1}-s_1)^{1/k_1},(u_2^{k_2}-s_2)^{1/k_2},m-m_1,\epsilon\right)P_2(\epsilon,im_1)\omega\left(s_1^{1/k_1},s_2^{1/k_2},m_1,\epsilon\right)\\
\hfill\left.\times \frac{1}{(u_1^{k_1}-s_1)}\frac{1}{s_1}\frac{1}{(u_2^{k_2}-s_2)}\frac{1}{s_2}dm_1ds_2ds_1\right]\exp\left(-\left(\frac{u_1}{T_1}\right)^{k_1}-\left(\frac{u_2}{T_2}\right)^{k_2}\right)e^{izm}\frac{du_2}{u_2}\frac{du_1}{u_1}dm.
\end{multline*}
\end{itemize}
\end{lemma}
\begin{proof}
The expressions (1) is a direct consequence of the derivation under the integral symbol. (2) and (3) are derived from Fubini theorem. The formula (4) is derived from Fubini's theorem and an application of the identities appearing in Proposition~\ref{prop1anexo} from Subsection~\ref{secapendice1}. We refer to Lemma 1 in~\cite{family2} for a comprehensive proof. 
\end{proof}

The actions of the operators described in Lemma~\ref{lema2} and the expression of the auxiliary equation (\ref{e234}) allow us to reduce the problem of finding $U(T_1,T_2,z,\epsilon)$ to that of searching a function $\omega(\tau_1,\tau_2,m,\epsilon)$, defined in appropriate domains, which solves the following convolution equation in the Borel-Fourier plane:
\begin{multline}\label{e309}
Q(im)\omega(\tau_1,\tau_2,m,\epsilon)=(k_1\tau_1^{k_1})^{\delta_1}(k_2\tau_2^{k_2})^{\delta_2}R(im)\omega(\tau_1,\tau_2,m,\epsilon)\\
=\sum_{\underline{\ell}=(\ell_1,\ell_2,\ell_3,\ell_4)\in I}\epsilon^{\Delta_{\underline{\ell}}-\ell_1+\ell_2-\ell_3+\ell_4}\frac{1}{(2\pi)^{1/2}}\int_{-\infty}^{\infty}C_{\underline{\ell}}(m-m_1,\epsilon)\left[\mathcal{A}_1(\tau_1,\tau_2,m_1,\epsilon)+\mathcal{A}_2(\tau_1,\tau_2,m_1,\epsilon)\right.\\
\left.+\mathcal{A}_3(\tau_1,\tau_2,m_1,\epsilon)+\mathcal{A}_4(\tau_1,\tau_2,m_1,\epsilon)\right]dm_1+\tau_1^{k_1}\tau_2^{k_2}\int_0^{\tau_1^{k_1}}\int_0^{\tau_2^{k_2}}\frac{1}{(2\pi)^{1/2}}\int_{-\infty}^{\infty}P_1(\epsilon,i(m-m_1))\\
\times\omega\left((\tau_1^{k_1}-s_1)^{1/k_1},(\tau_2^{k_2}-s_2)^{1/k_2},m-m_1,\epsilon\right)P_2(\epsilon,im_1)\omega(s_1^{1/k_1},s_2^{1/k_2},m_1,\epsilon)\\
\times\frac{1}{(\tau_1^{k_1}-s_1)}\frac{1}{s_1}\frac{1}{(\tau_2^{k_2}-s_2)}\frac{1}{s_2}ds_2ds_1dm_1+\Psi(\tau_1,\tau_2,m,\epsilon),
\end{multline}
where the expressions $\mathcal{A}_j$ for $1\le j\le 4$ are determined by
\begin{multline}\label{e317}
\mathcal{A}_1(\tau_1,\tau_2,m_1,\epsilon)=\frac{\tau_1^{k_1}}{\Gamma\left(\frac{d_{k_1,\ell_1,\ell_2}}{k_1}\right)}\frac{\tau_2^{k_2}}{\Gamma\left(\frac{d_{k_2,\ell_3,\ell_4}}{k_2}\right)}\int_0^{\tau_1^{k_1}}\int_0^{\tau_2^{k_2}}(\tau_1^{k_1}-s_1)^{\frac{d_{k_1,\ell_1,\ell_2}}{k_1}-1}
(\tau_2^{k_2}-s_2)^{\frac{d_{k_2,\ell_3,\ell_4}}{k_2}-1}\\
\times(k_1s_1)^{\ell_2}(k_2s_2)^{\ell_4}R_{\underline{\ell}}(im_1)\omega(s_1^{1/k_1},s_2^{1/k_2},m_1,\epsilon)\frac{ds_2}{s_2}\frac{ds_1}{s_1},
\end{multline}
\begin{multline}\label{e318}
\mathcal{A}_2(\tau_1,\tau_2,m_1,\epsilon)=\sum_{1\le h\le \ell_2-1}A_{\ell_2,h}\frac{\tau_1^{k_1}}{\Gamma\left(\frac{k_1(\ell_2-h)+d_{k_1,\ell_1,\ell_2}}{k_1}\right)}\frac{\tau_2^{k_2}}{\Gamma\left(\frac{d_{k_2,\ell_3,\ell_4}}{k_2}\right)}\\
\times\int_0^{\tau_1^{k_1}}\int_0^{\tau_2^{k_2}}(\tau_1^{k_1}-s_1)^{\frac{k_1(\ell_2-h)+d_{k_1,\ell_1,\ell_2}}{k_1}-1}
(\tau_2^{k_2}-s_2)^{\frac{d_{k_2,\ell_3,\ell_4}}{k_2}-1}(k_1s_1)^{h}(k_2s_2)^{\ell_4}R_{\underline{\ell}}(im_1)\\
\times\omega(s_1^{1/k_1},s_2^{1/k_2},m_1,\epsilon)\frac{ds_2}{s_2}\frac{ds_1}{s_1},
\end{multline}
\begin{multline}\label{e319}
\mathcal{A}_3(\tau_1,\tau_2,m_1,\epsilon)=\sum_{1\le q\le \ell_4-1}A_{\ell_4,q}\frac{\tau_1^{k_1}}
{\Gamma\left(\frac{d_{k_1,\ell_1,\ell_2}}{k_1}\right)}\frac{\tau_2^{k_2}}{\Gamma\left(\frac{k_2(\ell_4-q)+d_{k_2,\ell_3,\ell_4}}{k_2}\right)}\\
\times\int_0^{\tau_1^{k_1}}\int_0^{\tau_2^{k_2}}(\tau_1^{k_1}-s_1)^{\frac{d_{k_1,\ell_1,\ell_2}}{k_1}-1}
(\tau_2^{k_2}-s_2)^{\frac{k_2(\ell_4-q)+d_{k_2,\ell_3,\ell_4}}{k_2}-1}(k_1s_1)^{\ell_2}(k_2s_2)^{q}R_{\underline{\ell}}(im_1)\\
\times\omega(s_1^{1/k_1},s_2^{1/k_2},m_1,\epsilon)\frac{ds_2}{s_2}\frac{ds_1}{s_1},
\end{multline}
\begin{multline}\label{e320}
\mathcal{A}_4(\tau_1,\tau_2,m_1,\epsilon)=\sum_{1\le h\le \ell_2-1,1\le q\le \ell_4-1}A_{\ell_2,h}A_{\ell_4,q}\frac{\tau_1^{k_1}}
{\Gamma\left(\frac{k_1(\ell_2-h)+d_{k_1,\ell_1,\ell_2}}{k_1}\right)}\frac{\tau_2^{k_2}}{\Gamma\left(\frac{k_2(\ell_4-q)+d_{k_2,\ell_3,\ell_4}}{k_2}\right)}\\
\times\int_0^{\tau_1^{k_1}}\int_0^{\tau_2^{k_2}}(\tau_1^{k_1}-s_1)^{\frac{k_1(\ell_2-h)+d_{k_1,\ell_1,\ell_2}}{k_1}-1}
(\tau_2^{k_2}-s_2)^{\frac{k_2(\ell_4-q)+d_{k_2,\ell_3,\ell_4}}{k_2}-1}(k_1s_1)^{h}(k_2s_2)^{q}R_{\underline{\ell}}(im_1)\\
\times\omega(s_1^{1/k_1},s_2^{1/k_2},m_1,\epsilon)\frac{ds_2}{s_2}\frac{ds_1}{s_1}.
\end{multline}

\vspace{0.3cm}

\textbf{Remark:}  The term $\mathcal{A}_1$ in the previous expression corresponds to
$$T_1^{d_{k_1,\ell_1,\ell_2}}(T_1^{k_1+1}\partial_{T_1})^{\ell_2}T_2^{d_{k_2,\ell_3,\ell_4}}(T_2^{k_2+1}\partial_{T_2})^{\ell_4}R_{\underline{\ell}}(\partial_z)U(T_1,T_2,z,\epsilon).$$

The term $\mathcal{A}_2$ corresponds to
$$\sum_{1\le h\le \ell_2-1}A_{\ell_2,h}T_1^{d_{k_1,\ell_1,\ell_2}+k_1(\ell_2-h)}(T_1^{k_1+1}\partial_{T_1})^{h}T_2^{d_{k_2,\ell_3,\ell_4}}(T_2^{k_2+1}\partial_{T_2})^{\ell_4}R_{\underline{\ell}}(\partial_z)U(T_1,T_2,z,\epsilon).$$

The term $\mathcal{A}_3$ comes from 
$$\sum_{1\le q\le \ell_4-1}A_{\ell_4,q}T_1^{d_{k_1,\ell_1,\ell_2}}(T_1^{k_1+1}\partial_{T_1})^{\ell_2}T_2^{d_{k_2,\ell_3,\ell_4}+k_2(\ell_4-q)}(T_2^{k_2+1}\partial_{T_2})^{q}R_{\underline{\ell}}(\partial_z)U(T_1,T_2,z,\epsilon).$$

Finally, the term $\mathcal{A}_4$ emerges from
\begin{multline*}
\sum_{1\le h\le \ell_2-1,1\le q\le \ell_4-1}A_{\ell_2,h}A_{\ell_4,q}T_1^{d_{k_1,\ell_1,\ell_2}+k_1(\ell_2-h)}(T_1^{k_1+1}\partial_{T_1})^{h}T_2^{d_{k_2,\ell_3,\ell_4}+k_2(\ell_4-q)}(T_2^{k_2+1}\partial_{T_2})^{q}\\R_{\underline{\ell}}(\partial_z)U(T_1,T_2,z,\epsilon).
\end{multline*}

\subsection{Solution to an auxiliary problem. I}\label{secaux1}

In this section, we preserve the assumptions made considering the main problem (\ref{epral}), in Section~\ref{secmainproblem}. We seek for a solution to the auxiliary problem (\ref{e309}) belonging to certain Banach spaces of functions, described in Section~\ref{sec452}. This goal is attained by means of the estimation obtained in the next result.

For every $m\in\R$, we consider the polynomial
\begin{equation}\label{e377}
P_m(\tau_1,\tau_2)=Q(im)-R(im)(k_1\tau_1^{k_1})^{\delta_{1}}(k_2\tau_2^{k_2})^{\delta_{2}}.
\end{equation}

\begin{lemma}\label{lema365}
There exist $\rho_1,\rho_2>0$ and $C_1>0$ such that 
$$|P_m(\tau_1,\tau_2)|\ge C_1|R(im)|,\quad (\tau_1,\tau_2,m)\in D(0,\rho_1)\times D(0,\rho_2)\times \R.$$
\end{lemma}
\begin{proof}
In view of the assumptions made in (\ref{e120}) and (\ref{e121}) on the polynomials $Q,R$, we arrive at the existence of $r_{Q,R}>0$ such that 
$$\left|\frac{Q(im)}{R(im)}\right|\ge r_{Q,R},\qquad m\in\R.$$
Let us choose small enough $\rho_1,\rho_2>0$ such that 
$$(k_1\rho_1^{k_1})^{\delta_1}(k_2\rho_2^{k_2})^{\delta_2}\le\frac{r_{Q,R}}{2}.$$ 
This entails that 
$$\left|\frac{Q(im)}{R(im)}-(k_1\tau_1^{k_1})^{\delta_1}(k_2\tau_2^{k_2})^{\delta_2}\right|\ge \frac{r_{Q,R}}{2}$$
for every $m\in\R$.

For all $m\in\R$ and for all $\tau_j\in D(0,\rho_j)$, $j=1,2$, let us write
\begin{equation}\label{e378}
P_m(\tau_1,\tau_2)=R(im)\left[\frac{Q(im)}{R(im)}-(k_1\tau_1^{k_1})^{\delta_1}(k_2\tau_2^{k_2})^{\delta_2}\right],
\end{equation}
which allow us to conclude the result for $C_1=\frac{r_{Q,R}}{2}$.
\end{proof}

\begin{prop}\label{prop1}
Let $\rho_1,\rho_2$ be prescribed as in Lemma~\ref{lema365}. For every $\varpi>0$ there exists $\varsigma_{F,1}>0$ such that if
$$C_{\Psi}\le \varsigma_{F,1},$$
(see (\ref{e215}) for the value of $C_{\Psi}$) then the problem (\ref{e309}) admits a unique solution $\omega_{\rho_1,\rho_2}(\tau_1,\tau_2,m,\epsilon)$, continuous on $\R$ with respect to its third variable, and holomorphic with respect to $(\tau_1,\tau_2,\epsilon)$ on $D(0,\rho_1)\times D(0,\rho_2)\times D(0,\epsilon_0)$, such that 
\begin{equation}\label{e401}
|\omega_{\rho_1,\rho_2}(\tau_1,\tau_2,m,\epsilon)|\le \varpi\frac{1}{(1+|m|)^{\mu}}e^{-\beta|m|}\left|\tau_1\tau_2\right|,
\end{equation}
for every $(\tau_1,\tau_2,m,\epsilon)\in D(0,\rho_1)\times D(0,\rho_2)\times\R\times D(0,\epsilon_0)$.
\end{prop}
\begin{proof}
Let $\epsilon\in D(0,\epsilon_0)$ and consider the Banach space $B_{(\beta,\mu,\rho_1,\rho_2)}$ studied in Section~\ref{sec452}. We fix $\varpi>0$.

Let us denote by $\overline{B}(0,\varpi)\subseteq B_{(\beta,\mu,\rho_1,\rho_2)}$ the set of elements in $B_{(\beta,\mu,\rho_1,\rho_2)}$ whose norm is upper bounded by $\varpi$. We consider the operator 
\begin{multline*}
\mathcal{H}_{\epsilon}(\omega(\tau_1,\tau_2,m)):=\frac{1}{P_m(\tau_1,\tau_2)}\\
\times\left[\sum_{\underline{\ell}=(\ell_1,\ell_2,\ell_3,\ell_4)\in I}\epsilon^{\Delta_{\underline{\ell}}-\ell_1+\ell_2-\ell_3+\ell_4}\frac{1}{(2\pi)^{1/2}}\int_{-\infty}^{\infty}C_{\underline{\ell}}(m-m_1,\epsilon)\left[\mathcal{A}_1(\tau_1,\tau_2,m_1,\epsilon)+\mathcal{A}_2(\tau_1,\tau_2,m_1,\epsilon)\right.\right.\\
\left.+\mathcal{A}_3(\tau_1,\tau_2,m_1,\epsilon)+\mathcal{A}_4(\tau_1,\tau_2,m_1,\epsilon)\right]dm_1+\tau_1^{k_1}\tau_2^{k_2}\int_0^{\tau_1^{k_1}}\int_0^{\tau_2^{k_2}}\frac{1}{(2\pi)^{1/2}}\int_{-\infty}^{\infty}P_1(\epsilon,i(m-m_1))\\
\times\omega\left((\tau_1^{k_1}-s_1)^{1/k_1},(\tau_2^{k_2}-s_2)^{1/k_2},m-m_1\right)P_2(\epsilon,im_1)\omega(s_1^{1/k_1},s_2^{1/k_2},m_1)\\
\left.\times\frac{1}{(\tau_1^{k_1}-s_1)}\frac{1}{s_1}\frac{1}{(\tau_2^{k_2}-s_2)}\frac{1}{s_2}ds_2ds_1dm_1+\Psi(\tau_1,\tau_2,m,\epsilon)\right],
\end{multline*}
where $P_m$ is defined in (\ref{e377}), and $\mathcal{A}_j$ for $1\le j\le 4$ are defined in (\ref{e317})-(\ref{e320}), where the term $\omega\left(s_1^{1/k_1},s_2^{1/k_2},m_1,\epsilon\right)$ needs to be replaced by $\omega\left(s_1^{1/k_1},s_2^{1/k_2},m_1\right)$.

Regarding the assumptions made in (\ref{e121}), there exists $C_{R}>0$ such that
$$|R(im)|\ge C_{R}(1+|m|)^{\hbox{deg}(R)},\quad m\in\R.$$
Let $\omega(\tau_1,\tau_2,m)\in \overline{B}(0,\varpi)$. From the assumptions made in (\ref{e136}) and (\ref{e137}) together with (\ref{e111}), one can apply Lemma~\ref{lema365}, Proposition~\ref{prop477}, Proposition~\ref{prop441}, Proposition~\ref{prop487}, together with (\ref{e166}), one has that
\begin{multline*}
\left\|\mathcal{H}_{\epsilon}(\omega(\tau_1,\tau_2,m))\right\|_{(\beta,\mu,\rho_1,\rho_2)}\le\frac{1}{C_1C_R}\left[\sum_{\underline{\ell}=(\ell_1,\ell_2,\ell_3,\ell_4)\in I}\epsilon_0^{\Delta_{\underline{\ell}}-\ell_1+\ell_2-\ell_3+\ell_4}\frac{1}{(2\pi)^{1/2}}K\right.\\
\left.\times\tilde{C}_1\left[\mathfrak{C}_1+\mathfrak{C}_2+\mathfrak{C}_3+\mathfrak{C}_4\right]\varpi+\tilde{C}_2\frac{1}{(2\pi)^{1/2}}\varpi^2+ C_{\Psi}\right]
\end{multline*}
with $K>0$ being the constant involved in (\ref{e166}), $\tilde{C}_1$ and $\tilde{C}_2$ the positive constants appearing in Proposition~\ref{prop441} and Proposition~\ref{prop487}, respectively, and with $C_{\Psi}$ in (\ref{e215}), together with
$$\mathfrak{C}_1=\frac{1}{\Gamma\left(\frac{d_{k_1,\ell_1,\ell_2}}{k_1}\right)\Gamma\left(\frac{d_{k_2,\ell_3,\ell_4}}{k_2}\right)}k_1^{\ell_2}k_2^{\ell_4},$$
$$\mathfrak{C}_2=\sum_{1\le h\le \ell_2-1}|A_{\ell_2,h}|\frac{1}{\Gamma\left(\frac{k_1(\ell_2-h)+d_{k_1,\ell_1,\ell_2}}{k_1}\right)\Gamma\left(\frac{d_{k_2,\ell_3,\ell_4}}{k_2}\right)}k_1^{h}k_2^{\ell_4},$$
$$\mathfrak{C}_3=\sum_{1\le q\le \ell_4-1}|A_{\ell_4,q}|\frac{1}
{\Gamma\left(\frac{d_{k_1,\ell_1,\ell_2}}{k_1}\right)\Gamma\left(\frac{k_2(\ell_4-q)+d_{k_2,\ell_3,\ell_4}}{k_2}\right)}k_1^{\ell_2}k_2^{q},$$
$$\mathfrak{C}_4=\sum_{1\le h\le \ell_2-1,1\le q\le \ell_4-1}|A_{\ell_2,h}||A_{\ell_4,q}|\frac{1}
{\Gamma\left(\frac{k_1(\ell_2-h)+d_{k_1,\ell_1,\ell_2}}{k_1}\right)\Gamma\left(\frac{k_2(\ell_4-q)+d_{k_2,\ell_3,\ell_4}}{k_2}\right)}k_1^{h}k_2^{q}.
$$
It is worth remarking that the previous constants are linked to the conditions imposed on the parameters involved in the problem, (\ref{e167}) and (\ref{e114}). In this concern, the bounds $K\tilde{C}_1\mathfrak{C}_1$ are linked to Proposition~\ref{prop441} and the fact that 
$$k_1+\frac{d_{k_1,\ell_1,\ell_2}}{k_1}-1+\ell_2-1+\frac{1}{k_1}\ge0,$$
together with 
$$k_2+\frac{d_{k_2,\ell_3,\ell_4}}{k_2}-1+\ell_4-1+\frac{1}{k_2}\ge0.$$ 
In addition to this, the bounds $K\tilde{C}_1\mathfrak{C}_2$ are linked to Proposition~\ref{prop441} and the fact that 
$$k_1+\frac{k_1(\ell_2-h)+d_{k_1,\ell_1,\ell_2}}{k_1}-1+h-1+\frac{1}{k_1}\ge0,$$
for $1\le h\le \ell_2-1$, together with 
$$k_2+\frac{d_{k_2,\ell_3,\ell_4}}{k_2}-1+\ell_4-1+\frac{1}{k_2}\ge0.$$ 
Also, the bounds $K\tilde{C}_1\mathfrak{C}_3$ are linked to Proposition~\ref{prop441} and the fact that 
$$k_1+\frac{d_{k_1,\ell_1,\ell_2}}{k_1}-1+\ell_2-1+\frac{1}{k_1}\ge0,$$
and
$$k_2+\frac{k_2(\ell_4-q)+d_{k_2,\ell_3,\ell_4}}{k_2}-1+q-1+\frac{1}{k_2}\ge0,$$ 
for $1\le q\le \ell_4-1$. Finally, the bounds $K\tilde{C}_1\mathfrak{C}_4$ are linked to Proposition~\ref{prop441} and the fact that 
$$k_1+\frac{k_1(\ell_2-h)+d_{k_1,\ell_1,\ell_2}}{k_1}-1+h-1+\frac{1}{k_1}\ge0,$$
for $1\le h\le \ell_2-1$, with 
$$k_2+\frac{k_2(\ell_4-q)+d_{k_2,\ell_3,\ell_4}}{k_2}-1+q-1+\frac{1}{k_2}\ge0,$$ 
for $1\le q\le \ell_4-1$. 
 
Let us choose $\varsigma_{F,1}>0$, $\rho_1,\rho_2$ and $\epsilon_0>0$ small enough such that
\begin{equation}\label{e433}
\frac{1}{C_1C_{R}}\left(\sum_{\underline{\ell}=(\ell_1,\ell_2,\ell_3,\ell_4)\in I}\epsilon_0^{\Delta_{\underline{\ell}}-\ell_1+\ell_2-\ell_3+\ell_4}\frac{K}{(2\pi)^{1/2}}\tilde{C}_1\left[\mathfrak{C}_1+\mathfrak{C}_2+\mathfrak{C}_3+\mathfrak{C}_4\right]+2\tilde{C}_2\frac{1}{(2\pi)^{1/2}}\varpi\right)\le\frac{1}{2},
\end{equation}
and
$$\frac{1}{C_1C_R}\varsigma_{F,1}\le\frac{1}{2}\varpi.$$
Recall that $\tilde{C}_2$ can be taken close to 0 provided that $\rho_1,\rho_2$ are reduced, in view of Proposition~\ref{prop487}. 

The previous inequality entails that 
$$\mathcal{H}_{\epsilon}(\omega(\tau_1,\tau_2,m))\in \overline{B}(0,\varpi).$$

On the other hand, given $\omega_1,\omega_2\in \overline{B}(0,\varpi)\subseteq B_{(\beta,\mu,\rho_1,\rho_2)}$, it holds that

\begin{multline*}
\mathcal{H}_{\epsilon}(\omega_1(\tau_1,\tau_2,m))-\mathcal{H}_{\epsilon}(\omega_2(\tau_1,\tau_2,m)):=\frac{1}{P_m(\tau_1,\tau_2)}\\
\times\left[\sum_{\underline{\ell}=(\ell_1,\ell_2,\ell_3,\ell_4)\in I}\epsilon^{\Delta_{\underline{\ell}}-\ell_1+\ell_2-\ell_3+\ell_4}\frac{1}{(2\pi)^{1/2}}\int_{-\infty}^{\infty}C_{\underline{\ell}}(m-m_1,\epsilon)\left[\mathcal{A}^{\star}_1(\tau_1,\tau_2,m_1)\right.\right.\\
\left.+\mathcal{A}^{\star}_2(\tau_1,\tau_2,m_1)+\mathcal{A}^{\star}_3(\tau_1,\tau_2,m_1)+\mathcal{A}^{\star}_4(\tau_1,\tau_2,m_1)\right]dm_1
+\frac{\tau_1^{k_1}\tau_2^{k_2}}{(2\pi)^{1/2}}\int_0^{\tau_1^{k_1}}\int_0^{\tau_2^{k_2}}\int_{-\infty}^{\infty}\\
\left\{P_1(\epsilon,i(m-m_1))\omega_1\left((\tau_1^{k_1}-s_1)^{1/k_1},(\tau_2^{k_2}-s_2)^{1/k_2},m-m_1\right)P_2(\epsilon,im_1)\omega_1(s_1^{1/k_1},s_2^{1/k_2},m_1)\right.\\
\left.- P_1(\epsilon,i(m-m_1))\omega_2\left((\tau_1^{k_1}-s_1)^{1/k_1},(\tau_2^{k_2}-s_2)^{1/k_2},m-m_1\right)P_2(\epsilon,im_1)\omega_2(s_1^{1/k_1},s_2^{1/k_2},m_1)\right\}\\
\left.\times\frac{1}{\tau_1^{k_1}-s_1}\frac{1}{s_1}\frac{1}{\tau_2^{k_2}-s_2}\frac{1}{s_2}ds_2ds_1dm_1\right],
\end{multline*}
with
\begin{multline}\label{e317b}
\mathcal{A}^{\star}_1(\tau_1,\tau_2,m_1)=\frac{\tau_1^{k_1}}{\Gamma\left(\frac{d_{k_1,\ell_1,\ell_2}}{k_1}\right)}\frac{\tau_2^{k_2}}{\Gamma\left(\frac{d_{k_2,\ell_3,\ell_4}}{k_2}\right)}\int_0^{\tau_1^{k_1}}\int_0^{\tau_2^{k_2}}(\tau_1^{k_1}-s_1)^{\frac{d_{k_1,\ell_1,\ell_2}}{k_1}-1}
(\tau_2^{k_2}-s_2)^{\frac{d_{k_2,\ell_3,\ell_4}}{k_2}-1}\\
\times(k_1s_1)^{\ell_2}(k_2s_2)^{\ell_4}R_{\underline{\ell}}(im_1)(\omega_1(s_1^{1/k_1},s_2^{1/k_2},m_1)-\omega_2(s_1^{1/k_1},s_2^{1/k_2},m_1))\frac{ds_2}{s_2}\frac{ds_1}{s_1},
\end{multline}
\begin{multline}\label{e318b}
\mathcal{A}^{\star}_2(\tau_1,\tau_2,m_1)=\sum_{1\le h\le \ell_2-1}A_{\ell_2,h}\frac{\tau_1^{k_1}}{\Gamma\left(\frac{k_1(\ell_2-h)+d_{k_1,\ell_1,\ell_2}}{k_1}\right)}\frac{\tau_2^{k_2}}{\Gamma\left(\frac{d_{k_2,\ell_3,\ell_4}}{k_2}\right)}\\
\times\int_0^{\tau_1^{k_1}}\int_0^{\tau_2^{k_2}}(\tau_1^{k_1}-s_1)^{\frac{k_1(\ell_2-h)+d_{k_1,\ell_1,\ell_2}}{k_1}-1}
(\tau_2^{k_2}-s_2)^{\frac{d_{k_2,\ell_3,\ell_4}}{k_2}-1}(k_1s_1)^{h}(k_2s_2)^{\ell_4}R_{\underline{\ell}}(im_1)\\
\times(\omega_1(s_1^{1/k_1},s_2^{1/k_2},m_1)-\omega_2(s_1^{1/k_1},s_2^{1/k_2},m_1))\frac{ds_2}{s_2}\frac{ds_1}{s_1},
\end{multline}
\begin{multline}\label{e319b}
\mathcal{A}^{\star}_3(\tau_1,\tau_2,m_1)=\sum_{1\le q\le \ell_4-1}A_{\ell_4,q}\frac{\tau_1^{k_1}}
{\Gamma\left(\frac{d_{k_1,\ell_1,\ell_2}}{k_1}\right)}\frac{\tau_2^{k_2}}{\Gamma\left(\frac{k_2(\ell_4-q)+d_{k_2,\ell_3,\ell_4}}{k_2}\right)}\\
\times\int_0^{\tau_1^{k_1}}\int_0^{\tau_2^{k_2}}(\tau_1^{k_1}-s_1)^{\frac{d_{k_1,\ell_1,\ell_2}}{k_1}-1}
(\tau_2^{k_2}-s_2)^{\frac{k_2(\ell_4-q)+d_{k_2,\ell_3,\ell_4}}{k_2}-1}(k_1s_1)^{\ell_2}(k_2s_2)^{q}R_{\underline{\ell}}(im_1)\\
\times(\omega_1(s_1^{1/k_1},s_2^{1/k_2},m_1)-\omega_2(s_1^{1/k_1},s_2^{1/k_2},m_1))\frac{ds_2}{s_2}\frac{ds_1}{s_1},
\end{multline}
\begin{multline}\label{e320b}
\mathcal{A}^{\star}_4(\tau_1,\tau_2,m_1)=\sum_{1\le h\le \ell_2-1,1\le q\le \ell_4-1}A_{\ell_2,h}A_{\ell_4,q}\frac{\tau_1^{k_1}}
{\Gamma\left(\frac{k_1(\ell_2-h)+d_{k_1,\ell_1,\ell_2}}{k_1}\right)}\frac{\tau_2^{k_2}}{\Gamma\left(\frac{k_2(\ell_4-q)+d_{k_2,\ell_3,\ell_4}}{k_2}\right)}\\
\times\int_0^{\tau_1^{k_1}}\int_0^{\tau_2^{k_2}}(\tau_1^{k_1}-s_1)^{\frac{d_{k_1,\ell_1,\ell_2}+k_1(\ell_2-h)}{k_1}-1}
(\tau_2^{k_2}-s_2)^{\frac{k_2(\ell_4-q)+d_{k_2,\ell_3,\ell_4}}{k_2}-1}(k_1s_1)^{h}(k_2s_2)^{q}R_{\underline{\ell}}(im_1)\\
\times(\omega_1(s_1^{1/k_1},s_2^{1/k_2},m_1)-\omega_2(s_1^{1/k_1},s_2^{1/k_2},m_1))\frac{ds_2}{s_2}\frac{ds_1}{s_1}.
\end{multline}

Finally, we check that the expression
\begin{multline*}
P_1(\epsilon,i(m-m_1))\omega_1\left((\tau_1^{k_1}-s_1)^{1/k_1},(\tau_2^{k_2}-s_2)^{1/k_2},m-m_1\right)P_2(\epsilon,im_1)\omega_1(s_1^{1/k_1},s_2^{1/k_2},m_1)\\
- P_1(\epsilon,i(m-m_1))\omega_2\left((\tau_1^{k_1}-s_1)^{1/k_1},(\tau_2^{k_2}-s_2)^{1/k_2},m-m_1\right)P_2(\epsilon,im_1)\omega_2(s_1^{1/k_1},s_2^{1/k_2},m_1)
\end{multline*}
equals
\begin{multline*}
 P_1(\epsilon,i(m-m_1))\hfill\\
\times(\omega_1\left((\tau_1^{k_1}-s_1)^{1/k_1},(\tau_2^{k_2}-s_2)^{1/k_2},m-m_1\right)-\omega_2\left((\tau_1^{k_1}-s_1)^{1/k_1},(\tau_2^{k_2}-s_2)^{1/k_2},m-m_1\right))\\
\hfill\times P_2(\epsilon,im_1)\omega_1(s_1^{1/k_1},s_2^{1/k_2},m_1)\\
+ P_1(\epsilon,i(m-m_1))\omega_2((\tau_1^{k_1}-s_1)^{1/k_1},(\tau_2^{k_2}-s_2)^{1/k_2},m-m_1)\\
\times P_2(\epsilon,im_1)(\omega_1(s_1^{1/k_1},s_2^{1/k_2},m_1)-\omega_2(s_1^{1/k_1},s_2^{1/k_2},m_1)).
\end{multline*}

An analogous reasoning as before and (\ref{e433}) allow us to arrive at

\begin{multline*}
\left\|\mathcal{H}_{\epsilon}(\omega_1(\tau_1,\tau_2,m))-\mathcal{H}_{\epsilon}(\omega_2(\tau_1,\tau_2,m))\right\|_{(\beta,\mu,\rho_1,\rho_2)}\\
\le\frac{1}{C_1C_{R}}\left[ \sum_{\underline{\ell}=(\ell_1,\ell_2,\ell_3,\ell_4)\in I}\epsilon_0^{\Delta_{\underline{\ell}}-\ell_1+\ell_2-\ell_3+\ell_4}\frac{1}{(2\pi)^{1/2}}K\tilde{C}_1\left[\mathfrak{C}_1+\mathfrak{C}_2+\mathfrak{C}_3+\mathfrak{C}_4\right]\right.\\
\times\left\|\omega_1(\tau_1,\tau_2,m)-\omega_2(\tau_1,\tau_2,m)\right\|_{(\beta,\mu,\rho_1,\rho_2)}\hfill\\
+\tilde{C}_2\frac{1}{(2\pi)^{1/2}}\left\|\omega_1(\tau_1,\tau_2,m)-\omega_2(\tau_1,\tau_2,m)\right\|_{(\beta,\mu,\rho_1,\rho_2)}\left\|\omega_1(\tau_1,\tau_2,m)\right\|_{(\beta,\mu,\rho_1,\rho_2)}\\
\left.+ \tilde{C}_2\frac{1}{(2\pi)^{1/2}}\left\|\omega_2(\tau_1,\tau_2,m)\right\|_{(\beta,\mu,\rho_1,\rho_2)}\left\|\omega_1(\tau_1,\tau_2,m)-\omega_2(\tau_1,\tau_2,m)\right\|_{(\beta,\mu,\rho_1,\rho_2)}
\right]\\
\le \frac{1}{C_1C_{R}}\left[ \sum_{\underline{\ell}=(\ell_1,\ell_2,\ell_3,\ell_4)\in I}\epsilon_0^{\Delta_{\underline{\ell}}-\ell_1+\ell_2-\ell_3+\ell_4}\frac{1}{(2\pi)^{1/2}}K\tilde{C}_1\left[\mathfrak{C}_1+\mathfrak{C}_2+\mathfrak{C}_3+\mathfrak{C}_4\right]\right.\\
\left.+2\tilde{C}_2\frac{1}{(2\pi)^{1/2}}\varpi\right]\left\|\omega_1(\tau_1,\tau_2,m)-\omega_2(\tau_1,\tau_2,m)\right\|_{(\beta,\mu,\rho_1,\rho_2)}\\
\le \frac{1}{2}\left\|\omega_1(\tau_1,\tau_2,m)-\omega_2(\tau_1,\tau_2,m)\right\|_{(\beta,\mu,\rho_1,\rho_2)}.
\end{multline*}

The classical contractive mapping theorem can be applied on $\mathcal{H}_{\epsilon}:\overline{B}(0,\varpi)\to \overline{B}(0,\varpi)$, to arrive at the existence of a unique fixed point, say $\omega_{\epsilon}(\tau_1,\tau_2,m)$. This construction depends holomorphically on $\epsilon\in D(0,\epsilon_0)$. As a conclusion, one arrives at the existence of a function 
$$(\tau_1,\tau_2,m,\epsilon)\mapsto \omega_{\rho_1,\rho_2}(\tau_1,\tau_2,m,\epsilon)=\omega_{\epsilon}(\tau_1,\tau_2,m),$$ 
which is, by construction, a solution of (\ref{e309}), and satisfies the estimates (\ref{e401}). 
\end{proof}

\subsection{Solution to an auxiliary problem. II}\label{secaux2}

In this section, we still preserve the assumptions made regarding the main problem (\ref{epral}), in Section~\ref{secmainproblem} and search for solutions to the auxiliary problem (\ref{e309}) in the  Banach space of Section~\ref{sec453}. 

We start by providing alternative lower bounds on $P_m$ to those attained in Lemma~\ref{lema365}.

\begin{lemma}\label{lema365b}
There exist $d_1,d_2\in\R$ and $C_2>0$ such that 
$$|P_m(\tau_1,\tau_2)|\ge C_2|R(im)|(1+|k_1\tau_1^{k_1}|^{\delta_{1}}|k_2\tau_2^{k_2}|^{\delta_{2}}),$$
for every $(\tau_1,\tau_2,m)\in S_{d_1}\times S_{d_2}\times \R$.
Here, $P_m$ is the polynomial defined by (\ref{e377}), and $S_{d_j}$, for $j=1,2$ stands for some infinite sector centered at the origin, and bisecting direction $d_j$.
\end{lemma}
\begin{proof}
We recall from the proof of Lemma~\ref{lema365}, together with assumption (\ref{e133}) the existence of $r_{Q,R}>0$ such that
\begin{equation}\label{e538}
\nabla:=\left\{\frac{Q(im)}{R(im)}:m\in\R\right\}\subseteq S_{Q,R}\setminus D(0,r_{Q,R}).
\end{equation} 
Let $d_1,d_2\in\R$ be chosen such that 
$$d_{1,2}:=\delta_{1}k_1d_1+\delta_{2}k_2d_2\not\in\hbox{arg}(\nabla),$$
where $\hbox{arg}(\nabla)=\{\hbox{arg}(z):z\in\nabla\}$. We define $S_{d_1}$ (resp. $S_{d_2}$) an infinite sector with small opening, centered at the origin, and bisecting direction $d_1$ (resp. $d_2$). The infinite sector $S_{d_{1,2}}$ is defined accordingly. This entails in particular that $(k_1\tau_1^{k_1})^{\delta_1}(k_2\tau_2^{k_2})^{\delta_2}\not\in\nabla$ for $(\tau_1,\tau_2)\in S_{d_1}\times S_{d_2}$ when the opening of $S_{d_1}, S_{d_2}$ are close enough to 0. More precisely, for every $\xi\in S_{d_{1,2}}$ one guarantees the existence of a positive constant $C_2>0$ with
$$\left|\frac{Q(im)}{R(im)}-\xi\right|\ge C_2(1+|\xi|).$$
In particular, one has that
$$\left|\frac{Q(im)}{R(im)}-(k_1\tau_1^{k_1})^{\delta_1}(k_2\tau_2^{k_2})^{\delta_2}\right|\ge C_2(1+|k_1\tau_1^{k_1}|^{\delta_1}|k_2\tau_2^{k_2}|^{\delta_2}).$$
The result follows from the decomposition (\ref{e378}).
\end{proof}


\begin{prop}
Let $\nu_1,\nu_2>0$ and choose $d_1,d_2\in\R$ as in Lemma~\ref{lema365b}. Then, there exists $D_{\Psi}>0$ such that
$$\sup_{(\tau_1,\tau_2,m)\in S_{d_1}\times S_{d_2}\times \R}\left|\frac{\Psi(\tau_1,\tau_2,m,\epsilon)}{P_m(\tau_1,\tau_2)}\right|(1+|m|)^{\mu}e^{\beta|m|}\frac{1+|\tau_1|^{2k_1}}{|\tau_1|}\frac{1+|\tau_2|^{2k_2}}{|\tau_2|}\exp\left(-\nu_1|\tau_1|^{k_1}-\nu_2|\tau_2|^{k_2}\right)\le D_{\Psi},$$
valid for all $\epsilon\in D(0,\epsilon_0)$.
\end{prop}
\begin{proof}
Taking into account the definition of $\Psi$, together with Lemma~\ref{lema365b} and (\ref{e121}), one has that
\begin{multline*}
\sup_{(\tau_1,\tau_2,m)\in S_{d_1}\times S_{d_2}\times \R}\left|\frac{\Psi(\tau_1,\tau_2,m,\epsilon)}{P_m(\tau_1,\tau_2)}\right|(1+|m|)^{\mu}e^{\beta|m|}\frac{1+|\tau_1|^{2k_1}}{|\tau_1|}\frac{1+|\tau_2|^{2k_2}}{|\tau_2|}\exp\left(-\nu_1|\tau_1|^{k_1}-\nu_2|\tau_2|^{k_2}\right)\\
\le \frac{\tilde{K}}{C_2}\frac{1}{\min_{m\in\R}|R(im)|}\sum_{n_1\in N_1}\frac{1}{\Gamma(n_1/k_1)}\left(\sup_{\tau_1\in S_{d_1}}|\tau_1|^{n_1-1}(1+|\tau_1|^{2k_1})\exp(-\nu_1|\tau_1|^{k_1})\right)\\
\times \sum_{n_2\in N_2}\frac{1}{\Gamma(n_2/k_2)}\left(\sup_{\tau_2\in S_{d_2}}|\tau_2|^{n_2-1}(1+|\tau_2|^{2k_2})\exp(-\nu_2|\tau_2|^{k_2})\right),
\end{multline*}
where $\tilde{K}$ is the positive constant appearing in (\ref{e280}). The last expression is finite and does not depend on $\epsilon\in D(0,\epsilon_0)$.
\end{proof}

\vspace{0.3cm}

\textbf{Remark:}
In terms of the Banach space described in Section~\ref{sec453}, the previous result can be read as follows: there exists $D_{\Psi}>0$ such that for all $\epsilon\in D(0,\epsilon_0)$, the function $\frac{\Psi(\tau_1,\tau_2,m,\epsilon)}{P_m(\tau_1,\tau_2)}$ belongs to $E_{(\beta,\mu,\nu_1,\nu_2,S_{d_1,S_{d_2}})}$, and 
$$\sup_{\epsilon\in D(0,\epsilon_0)}\left\|\frac{\Psi(\tau_1,\tau_2,m,\epsilon)}{P_m(\tau_1,\tau_2)}\right\|_{(\beta,\mu,\nu_1,\nu_2,S_{d_1,S_{d_2}})}\le D_{\Psi}.$$
Observe from the proof of the previous result that the constant $D_{\Psi}$ tends to 0 when $\tilde{K}$ approaches to 0.

\vspace{0.3cm}

\begin{prop}\label{prop2}
For every $\varpi>0$ there exists $\varsigma_{F,2}>0$ such that if
$$D_{\Psi}\le \varsigma_{F,2},$$
provided that the constants $C_{P_1}, C_{P_2}$ appearing in (\ref{e706}) are small enough, then the problem (\ref{e309}) admits a unique solution $\omega_{S_{d_1},S_{d_2}}(\tau_1,\tau_2,m,\epsilon)$, continuous on $\R$ with respect to its third variable, and holomorphic with respect to $(\tau_1,\tau_2,\epsilon)$ on $S_{d_1}\times S_{d_2}\times D(0,\epsilon_0)$, such that 
\begin{multline}\label{e401b}
|\omega_{S_{d_1},S_{d_2}}(\tau_1,\tau_2,m,\epsilon)|\le \varpi\frac{1}{(1+|m|)^{\mu}}\\
\times e^{-\beta|m|}\frac{\left|\tau_1\right|}{1+\left|\tau_1\right|^{2k_1}}\frac{\left|\tau_2\right|}{1+\left|\tau_2\right|^{2k_2}}\exp\left(\nu_1\left|\tau_1\right|^{k_1}+\nu_2\left|\tau_2\right|^{k_2}\right),
\end{multline}
for every $(\tau_1,\tau_2,m,\epsilon)\in S_{d_1}\times S_{d_2} \times \R \times D(0,\epsilon_0)$.
\end{prop}
\begin{proof}
The proof of Proposition~\ref{prop2} follows the same line of arguments and similar notations as the one of Proposition~\ref{prop1} owing to the conditions (\ref{e167}). Indeed, an inequality for $\left\|\mathcal{H}_{\epsilon}(\omega(\tau_1,\tau_2,m))\right\|_{(\beta,\mu,\nu_1,\nu_2,S_{d_1},S_{d_2})}$ is reached similar to the one obtained, where $C_{\Psi}$ is replaced by $D_{\Psi}$, and the constant $C_1$ is replaced by $C_2$, together with $K\tilde{C}_1$ and $\tilde{C}_2$ appearing in Proposition~\ref{prop441b} and Proposition~\ref{prop487b}. Furthermore, similar inequalities to (\ref{e433}) and for the difference $\left\|\mathcal{H}_{\epsilon}(\omega_1)-\mathcal{H}_{\epsilon}(\omega_2)\right\|_{(\beta,\mu,\nu_1,\nu_2,S_{d_1},S_{d_2})}$ are achieved, provided that $\epsilon_0>0$ is taken close to 0 (bearing in mind the condition (\ref{e111})) and that $C_{\rho_1},C_{\rho_2}>0$ are small enough, according to the remark after Proposition~\ref{prop441b}.
\end{proof}

\subsection{Solution to an auxiliary problem, III. Analytic continuation}\label{secaux3}

We maintain the assumptions made on the elements involved in the main equation (\ref{epral}),  now searching for solutions in a third Banach space, and exploring its analytic continuation by means of the results obtained so far. 

Taking into account the lower estimates obtained in Lemma~\ref{lema365}, and an analogous line of arguments as those followed in Proposition~\ref{prop1}, we arrive at the following result.

\begin{prop}\label{prop3}
For every $\varpi>0$ there exists $\varsigma_{F,3}>0$ such that if
$$C_{\Psi}\le \varsigma_{F,3},$$
then the auxiliary problem (\ref{e309}) admits a unique solution $\omega_{\rho_1,\rho_2,S_{d_1},S_{d_2}}(\tau_1,\tau_2,m,\epsilon)$, continuous on $\R$ with respect to its third variable, and holomorphic with respect to $(\tau_1,\tau_2,\epsilon)$ on $(S_{d_1}\cap D(0,\rho_1))\times (S_{d_2}\cap D(0,\rho_2))\times D(0,\epsilon_0)$, such that 
\begin{equation}\label{e401c}
|\omega_{\rho_1,\rho_2,S_{d_1},S_{d_2}}(\tau_1,\tau_2,m,\epsilon)|\le \varpi\frac{1}{(1+|m|)^{\mu}} e^{-\beta|m|}\left|\tau_1\tau_2\right|,
\end{equation}
for every $(\tau_1,\tau_2,m,\epsilon)\in (S_{d_1}\cap D(0,\rho_1))\times (S_{d_2}\cap D(0,\rho_2))\times \R\times D(0,\epsilon_0)$.
\end{prop}

The next result proves that the solutions of (\ref{e309}) obtained in Proposition~\ref{prop1} and Proposition~\ref{prop2} are related by analytic continuation, by means of the solution constructed in Proposition~\ref{prop3}.

\begin{prop}\label{propauxpral}
Let $\varpi>0$. There exists $\varsigma_{F}>0$ such that if $C_{\Psi}\le \varsigma_{F}$ and $D_{\Psi}\le \varsigma_{F}$, then the following statements hold:
\begin{itemize}
\item For every fixed $\tau_1\in S_{d_1}\cap D(0,\rho_1)$, $m\in\R$ and $\epsilon\in D(0,\epsilon_0)$, the map
$$\tau_2\mapsto \omega_{S_{d_1},S_{d_2}}(\tau_1,\tau_2,m,\epsilon)$$
defined  on $S_{d_2}$ (see Proposition~\ref{prop2}) has an analytic continuation on $D(0,\rho_2)$, which is $\tau_2\mapsto \omega_{\rho_1,\rho_2}(\tau_1,\tau_2,m,\epsilon)$.
\item For every fixed $\tau_2\in S_{d_2}\cap D(0,\rho_2)$, $m\in\R$ and $\epsilon\in D(0,\epsilon_0)$, the map
$$\tau_1\mapsto \omega_{S_{d_1},S_{d_2}}(\tau_1,\tau_2,m,\epsilon)$$
defined  on $S_{d_1}$ (see Proposition~\ref{prop2}) has an analytic continuation on $D(0,\rho_1)$, which is $\tau_1\mapsto \omega_{\rho_1,\rho_2}(\tau_1,\tau_2,m,\epsilon)$.
\end{itemize}
\end{prop}
\begin{proof}
We fix
$$\Delta:=\max_{(\tau_1,\tau_2)\in D(0,\rho_1)\times D(0,\rho_2)}\frac{1}{1+\left|\tau_1\right|^{2k_1}}\frac{1}{1+\left|\tau_2\right|^{2k_1}}\exp(\nu_1 \left|\tau_1\right|^{k_1}+\nu_2 \left|\tau_2\right|^{k_2}).$$
We observe that $\Delta>0$. Let $(\tau_1,\tau_2,m,\epsilon)\mapsto\omega_{S_{d_1},S_{d_2}}(\tau_1,\tau_2,m,\epsilon)$ be the function obtained in Proposition~\ref{prop2}, which solves (\ref{e309}) and such that (\ref{e401b}) holds for some given $\varpi$, if $D_{\Psi}\le \varsigma_{F,2}$. Then, for $(\tau_1,\tau_2,m,\epsilon)\in (S_{d_1}\cap D(0,\rho_1))\times (S_{d_2}\cap D(0,\rho_2))\times \R\times D(0,\epsilon_0)$, one has that
\begin{multline*}
|\omega_{S_{d_1},S_{d_2}}(\tau_1,\tau_2,m,\epsilon)|\le \varpi\frac{1}{(1+|m|)^{\mu}}\\
\times e^{-\beta|m|}\frac{\left|\tau_1\right|}{1+\left|\tau_1\right|^{2k_1}}\frac{\left|\tau_2\right|}{1+\left|\tau_2\right|^{2k_2}}\exp\left(\nu_1\left|\tau_1\right|^{k_1}+\nu_2\left|\tau_2\right|^{k_2}\right)\\
\le \varpi \Delta\frac{1}{(1+|m|)^{\mu}}e^{-\beta|m|}\left|\tau_1\right|\left|\tau_2\right|.
\end{multline*}
This entails that for all $\epsilon\in D(0,\epsilon_0)$, the function 
$$(S_{d_1}\cap D(0,\rho_1))\times (S_{d_2}\cap D(0,\rho_2))\times \R \ni(\tau_1,\tau_2,m)\mapsto\omega_{S_{d_1},S_{d_2}}(\tau_1,\tau_2,m,\epsilon)$$
is the fixed point of the operator $\mathcal{H}_{\epsilon}$ when defined on the closed ball  $\overline{B}(0,\Delta\varpi)$  of the Banach space $F_{(\beta,\mu,\rho_1,\rho_2,S_{d_1},S_{d_2})}$ of Section~\ref{sec924}. From unicity of the fixed point for $\mathcal{H}_{\epsilon}$ in such ball obtained in Proposition~\ref{prop3}, $\omega_{\rho_1,\rho_2,S_{d_1},S_{d_2}}$, we conclude that 
$$\omega_{S_{d_1},S_{d_2}}(\tau_1,\tau_2,m,\epsilon)=\omega_{\rho_1,\rho_2,S_{d_1},S_{d_2}}(\tau_1,\tau_2,m,\epsilon)$$
for every $(\tau_1,\tau_2,m)\in (S_{d_1}\cap D(0,\rho_1))\times (S_{d_2}\cap D(0,\rho_2))\times \R$.  An analogous reasoning leads us to 
$$\omega_{\rho_1,\rho_2}(\tau_1,\tau_2,m,\epsilon)=\omega_{\rho_1,\rho_2,S_{d_1},S_{d_2}}(\tau_1,\tau_2,m,\epsilon),$$
for every $(\tau_1,\tau_2,m)\in (S_{d_1}\cap D(0,\rho_1))\times (S_{d_2}\cap D(0,\rho_2))\times \R$, with $\omega_{\rho_1,\rho_2}$ being the fixed point of $\mathcal{H}_{\epsilon}$, defined on the ball $\overline{B}(0,\Delta\varsigma)$ of the Banach space $B_{(\beta,\mu,\rho_1,\rho_2)}$, obtained in Proposition~\ref{prop1}. 

The result follows from the variation of $\epsilon\in D(0,\epsilon_0)$.
\end{proof}

At this point, we can state the main result of the present section, summarizing all the previous results.

\begin{theo}\label{teo1}
Under the assumptions of Section~\ref{secmainproblem}, we consider the Cauchy problem (\ref{epral}). Let $\rho_1,\rho_2>0$ determined in Section~\ref{secaux1}, and $d_1,d_2\in\R$ chosen in Section~\ref{secaux2}.

Let $\mathcal{E}\subseteq D(0,\epsilon_0)$ and $\mathcal{T}_1,\mathcal{T}_2\subseteq D(0,r_{\mathcal{T}})$ for some small enough $r_{\mathcal{T}}>0$ be three bounded sectors with vertex at the origin, chosen in such a way that
\begin{itemize}
\item There exists $\Delta_1>0$ with $\cos(k_1(d_1-\hbox{arg}(\epsilon t_1)))>\Delta_1$,
for all $\epsilon\in\mathcal{E}$ and $t_1\in\mathcal{T}_1$.
\item There exists $\Delta_2>0$ with $\cos(k_2(d_2-\hbox{arg}(\epsilon t_2)))>\Delta_2$,
for all $\epsilon\in\mathcal{E}$ and $t_1\in\mathcal{T}_2$.
\end{itemize}

Then, provided that $\epsilon_0$, the quantity $\tilde{K}$ from (\ref{e280}), together with the constants $C_{P_1},C_{P_2}>0$ from (\ref{e706}) are taken small enough, and for every $0<\beta'<\beta$, the problem (\ref{epral}), under null initial data $u(0,t_2,z,\epsilon)=u(t_1,0,z,\epsilon)=0$, admits an analytic solution $u_{d_1,d_2}(t_1,t_2,z,\epsilon)\in\mathcal{O}_b(\mathcal{T}_1\times\mathcal{T}_2\times H_{\beta'}\times \mathcal{E})$.
\end{theo}

\begin{proof}
Let $\omega_{S_{d_1},S_{d_2}}(\tau_1,\tau_2,m,\epsilon)$ be the solution of the auxiliary equation (\ref{e309}), obtained in Proposition~\ref{prop2}. Regarding (\ref{e401b}), the choice of $d_1,d_2$ at the statement of the result guarantees that the function 
\begin{multline}
u_{d_1,d_2}(t_1,t_2,z,\epsilon)=\frac{k_1k_2}{(2\pi)^{1/2}}\int_{-\infty}^{\infty}\int_{L_{d_1}}\int_{L_{d_2}}\omega_{S_{d_1},S_{d_2}}(u_1,u_2,m,\epsilon)\\
\times \exp\left(-\left(\frac{u_1}{\epsilon t_1}\right)^{k_1}-\left(\frac{u_2}{\epsilon t_2}\right)^{k_2}\right)e^{izm}\frac{du_2}{u_2}\frac{du_1}{u_1}dm,\label{e776}
\end{multline}
is well-defined, holomorphic and bounded on $\mathcal{T}_1\times\mathcal{T}_2\times H_{\beta'}\times \mathcal{E}$, provided that $r_{\mathcal{T}},\epsilon_0>0$ are small enough, and $0<\beta'<\beta$. Indeed, observe that for all $(t_1,t_2,z,\epsilon)\in\mathcal{T}_1\times\mathcal{T}_2\times H_{\beta'}\times \mathcal{E}$, one has that
\begin{multline*}
|u_{d_1,d_2}(t_1,t_2,z,\epsilon)|\le \varpi\frac{k_1k_2}{(2\pi)^{1/2}}\left[\int_{-\infty}^{\infty}\frac{1}{(1+|m|)^{\mu}}e^{-|m|(\beta-|\hbox{Im}(z)|)}dm\right]\\
\times\prod_{j=1}^{2}\left(\int_{0}^{\infty}\frac{r_j}{1+r_j^{2k_j}}\exp\left(r_j^{k_j}(\nu_j-\frac{\Delta_j}{|\epsilon t_j|^{k_j}})\right)dr_j\right),
\end{multline*}
for some $\varpi>0$, after the parametrization $u_j=r_je^{id_j}$ for $r_j\in[0,\infty)$ and $j=1,2$. Assuming that $\epsilon_0 r_{\mathcal{T}}<(\Delta_j/\nu_j)^{1/k_j}$, for $j=1,2$, the previous integrals converge.
\end{proof}

\section{Parametric Gevrey series expansion}\label{secpar}

In this section, we provide an asymptotic representation of the analytic solution, obtained in the previous section. We maintain the assumptions made on the elements in the construction of the main problem (\ref{epral}), stated in Section~\ref{secmainproblem}.

We split the integral representation of the solution to problem (\ref{epral}) obtained in Theorem~\ref{teo1} as the sum 
\begin{equation}\label{e792}
u_{d_1,d_2}(t_1,t_2,z,\epsilon)=J_1(t_1,t_2,z,\epsilon)+J_2(t_1,t_2,z,\epsilon)+J_3(t_1,t_2,z,\epsilon),
\end{equation}
with
\begin{multline*}
J_1(t_1,t_2,z,\epsilon):=\frac{k_1k_2}{(2\pi)^{1/2}}\int_{-\infty}^{\infty}\int_{L_{d_1},\rho_1/2}\int_{L_{d_2},\rho_2/2}\omega_{S_{d_1},S_{d_2}}(u_1,u_2,m,\epsilon)\nonumber\\
\times \exp\left(-\left(\frac{u_1}{\epsilon t_1}\right)^{k_1}-\left(\frac{u_2}{\epsilon t_2}\right)^{k_2}\right)e^{izm}\frac{du_2}{u_2}\frac{du_1}{u_1}dm,
\end{multline*}
\begin{multline*}
J_2(t_1,t_2,z,\epsilon):=\frac{k_1k_2}{(2\pi)^{1/2}}\int_{-\infty}^{\infty}\int_{L_{d_1},\rho_1/2}\int_{L_{d_2},\rho_2/2,\infty}\omega_{S_{d_1},S_{d_2}}(u_1,u_2,m,\epsilon)\nonumber\\
\times \exp\left(-\left(\frac{u_1}{\epsilon t_1}\right)^{k_1}-\left(\frac{u_2}{\epsilon t_2}\right)^{k_2}\right)e^{izm}\frac{du_2}{u_2}\frac{du_1}{u_1}dm,
\end{multline*}
and
\begin{multline*}
J_3(t_1,t_2,z,\epsilon):=\frac{k_1k_2}{(2\pi)^{1/2}}\int_{-\infty}^{\infty}\int_{L_{d_1},\rho_1/2,\infty}\int_{L_{d_2}}\omega_{S_{d_1},S_{d_2}}(u_1,u_2,m,\epsilon)\nonumber\\
\times \exp\left(-\left(\frac{u_1}{\epsilon t_1}\right)^{k_1}-\left(\frac{u_2}{\epsilon t_2}\right)^{k_2}\right)e^{izm}\frac{du_2}{u_2}\frac{du_1}{u_1}dm,
\end{multline*}

where $L_{d_1,\rho_1/2}=[0,\frac{\rho_1}{2}]e^{\sqrt{-1}d_1}$, $L_{d_2,\rho_2/2}=[0,\frac{\rho_2}{2}]e^{\sqrt{-1}d_2}$, $L_{d_1,\rho_1/2,\infty}=[\frac{\rho_1}{2},\infty)e^{\sqrt{-1}d_1}$, $L_{d_2,\rho_2/2,\infty}=[\frac{\rho_2}{2},\infty)e^{\sqrt{-1}d_2}$.

Our objective is to study the asymptotic Gevrey related to each piece in the previous decomposition.

\subsection{Gevrey expansions for $J_1$}\label{secj1}

Let us recall the notion of a good covering, which will be essential in our reasoning. 

\begin{defin}\label{defi662}
Let $\varsigma\ge2$ be an integer. Let $\underline{\mathcal{E}}=(\mathcal{E}_p)_{0\le p\le \varsigma-1}$ be a set of bounded sectors with vertex at the origin, $\mathcal{E}_p\subseteq D(0,\epsilon_0)$, for $0\le p\le \varsigma-1$ such that $\mathcal{E}_{p}\cap \mathcal{E}_{p+1}\neq\emptyset$ for $0\le p\le\varsigma-1$ (with the notation $\mathcal{E}_\varsigma:=\mathcal{E}_0$), which are three by three disjoint, i.e. $\mathcal{E}_{p_1}\cap \mathcal{E}_{p_2}\cap \mathcal{E}_{p_3}=\emptyset$ for all $0\le p_1,p_2,p_3\le\varsigma-1$, with $p_1\neq p_2\neq p_3$ and $p_1\neq p_3$. In addition to this, there exists a neighborhood of the origin $\mathcal{U}$ such that $\mathcal{U}\setminus\{0\}=\cup_{p=0}^{\varsigma-1}\mathcal{E}_p$.

In this situation, we say the family $\underline{\mathcal{E}}$ determines a good covering in $\C^{\star}$.

\end{defin}

Let us depart from given bounded open sectors $\mathcal{T}_1,\mathcal{T}_2,\mathcal{E}$ with vertex at the origin, and $\rho_1,\rho_2>0$ and $d_1,d_2\in\R$, under the hypotheses of Theorem~\ref{teo1}. We choose a good covering $\underline{\mathcal{E}}=(\mathcal{E}_p)_{0\le p\le \varsigma-1}$ such that $\mathcal{E}_0:=\mathcal{E}$. In addition to this, we choose the real numbers $\mathfrak{d}_p,\tilde{\mathfrak{d}}_p$, for $0\le p\le \varsigma-1$ with
$$\mathfrak{d}_0:=d_1,\qquad \tilde{\mathfrak{d}}_0=d_2,$$
in such a way that the following conditions hold:
\begin{itemize}
\item For every $0\le p\le \varsigma-1$ there exists $\nabla_p>0$ such that $\cos(k_1(\mathfrak{d}_p-\hbox{arg}(\epsilon t_1)))>\nabla_p$, for $\epsilon\in\mathcal{E}_p$, $t_1\in\mathcal{T}_1$, and
\item For every $0\le p\le \varsigma-1$ there exists $\tilde{\nabla}_p>0$ such that $\cos(k_2(\tilde{\mathfrak{d}}_p-\hbox{arg}(\epsilon t_2)))>\tilde{\nabla}_p$, for $\epsilon\in\mathcal{E}_p$, $t_2\in\mathcal{T}_2$. 
\end{itemize}

Observe that one can choose $\nabla_0=\Delta_1$ and $\tilde{\nabla}_0=\Delta_2$, where $\Delta_1,\Delta_2$ are the constants in the statements of Theorem~\ref{teo1}.

For every $0\le p\le \varsigma-1$, we construct the function
\begin{multline}
J_{1,p}(t_1,t_2,z,\epsilon):=\frac{k_1k_2}{(2\pi)^{1/2}}\int_{-\infty}^{\infty}\int_{L_{\mathfrak{d}_p},\rho_1/2}\int_{L_{\tilde{\mathfrak{d}}_p},\rho_2/2}\omega_{\rho_1,\rho_2}(u_1,u_2,m,\epsilon)\\
\hfill\times \exp\left(-\left(\frac{u_1}{\epsilon t_1}\right)^{k_1}-\left(\frac{u_2}{\epsilon t_2}\right)^{k_2}\right)e^{izm}\frac{du_2}{u_2}\frac{du_1}{u_1}dm,\label{e685}
\end{multline}
with $L_{\mathfrak{d}_p,\rho_1/2}=[0,\rho_1/2]e^{i\mathfrak{d}_p}$, $L_{\tilde{\mathfrak{d}}_p,\rho_2/2}=[0,\rho_2/2]e^{i\tilde{\mathfrak{d}}_p}$, which turns out to be an analytic and bounded function on $\mathcal{T}_1\times\mathcal{T}_2\times H_{\beta'}\times\mathcal{E}_p$, for every $0<\beta'<\beta$.

Observe that $J_{1,0}(t_1,t_2,z,\epsilon)=J_{1}(t_1,t_2,z,\epsilon)$ for all $\epsilon\in\mathcal{E}=\mathcal{E}_0$, $t_1\in\mathcal{T}_1$, $t_2\in\mathcal{T}_2$, $z\in H_{\beta'}$.

The next Proposition provides bounds for the differences of consecutive maps $J_{1,p}$.

\begin{prop}\label{prop689}
Under the previous assumptions, the following statements hold for every $0\le p\le \varsigma-1$:
\begin{itemize}
\item Case 1: $\mathfrak{d}_p=\mathfrak{d}_{p+1}$, and $\tilde{\mathfrak{d}}_p\neq\tilde{\mathfrak{d}}_{p+1}$. There exist $C_{p,1},C_{p,2}>0$ such that
$$|J_{1,p+1}(t_1,t_2,z,\epsilon)-J_{1,p}(t_1,t_2,z,\epsilon)|\le C_{p,1}\exp\left(-\frac{C_{p,2}}{|\epsilon|^{k_2}}\right),$$
for all $(t_1,t_2,z,\epsilon)\in \mathcal{T}_1\times\mathcal{T}_2\times H_{\beta'}\times (\mathcal{E}_p\cap\mathcal{E}_{p+1})$.
\item Case 2: $\mathfrak{d}_p\neq\mathfrak{d}_{p+1}$, and $\tilde{\mathfrak{d}}_p=\tilde{\mathfrak{d}}_{p+1}$. There exist $C_{p,3},C_{p,4}>0$ such that
$$|J_{1,p+1}(t_1,t_2,z,\epsilon)-J_{1,p}(t_1,t_2,z,\epsilon)|\le C_{p,3}\exp\left(-\frac{C_{p,4}}{|\epsilon|^{k_1}}\right),$$
for all $(t_1,t_2,z,\epsilon)\in \mathcal{T}_1\times\mathcal{T}_2\times H_{\beta'}\times (\mathcal{E}_p\cap\mathcal{E}_{p+1})$.
\item Case 3: $\mathfrak{d}_p\neq\mathfrak{d}_{p+1}$, and $\tilde{\mathfrak{d}}_p\neq\tilde{\mathfrak{d}}_{p+1}$. There exist $C_{p,5},C_{p,6}>0$ such that
$$|J_{1,p+1}(t_1,t_2,z,\epsilon)-J_{1,p}(t_1,t_2,z,\epsilon)|\le C_{p,5}\exp\left(-\frac{C_{p,6}}{|\epsilon|^{k_2}}\right),$$
for all $(t_1,t_2,z,\epsilon)\in \mathcal{T}_1\times\mathcal{T}_2\times H_{\beta'}\times (\mathcal{E}_p\cap\mathcal{E}_{p+1})$.
\end{itemize}
\end{prop}
\begin{proof}

Let $0\le p\le \varsigma-1$.

First, assume that $\mathfrak{d}_p=\mathfrak{d}_{p+1}$ and $\tilde{\mathfrak{d}}_p\neq\tilde{\mathfrak{d}}_{p+1}$. Choose $0<\beta'<\beta$. Then, one has that for all $(t_1,t_2,z,\epsilon)\in \mathcal{T}_1\times\mathcal{T}_2\times H_{\beta'}\times (\mathcal{E}_p\cap \mathcal{E}_{p+1})$,
$$J_{1,p+1}(t_1,t_2,z,\epsilon)-J_{1,p}(t_1,t_2,z,\epsilon)=\frac{k_1k_2}{(2\pi)^{1/2}}\int_{-\infty}^{\infty}\int_{L_{\mathfrak{d}_p},\rho_1/2}I_1(u_1,m,\epsilon)\exp\left(-\left(\frac{u_1}{\epsilon t_1}\right)^{k_1}\right)e^{izm}\frac{du_1}{u_1}dm,$$

\begin{multline}\label{e702}
I_1(u_1,m,\epsilon)=\int_{L_{\tilde{\mathfrak{d}}_{p+1},\rho_2/2}}\omega_{\rho_1,\rho_2}(u_1,u_2,m,\epsilon)\exp\left(-\left(\frac{u_2}{\epsilon t_2}\right)^{k_2}\right)\frac{d u_2}{u_2}\\
- \int_{L_{\tilde{\mathfrak{d}}_{p},\rho_2/2}}\omega_{\rho_1,\rho_2}(u_1,u_2,m,\epsilon)\exp\left(-\left(\frac{u_2}{\epsilon t_2}\right)^{k_2}\right)\frac{d u_2}{u_2}.
\end{multline}

We recall from Proposition~\ref{prop1} that the map $u_2\mapsto \omega_{\rho_1,\rho_2}(u_1,u_2,m,\epsilon)$ is holomorphic on the disc $D(0,\rho_2)$. Therefore, one can apply Cauchy's theorem to deform the integration path in $I(u_1,m,\epsilon)$ to arrive at
$$I_1(u_1,m,\epsilon)=\int_{C_{\tilde{\mathfrak{d}}_{p},\tilde{\mathfrak{d}}_{p+1},\rho_2/2}}\omega_{\rho_1,\rho_2}(u_1,u_2,m,\epsilon)\exp\left(-\left(\frac{u_2}{\epsilon t_2}\right)^{k_2}\right)\frac{d u_2}{u_2},$$
where $C_{\tilde{\mathfrak{d}}_{p},\tilde{\mathfrak{d}}_{p+1},\rho_2/2}$ is the arc of circle joining the points $\rho_2/2e^{\sqrt{-1}\tilde{\mathfrak{d}}_{p}}$ and $\rho_2/2e^{\sqrt{-1}\tilde{\mathfrak{d}}_{p+1}}$. Taking into account the bounds in Proposition~\ref{prop1}, we have that
\begin{align*}
|I_1(u_1,m,\epsilon)|&\le\left|u_1\right|\frac{1}{(1+|m|)^{\mu}}e^{-\beta|m|}\varpi_{\rho_1,\rho_2}\left|\int_{\tilde{\mathfrak{d}}_p}^{\tilde{\mathfrak{d}}_{p+1}}\exp\left(-\left(\frac{\rho_2/2}{|\epsilon t_2|}\right)^{k_2}\tilde{\Delta}_{p,p+1}\right)\frac{\rho_2}{2}d\theta\right|\\
&\le \left|u_1\right|\frac{1}{(1+|m|)^{\mu}}e^{-\beta|m|}\varpi_{\rho_1,\rho_2}\frac{\rho_2}{2}|\tilde{\mathfrak{d}}_{p+1}-\tilde{\mathfrak{d}}_{p}|\exp\left(-\left(\frac{\rho_2/2}{r_{\mathcal{T}_2}}\right)^{k_2}\tilde{\Delta}_{p,p+1}\frac{1}{|\epsilon|^{k_2}}\right),
\end{align*}
for some $\varpi_{\rho_1,\rho_2},\tilde{\Delta}_{p,p+1}>0$, valid for all $\epsilon\in \mathcal{E}_{p}\cap\mathcal{E}_{p+1}$, $u_1\in L_{\mathfrak{d}_p,\rho_1/2}$ and $m\in\R$. Here, $r_{\mathcal{T}_2}>0$ stands for the radius of $\mathcal{T}_2$.

This entails that 
\begin{multline*}
|J_{1,p+1}(t_1,t_2,z,\epsilon)-J_{1,p}(t_1,t_2,z,\epsilon)|\le\frac{k_1k_2}{(2\pi)^{1/2}}\left(\int_{-\infty}^{\infty}e^{-(\beta-\beta')|m|}\frac{1}{(1+|m|)^{\mu}}dm\right)\varpi_{\rho_1,\rho_2}\\
\times E_1(|\epsilon t_1|)\frac{\rho_2}{2}|\tilde{\mathfrak{d}}_{p+1}-\tilde{\mathfrak{d}}_{p}|\exp\left(-\left(\frac{\rho_2/2}{r_{\mathcal{T}_2}}\right)^{k_2}\tilde{\Delta}_{p,p+1}\frac{1}{|\epsilon|^{k_2}}\right),
\end{multline*}
with
$$E_{1}(|\epsilon t_1|)=\int_0^{\rho_1/2}\exp\left(-\left(\frac{r_1}{|\epsilon t_1|}\right)^{k_1}\right)dr_1\le \frac{\rho_1}{2}.$$

As a result, one concludes the existence of two constants $C_{p,1},C_{p,2}>0$ such that
$$|J_{1,p+1}(t_1,t_2,z,\epsilon)-J_{1,p}(t_1,t_2,z,\epsilon)|\le C_{p,1}\exp\left(-\frac{C_{p,2}}{|\epsilon|^{k_2}}\right),$$
for all $(t_1,t_2,z,\epsilon)\in \mathcal{T}_1\times\mathcal{T}_2\times H_{\beta'}\times (\mathcal{E}_p\cap\mathcal{E}_{p+1})$.

The second case can be handled in an analogous way as the first one. Indeed, one has that for all $(t_1,t_2,z,\epsilon)\in \mathcal{T}_1\times\mathcal{T}_2\times H_{\beta'}\times (\mathcal{E}_p\cap \mathcal{E}_{p+1})$,
$$J_{1,p+1}(t_1,t_2,z,\epsilon)-J_{1,p}(t_1,t_2,z,\epsilon)=\frac{k_1k_2}{(2\pi)^{1/2}}\int_{-\infty}^{\infty}\int_{L_{\tilde{\mathfrak{d}}_2},\rho_2/2}I_2(u_2,m,\epsilon)\exp\left(-\left(\frac{u_2}{\epsilon t_2}\right)^{k_2}\right)e^{izm}\frac{du_2}{u_2}dm,$$
where
\begin{multline*}
I_2(u_2,m,\epsilon)=\int_{L_{\mathfrak{d}_{p+1},\rho_1/2}}\omega_{\rho_1,\rho_2}(u_1,u_2,m,\epsilon)\exp\left(-\left(\frac{u_1}{\epsilon t_1}\right)^{k_1}\right)\frac{d u_1}{u_1}\\
- \int_{L_{\mathfrak{d}_{p},\rho_1/2}}\omega_{\rho_1,\rho_2}(u_1,u_2,m,\epsilon)\exp\left(-\left(\frac{u_1}{\epsilon t_1}\right)^{k_1}\right)\frac{d u_1}{u_1}.
\end{multline*}

Once again, one can deform the integration path appearing in $I_2(u_2,m,\epsilon)$  due to $u_1\mapsto \omega_{\rho_1,\rho_2}(u_1,u_2,m,\epsilon)$ is a holomorphic function on $D(0,\rho_1)$. We get the representation
$$I_2(u_2,m,\epsilon)=\int_{C_{\mathfrak{d}_{p},\mathfrak{d}_{p+1},\rho_1/2}}\omega_{\rho_1,\rho_2}(u_1,u_2,m,\epsilon)\exp\left(-\left(\frac{u_1}{\epsilon t_1}\right)^{k_1}\right)\frac{d u_1}{u_1}.$$
Proposition~\ref{prop1} yields
\begin{align*}
|I_2(u_2,m,\epsilon)|&\le\left|u_2\right|\frac{1}{(1+|m|)^{\mu}}e^{-\beta|m|}\varpi_{\rho_1,\rho_2}\left|\int_{\mathfrak{d}_p}^{\mathfrak{d}_{p+1}}\exp\left(-\left(\frac{\rho_1/2}{|\epsilon t_1|}\right)^{k_1}\Delta_{p,p+1}\right)\frac{\rho_1}{2}d\theta\right|\\
&\le \left|u_2\right|\frac{1}{(1+|m|)^{\mu}}e^{-\beta|m|}\varpi_{\rho_1,\rho_2}\frac{\rho_1}{2}|\mathfrak{d}_{p+1}-\mathfrak{d}_{p}|\exp\left(-\left(\frac{\rho_1/2}{r_{\mathcal{T}_1}}\right)^{k_1}\Delta_{p,p+1}\frac{1}{|\epsilon|^{k_1}}\right),
\end{align*}
with $\epsilon\in \mathcal{E}_{p}\cap\mathcal{E}_{p+1}$, $u_2\in L_{\tilde{\mathfrak{d}}_p,\rho_2/2}$ and $m\in\R$. Here, $r_{\mathcal{T}_1}>0$ stands for the radius of $\mathcal{T}_1$.

\begin{multline*}
|J_{1,p+1}(t_1,t_2,z,\epsilon)-J_{1,p}(t_1,t_2,z,\epsilon)|\le\frac{k_1k_2}{(2\pi)^{1/2}}\left(\int_{-\infty}^{\infty}e^{-(\beta-\beta')|m|}\frac{1}{(1+|m|)^{\mu}}dm\right)\varpi_{\rho_1,\rho_2}\\
\times E_2(|\epsilon t_2|)\frac{\rho_1}{2}|\mathfrak{d}_{p+1}-\mathfrak{d}_{p}|\exp\left(-\left(\frac{\rho_1/2}{r_{\mathcal{T}_1}}\right)^{k_2}\Delta_{p,p+1}\frac{1}{|\epsilon|^{k_1}}\right),
\end{multline*}
with
$$E_{2}(|\epsilon t_2|)=\int_0^{\rho_2/2}\exp\left(-\left(\frac{r_2}{|\epsilon t_2|}\right)^{k_2}\right)dr_2\le\frac{\rho_2}{2}.$$

This leads to the existence of $C_{p,3},C_{p,4}>0$ such that
$$|J_{1,p+1}(t_1,t_2,z,\epsilon)-J_{1,p}(t_1,t_2,z,\epsilon)|\le C_{p,3}\exp\left(-\frac{C_{p,4}}{|\epsilon|^{k_1}}\right),$$
for all $(t_1,t_2,z,\epsilon)\in \mathcal{T}_1\times\mathcal{T}_2\times H_{\beta'}\times (\mathcal{E}_p\cap\mathcal{E}_{p+1})$.

The third case deals with $\mathfrak{d}_p\neq\mathfrak{d}_{p+1}$, and $\tilde{\mathfrak{d}}_p\neq\tilde{\mathfrak{d}}_{p+1}$. Let us consider the following deformation of the integration paths involved in the integrals defining the analytic maps $J_{1,p}$ and $J_{1,p+1}$, with respect to their first variable, $t_1$.

Let $j\in\{p,p+1\}$. We define $\theta_p=(\mathfrak{d}_p+\mathfrak{d}_{p+1})/2$, and consider the concatenation of paths 
$$L_{\theta_{p,p+1},0,\rho_1/4}+C_{\theta_{p,p+1},\mathfrak{d}_j,\rho_1/4}+L_{\mathfrak{d}_j,\rho_1/4,\rho_1/2},$$
(see Figure~\ref{fig1}), with $L_{\theta_{p,p+1},0,\rho_1/4}=[0,\rho_1/4]e^{i\theta_{p,p+1}}$, $C_{\theta_{p,p+1},\mathfrak{d}_j,\rho_1/4}$ is the path parametrized by $[\theta_{p,p+1},\mathfrak{d}_j]\ni\theta\mapsto \rho_1/4 e^{i\theta}$, and $L_{\mathfrak{d}_j,\rho_1/4,\rho_1/2}=[\rho_1/4,\rho_1/2]e^{i\mathfrak{d}_j}$.

\begin{figure}
	\centering
		\includegraphics[width=0.3\textwidth]{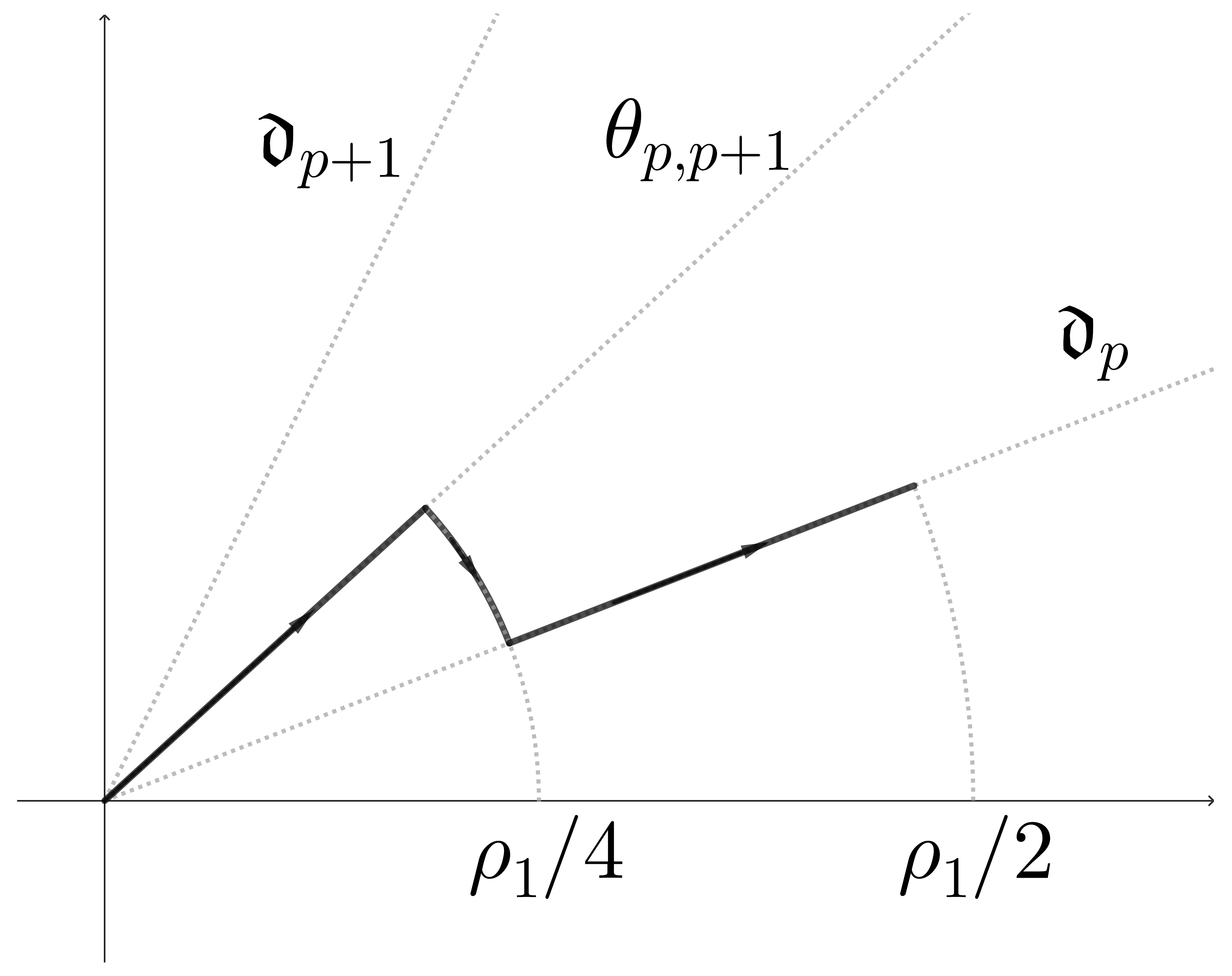}
		\caption{Concatenation of paths, for $j=p$}
		\label{fig1}
\end{figure}

We define the paths $L_{\tilde{\theta}_{p,p+1},0,\rho_1/4}$, $C_{\tilde{\theta}_{p,p+1},\tilde{\mathfrak{d}}_j,\rho_1/4}$, and $L_{\tilde{\mathfrak{d}}_j,\rho_1/4,\rho_1/2}$ in an analogous manner, and consider the concatenation of paths
$$L_{\tilde{\theta}_{p,p+1},0,\rho_1/4}+C_{\tilde{\theta}_{p,p+1},\tilde{\mathfrak{d}}_j,\rho_1/4}+L_{\tilde{\mathfrak{d}}_j,\rho_1/4,\rho_1/2}.$$

We observe that, by means of Cauchy theorem, one can write the difference of two consecutive solutions, by splitting the integrals involved, in the form
$$J_{1,p+1}(t_1,t_2,z,\epsilon)-J_{1,p}(t_1,t_2,z,\epsilon)=I_1+I_2+I_3-I_4-I_5-I_6,$$
after deforming $L_{\mathfrak{d}_j,\rho_1/2}$ and $L_{\tilde{\mathfrak{d}}_j,\rho_1/2}$ arriving to the previous two concatenations of paths. Here,
\begin{multline*}
I_1=\frac{k_1k_2}{(2\pi)^{1/2}}\int_{-\infty}^{\infty}\int_{L_{\mathfrak{d}_{p+1},\rho_1/4,\rho_1/2}}\int_{L_{\tilde{\mathfrak{d}}_{p+1},\rho_2/2}}\omega_{\rho_1,\rho_2}(u_1,u_2,m,\epsilon)\\
\times \exp\left(-\left(\frac{u_1}{\epsilon t_1}\right)^{k_1}-\left(\frac{u_2}{\epsilon t_2}\right)^{k_2}\right)e^{izm}\frac{du_2}{u_2}\frac{du_1}{u_1}dm,
\end{multline*}

\begin{multline*}
I_2=\frac{k_1k_2}{(2\pi)^{1/2}}\int_{-\infty}^{\infty}\int_{C_{\theta_{p,p+1},\mathfrak{d}_{p+1},\rho_1/4}}\int_{L_{\tilde{\mathfrak{d}}_{p+1},\rho_2/2}}\omega_{\rho_1,\rho_2}(u_1,u_2,m,\epsilon)\\
\times \exp\left(-\left(\frac{u_1}{\epsilon t_1}\right)^{k_1}-\left(\frac{u_2}{\epsilon t_2}\right)^{k_2}\right)e^{izm}\frac{du_2}{u_2}\frac{du_1}{u_1}dm,
\end{multline*}

\begin{multline*}
I_3=\frac{k_1k_2}{(2\pi)^{1/2}}\int_{-\infty}^{\infty}\int_{L_{\theta_{p,p+1},0,\rho_1/4}}\int_{L_{\tilde{\mathfrak{d}}_{p+1},\rho_2/2}}\omega_{\rho_1,\rho_2}(u_1,u_2,m,\epsilon)\\
\times \exp\left(-\left(\frac{u_1}{\epsilon t_1}\right)^{k_1}-\left(\frac{u_2}{\epsilon t_2}\right)^{k_2}\right)e^{izm}\frac{du_2}{u_2}\frac{du_1}{u_1}dm,
\end{multline*}

\begin{multline*}
I_4=\frac{k_1k_2}{(2\pi)^{1/2}}\int_{-\infty}^{\infty}\int_{L_{\mathfrak{d}_{p},\rho_1/4,\rho_1/2}}\int_{L_{\tilde{\mathfrak{d}}_{p},\rho_2/2}}\omega_{\rho_1,\rho_2}(u_1,u_2,m,\epsilon)\\
\times \exp\left(-\left(\frac{u_1}{\epsilon t_1}\right)^{k_1}-\left(\frac{u_2}{\epsilon t_2}\right)^{k_2}\right)e^{izm}\frac{du_2}{u_2}\frac{du_1}{u_1}dm,
\end{multline*}

\begin{multline*}
I_5=\frac{k_1k_2}{(2\pi)^{1/2}}\int_{-\infty}^{\infty}\int_{C_{\theta_{p,p+1},\mathfrak{d}_{p},\rho_1/4}}\int_{L_{\tilde{\mathfrak{d}}_{p},\rho_2/2}}\omega_{\rho_1,\rho_2}(u_1,u_2,m,\epsilon)\\
\times \exp\left(-\left(\frac{u_1}{\epsilon t_1}\right)^{k_1}-\left(\frac{u_2}{\epsilon t_2}\right)^{k_2}\right)e^{izm}\frac{du_2}{u_2}\frac{du_1}{u_1}dm,
\end{multline*}

\begin{multline*}
I_6=\frac{k_1k_2}{(2\pi)^{1/2}}\int_{-\infty}^{\infty}\int_{L_{\theta_{p,p+1},0,\rho_1/4}}\int_{L_{\tilde{\mathfrak{d}}_{p},\rho_2/2}}\omega_{\rho_1,\rho_2}(u_1,u_2,m,\epsilon)\\
\times \exp\left(-\left(\frac{u_1}{\epsilon t_1}\right)^{k_1}-\left(\frac{u_2}{\epsilon t_2}\right)^{k_2}\right)e^{izm}\frac{du_2}{u_2}\frac{du_1}{u_1}dm.
\end{multline*}

We proceed to upper estimate $I_1,\ldots,I_6$.

Let us first consider $I_3-I_6$. We write
$$I_3-I_6=\frac{k_1k_2}{(2\pi)^{1/2}}\int_{-\infty}^{\infty}\int_{L_{\theta_{p,p+1},0,\rho_1/4}} K_2(u_1,t_2,m,\epsilon) \exp\left(-\left(\frac{u_1}{\epsilon t_1}\right)^{k_1}\right)e^{izm}\frac{du_1}{u_1}dm,$$
with
\begin{multline*}
K_2(u_1,t_2,m,\epsilon)=\int_{L_{\tilde{\mathfrak{d}}_{p+1},\rho_2/2}}\omega_{\rho_1,\rho_2}(u_1,u_2,m,\epsilon)\exp\left(-\left(\frac{u_2}{\epsilon t_2}\right)^{k_2}\right)\frac{du_2}{u_2}\\
-\int_{L_{\tilde{\mathfrak{d}}_{p},\rho_2/2}}\omega_{\rho_1,\rho_2}(u_1,u_2,m,\epsilon)\exp\left(-\left(\frac{u_2}{\epsilon t_2}\right)^{k_2}\right)\frac{du_2}{u_2}.
\end{multline*}

A deformation of the integration path allows us to substitute the integrals in the expression of $K_2(u_1,t_2, m,\epsilon)$ as for $I_1(u_1,m,\epsilon)$ (see~(\ref{e702})) to provide similar upper estimates, and arrive at

$$|K_2(u_1,t_2,m,\epsilon)|\le |u_1|\frac{1}{(1+|m|)^{\mu}}e^{-\beta|m|}\varpi_{\rho_1,\rho_2}\frac{\rho_2}{2}|\tilde{\mathfrak{d}}_{p+1}-\tilde{\mathfrak{d}}_{p}|\exp\left(-\left(\frac{\rho_2/2}{r_{\mathcal{T}_2}}\right)^{k_2}\tilde{\nabla}_{p,p+1}\frac{1}{|\epsilon|^{k_2}}\right),$$
for some constant $\tilde{\nabla}_{p,p+1}>0$, for all $\epsilon\in \mathcal{E}_p\cap\mathcal{E}_{p+1}$, $u_1\in L_{\theta_{p,p+1},0, \rho_1/4}$ and $m\in\R$, where $r_{\mathcal{T}_2}:=\sup\{|t_2|:t_2\in\mathcal{T}_2\}$. 

As a result, 
\begin{multline*}
|I_3-I_6|\le\frac{k_1k_2}{(2\pi)^{1/2}}\left(\int_{-\infty}^{\infty}e^{-(\beta-\beta')|m|}\frac{1}{(1+|m|)^{\mu}}dm\right)\varpi_{\rho_1,\rho_2}\\
\times E_3(|\epsilon t_1|)\frac{\rho_2}{2}|\tilde{\mathfrak{d}}_{p+1}-\tilde{\mathfrak{d}}_{p}|\exp\left(-\left(\frac{\rho_2/2}{r_{\mathcal{T}_2}}\right)^{k_2}\tilde{\nabla}_{p,p+1}\frac{1}{|\epsilon|^{k_2}}\right),
\end{multline*}
with
$$E_{3}(|\epsilon t_1|)=\int_0^{\rho_1/4}\exp\left(-\left(\frac{r_1}{|\epsilon t_1|}\right)^{k_1}\right)dr_1\le\frac{\rho_1}{4}.$$

As a consequence, 
\begin{equation}\label{e822}
|I_3-I_6|\le K_{p,1}\exp\left(-\frac{K_{p,2}}{|\epsilon|^{k_2}}\right),
\end{equation}
for every $(t_1,t_2,z,\epsilon)\in \mathcal{T}_1\times \mathcal{T}_2\times H_{\beta'}\times (\mathcal{E}_p\cap \mathcal{E}_{p+1})$.

Let us provide upper estimates for $I_1$. The choice made on the domains involved in the construction of the map $J_{1,p}$ allows us to estimate

$$|I_1|\le\frac{k_1k_2}{(2\pi)^{1/2}}\left(\int_{-\infty}^{\infty}e^{-(\beta-\beta')|m|}\frac{1}{(1+|m|)^{\mu}}dm\right)\varpi_{\rho_1,\rho_2}E_{4}(|\epsilon t_1|)E_{5}(|\epsilon t_2|),$$
with
$$E_{4}(|\epsilon t_1|)=\int_{\rho_1/4}^{\rho_1/2}\exp\left(-\left(\frac{r_1}{|\epsilon t_1|}\right)^{k_1}\nabla_{p+1}\right)dr_1,$$
$$E_{5}(|\epsilon t_2|)=\int_0^{\rho_2/2}\exp\left(-\left(\frac{r_2}{|\epsilon t_2|}\right)^{k_2}\tilde{\nabla}_{p+1}\right)dr_2.$$

Observe that 
\begin{multline*}
E_4(|\epsilon t_1|)=\int_{\rho_1/4}^{\rho_1/2}\frac{k_1r_1^{k_1-1}}{|\epsilon t_1|^{k_1}}\nabla_{p+1}\exp\left(-\left(\frac{r_1}{|\epsilon t_1|}\right)^{k_1}\nabla_{p+1}\right)\left[\frac{1}{\nabla_{p+1}}\frac{|\epsilon t_1|^{k_1}}{k_1 r_1^{k_1-1}}\right]dr_1\\
\le |\epsilon t_1|^{k_1} \frac{1}{\nabla_{p+1}}\frac{1}{k_1(\rho_1/4)^{k_1-1}}\left[-\exp\left(-\left(\frac{r_1}{|\epsilon t_1|}\right)^{k_1}\nabla_{p+1}\right)\right]_{r_1=\rho_1/4}^{r_1=\rho_1/2}\\
\le \epsilon_0^{k_1}r_{\mathcal{T}_1}^{k_1} \frac{1}{\nabla_{p+1}}\frac{1}{k_1(\rho_1/4)^{k_1-1}}\exp\left(-\left(\frac{\rho_1/4}{r_{\mathcal{T}_2}}\right)^{k_1}\nabla_{p+1}\frac{1}{|\epsilon|^{k_1}}\right).
\end{multline*}

On the other hand, analogous estimates as before yield
$$E_5(|\epsilon t_2|)\le \int_0^{\rho_2/2}\exp\left(-\left(\frac{r_2}{r_{\mathcal{T}_2}|\epsilon|}\right)^{k_2}\tilde{\nabla}_{p+1}\right)dr_2\le\frac{\rho_2}{2}.$$

Therefore, there exist $K_{p,3},K_{p,4}>0$ such that 
\begin{equation}\label{e822b}
|I_1|\le K_{p,3}\exp\left(-\frac{K_{p,4}}{|\epsilon|^{k_1}}\right),
\end{equation}
for every $(t_1,t_2,z,\epsilon)\in \mathcal{T}_1\times \mathcal{T}_2\times H_{\beta'}\times (\mathcal{E}_p\cap \mathcal{E}_{p+1})$.

An analogous argument guarantees the existence of $K_{p,5},K_{p,6}>0$ such that 
\begin{equation}\label{e822c}
|I_4|\le K_{p,5}\exp\left(-\frac{K_{p,6}}{|\epsilon|^{k_1}}\right),
\end{equation}
for every $(t_1,t_2,z,\epsilon)\in \mathcal{T}_1\times \mathcal{T}_2\times H_{\beta'}\times (\mathcal{E}_p\cap \mathcal{E}_{p+1})$.

Let us provide upper estimates for $|I_2|$. Regarding the assumptions made on the domains involved in the construction of the map $J_{1,p}$, one has the existence of $\nabla_{p+1,\theta}>0$ such that
$$|I_2|\le\frac{k_1k_2}{(2\pi)^{1/2}}\left(\int_{-\infty}^{\infty}e^{-(\beta-\beta')|m|}\frac{1}{(1+|m|)^{\mu}}dm\right)\varpi_{\rho_1,\rho_2}E_{6}(|\epsilon t_1|)E_{7}(|\epsilon t_2|),$$
with
$$E_{6}(|\epsilon t_1|)=\left|\int_{\theta_{p,p+1}}^{\mathfrak{d}_{p+1}}\exp\left(-\left(\frac{\rho_1/4}{|\epsilon t_1|}\right)^{k_1}\nabla_{p+1,\theta}\right)\frac{\rho_1}{4}dh\right|,$$
$$E_{7}(|\epsilon t_2|)=\int_0^{\rho_2/2}\exp\left(-\left(\frac{r_2}{|\epsilon t_2|}\right)^{k_2}\tilde{\nabla}_{p+1}\right)dr_2.$$

Analogous estimates as for $E_4$ and $E_5$ guarantee the existence of upper bounds of $E_{6}(|\epsilon t_1|)$ and of $E_{7}(|\epsilon t_2|)$.

This entails there exist of $K_{p,7},K_{p,8}>0$ such that 
\begin{equation}\label{e822d}
|I_2|\le K_{p,7}\exp\left(-\frac{K_{p,8}}{|\epsilon|^{k_1}}\right),
\end{equation}
for every $(t_1,t_2,z,\epsilon)\in \mathcal{T}_1\times \mathcal{T}_2\times H_{\beta'}\times (\mathcal{E}_p\cap \mathcal{E}_{p+1})$.

Analogous bounds can be followed to arrive at
\begin{equation}\label{e822e}
|I_5|\le K_{p,9}\exp\left(-\frac{K_{p,10}}{|\epsilon|^{k_1}}\right),
\end{equation}
for some positive constants $K_{p,9}$ and $K_{p,10}$, valid for all $(t_1,t_2,z,\epsilon)\in \mathcal{T}_1\times \mathcal{T}_2\times H_{\beta'}\times (\mathcal{E}_p\cap \mathcal{E}_{p+1})$.

From (\ref{e822}), (\ref{e822b}), (\ref{e822c}), (\ref{e822d}) and (\ref{e822e}) we conclude the existence of $M_{p,1},M_{p,2},M_{p,3},M_{p,4}>0$ such that
$$|J_{1,p+1}(t_1,t_2,z,\epsilon)-J_{1,p}(t_1,t_2,z,\epsilon)|\le M_{p,1}\exp\left(-\frac{M_{p,2}}{|\epsilon|^{k_1}}\right)+M_{p,3}\exp\left(-\frac{M_{p,4}}{|\epsilon|^{k_2}}\right),$$
for $(t_1,t_2,z,\epsilon)\in \mathcal{T}_1\times \mathcal{T}_2\times H_{\beta'}\times (\mathcal{E}_p\cap \mathcal{E}_{p+1})$. Recall that $k_1> k_2$, which allows us to conclude the bounds provided at the statements of the third situation.

\end{proof}

At this point, we can apply Theorem (RS), see Section~\ref{secrs}. For $0<\beta'<\beta$, we write $\mathbb{E}$ for the Banach space of bounded holomorphic functions defined on $\mathcal{T}_1\times\mathcal{T}_2\times H_{\beta'}$ equipped with the norm of the supremum, i.e. $\mathbb{E}:=\mathcal{O}_b(\mathcal{T}_1\times\mathcal{T}_2\times H_{\beta'})$.

\begin{prop}\label{propteo2}
Under the previous assumptions, for every $0\le p\le \varsigma-1$, the solution (\ref{e685}) of (\ref{epral}) admits a splitting of the form
$$J_{1,p}(t_1,t_2,z,\epsilon)=a(t_1,t_2,z,\epsilon)+J_{1,1,p}(t_1,t_2,z,\epsilon)+J_{1,2,p}(t_1,t_2,z,\epsilon),$$
where $a(t_1,t_2,z,\epsilon)\in\mathbb{E}\{\epsilon\}$, and $J_{1,j,p}\in\mathcal{O}_b(\mathcal{T}_1\times\mathcal{T}_2\times H_{\beta'}\times \mathcal{E}_p)$, for all $0\le p\le \varsigma-1$ and $j=1,2$. Moreover, there exist two formal power series in $\epsilon$ with coefficients in $\mathbb{E}$, say $\hat{J}_{1,j}\in\mathbb{E}[[\epsilon]]$, for $j=1,2$ which satisfy that $J_{1,j,p}$ admits $\hat{J}_{1,j}$ as its common Gevrey asymptotic expansion of order $1/k_j$ with respect to $\epsilon$ on $\mathcal{E}_p$, for all $0\le p\le \varsigma-1$ (see Section~\ref{secrs} for its precise definition). 
\end{prop}
\begin{proof}
Let us split the set $\{0,\ldots,\varsigma-1\}$ into the set $I_1$ of indices such that Case 1 or Case 3 of Proposition~\ref{prop689} hold, and $I_2=\{0,\ldots,\varsigma-1\}\setminus I_1$ (i.e. the set of indices for which Case 2 of Proposition~\ref{prop689} holds). Multilevel Ramis-Sibuya Theorem (RS) can be applied to the functions $G_p:\mathcal{E}_p\to\mathcal{O}_b(\mathcal{T}_1\times\mathcal{T}_2\times H_{\beta'})$ defined by 
$$G_p(\epsilon):=J_{1,p}(t_1,t_2,z,\epsilon),\quad \epsilon\in\mathcal{E}_p,$$
for $0\le p\le \varsigma-1$, and where $\mathbb{E}$ stands for the Banach space of holomorphic and bounded functions defined on $\mathcal{T}_1\times\mathcal{T}_2\times H_{\beta'}$ equipped with the norm of the supremum, for some fixed $0<\beta'<\beta$. This is a consequence of the different exponential decays in the perturbation parameter, uniform on the rest of variables, showed in Proposition~\ref{prop689}.

\end{proof}

The particularization of the previous result to the first index in the good covering allows us to conclude.

\begin{corol}
In the situation of Proposition~\ref{propteo2}, the analytic map $J_{1}(t_1,t_2,z,\epsilon)$, defined in (\ref{e792}), defined in $\mathcal{T}_1\times\mathcal{T}_2\times H_{\beta'}\times\mathcal{E}$ admits a splitting of the form
$$J_1(t_1,t_2,z,\epsilon)=a(t_1,t_2,z,\epsilon)+J_{1,1,0}(t_1,t_2,z,\epsilon)+J_{1,2,0}(t_1,t_2,z,\epsilon),$$
where $a(t_1,t_2,z,\epsilon)\in\mathbb{E}\{\epsilon\}$, and $J_{1,j,0}\in\mathcal{O}_b(\mathcal{T}_1\times\mathcal{T}_2\times H_{\beta'}\times \mathcal{E})$ for $j=1,2$. Moreover, the formal power series $\hat{J}_{1,j}\in\mathbb{E}[[\epsilon]]$, for $j=1,2$ satisfy that $J_{1,j,0}$ admits $\hat{J}_{1,j}$ as its common Gevrey asymptotic expansion of order $1/k_j$ with respect to $\epsilon$ on $\mathcal{E}$. 
\end{corol}

\subsection{Gevrey bounds for $J_2$ and $J_3$}\label{secj23}

We recall the next lemma from~\cite{malek20}, which will be crucial in the next two propositions.

\begin{lemma}[Lemma 14. 1),~\cite{malek20}]\label{lema1131}
Let $k'\ge1$ be an integer number, and let $M>0$ a real number. There exists $C_{k'}>0$ only depending on $k'$ such that the next inequality 
$$\left(\frac{1}{r}\right)^{N}\exp\left(-\frac{M}{r^{k'}}\right)\le C_{k'}A_{k'}^{N}\left(\frac{N}{k'}\right)^{1/2}\Gamma\left(\frac{N}{k'}\right)$$
holds for all integer $N\ge1$ and any real number $r>0$, and where $A_{k'}=(1/M)^{1/k'}$.
\end{lemma}

The next result provides bounds for $J_{2}$.

\begin{prop}\label{prop8}
There exist $C_{J_2},K_{J_2}>0$ such that
\begin{equation}\label{e1141}
|J_{2}(t_1,t_2,z,\epsilon)|\le C_{J_2}K_{J_2}^N\left(\frac{N}{k_2}\right)^{1/2}\Gamma\left(\frac{N}{k_2}\right)|\epsilon|^N,
\end{equation}
for all positive integer $N$, all $t_1\in\mathcal{T}_1$, $t_2\in\mathcal{T}_2$, $z\in H_{\beta'}$ and $\epsilon\in\mathcal{E}$.
\end{prop}
\begin{proof}
In view of the estimates for $\omega_{S_{d_1},S_{d_2}}$ determined in (\ref{e401b}) we deduce that
\begin{multline}
|J_2(t_1,t_2,z,\epsilon)|\le \frac{k_1k_2}{(2\pi)^{1/2}}\left(\int_{-\infty}^{\infty}(1+|m|)^{-\mu}e^{-\beta|m|}e^{\beta'|m|}dm\right)\varpi_{S_{d_1},S_{d_2}}\\
\times \int_{0}^{\rho_1/2}\exp\left(\nu_1 r_1^{k_1}\right)\exp\left(-\frac{r_1^{k_1}}{|\epsilon t_1|^{k_1}}\Delta_1\right)dr_1\times \int_{\rho_2/2}^{\infty}\exp\left(\nu_2 r_2^{k_2}\right)\exp\left(-\frac{r_2^{k_2}}{|\epsilon t_2|^{k_2}}\Delta_2\right)dr_2,\label{e1152a}
\end{multline}
where $\Delta_1,\Delta_2>0$ are the constants in Theorem~\ref{teo1}, for all $t_1\in\mathcal{T}_1$, $t_2\in\mathcal{T}_2$, $z\in H_{\beta'}$ and $\epsilon\in\mathcal{E}$.

Observe that
\begin{equation}\label{e1153}
\int_{0}^{\rho_1/2}\exp\left(\nu_1 r_1^{k_1}\right)\exp\left(-\frac{r_1^{k_1}}{|\epsilon t_1|^{k_1}}\Delta_1\right)dr_1\le \exp\left(\nu_1 \left(\frac{\rho_1}{2}\right)^{k_1}\right).
\end{equation}
In addition, one has that 
\begin{multline}
\int_{\rho_2/2}^{\infty}\exp\left(\nu_2 r_2^{k_2}\right)\exp\left(-\frac{r_2^{k_2}}{|\epsilon t_2|^{k_2}}\Delta_2\right)dr_2\\
= \int_{\rho_2/2}^{\infty}\exp\left(\nu_2 r_2^{k_2}\right)\exp\left(-\frac{r_2^{k_2}}{2|\epsilon t_2|^{k_2}}\Delta_2\right)\exp\left(-\frac{r_2^{k_2}}{2|\epsilon t_2|^{k_2}}\Delta_2\right)dr_2\\
\le \exp\left(-\frac{(\rho_2/2)^{k_2}}{2|\epsilon|^{k_2}r_{\mathcal{T}_2}^{k_2}}\Delta_2\right)\left(\int_{\rho_2/2}^{\infty}\exp\left(\nu_2 r_2^{k_2}\right)\exp\left(-\frac{r_2^{k_2}}{2\epsilon_0^{k_2} r_{\mathcal{T}_2}^{k_2}}\Delta_2\right)dr_2\right),\label{e1158}
\end{multline}
the last integral appearing above being convergent provided that $r_{\mathcal{T}_2}>0$ is chosen to be small enough. We apply Lemma~\ref{lema1131} to conclude that 
\begin{equation}\label{e1162}
\exp\left(-\frac{(\rho_2/2)^{k_2}}{2|\epsilon|^{k_2}r_{\mathcal{T}_2}^{k_2}}\Delta_2\right)\le C_{k_2}A_{k_2}^N\left(\frac{N}{k_2}\right)^{1/2}\Gamma\left(\frac{N}{k_2}\right)|\epsilon|^N,
\end{equation}
for all positive integer $N\ge 1$, $\epsilon\in\mathcal{E}$, where $C_{k_2}$ is a constant depending on $k_2$, and 
$$A_{k_2}=\left(\frac{2r_{\mathcal{T}_2}^{k_2}}{\Delta_2(\rho_2/2)^{k_2}}\right)^{\frac{1}{k_2}}.$$
The estimates (\ref{e1141}) are deduced from the inequalities (\ref{e1152a}), (\ref{e1153}), (\ref{e1158}), (\ref{e1162}).

\end{proof}



Regarding the estimates for $J_3$, one arrives at the following result.

\begin{prop}\label{prop9}
There exist $C_{J_3},K_{J_3}>0$ such that
\begin{equation}\label{e1141b}
|J_{3}(t_1,t_2,z,\epsilon)|\le C_{J_3}K_{J_3}^N\left(\frac{N}{k_1}\right)^{1/2}\Gamma\left(\frac{N}{k_1}\right)|\epsilon|^N,
\end{equation}
for all positive integer $N$, all $t_1\in\mathcal{T}_1$, $t_2\in\mathcal{T}_2$, $z\in H_{\beta'}$ and $\epsilon\in\mathcal{E}$.
\end{prop}
\begin{proof}
In view of the estimates for $\omega_{S_{d_1},S_{d_2}}$ determined in (\ref{e401b}) we deduce that
\begin{multline}
|J_3(t_1,t_2,z,\epsilon)|\le \frac{k_1k_2}{(2\pi)^{1/2}}\left(\int_{-\infty}^{\infty}(1+|m|)^{-\mu}e^{-\beta|m|}e^{\beta'|m|}dm\right)\varpi_{S_{d_1},S_{d_2}}\\
\times \int_{\rho_1/2}^{\infty}\exp\left(\nu_1 r_1^{k_1}\right)\exp\left(-\frac{r_1^{k_1}}{|\epsilon t_1|^{k_1}}\Delta_1\right)dr_1\times \int_{0}^{\infty}\exp\left(\nu_2 r_2^{k_2}\right)\exp\left(-\frac{r_2^{k_2}}{|\epsilon t_2|^{k_2}}\Delta_2\right)dr_2,\label{e1152}
\end{multline}
where $\Delta_1,\Delta_2>0$ are the constants in Theorem~\ref{teo1}, for all $t_1\in\mathcal{T}_1$, $t_2\in\mathcal{T}_2$, $z\in H_{\beta'}$ and $\epsilon\in\mathcal{E}$.

A similar reasoning made in (\ref{e1158}) can be followed to arrive at
\begin{multline}\int_{\rho_1/2}^{\infty}\exp\left(\nu_1 r_1^{k_1}\right)\exp\left(-\frac{r_1^{k_1}}{|\epsilon t_1|^{k_1}}\Delta_1\right)dr_1\le \left(\int_{\rho_1/2}^{\infty}\exp\left(\nu_1 r_1^{k_1}\right)\exp\left(-\frac{r_1^{k_1}}{2\epsilon_0^{k_1} r_{\mathcal{T}_1}^{k_1}}\Delta_1\right)dr_1\right)\\
\times  C_{k_1}A_{k_1}^N\left(\frac{N}{k_1}\right)^{1/2}\Gamma\left(\frac{N}{k_1}\right)|\epsilon|^N,\label{e1203}
\end{multline}
for all positive integer $N\ge1$, $\epsilon\in\mathcal{E}$, where $C_{k_1}$ is a constant depending on $k_1$, and 
$$A_{k_1}=\left(\frac{2r_{\mathcal{T}_1}^{K_1}}{\Delta_1(\rho_1/2)^{k_1})}\right)^{\frac{1}{k_1}},$$
where the last integral appearing in (\ref{e1203}) is finite provided that $r_{\mathcal{T}_1}$ is chosen small enough.

Regarding the last integral in (\ref{e1152}) we observe that
\begin{equation}\label{e1210}
\int_{0}^{\infty}\exp\left(\nu_2 r_2^{k_2}\right)\exp\left(-\frac{r_2^{k_2}}{|\epsilon t_2|^{k_2}}\Delta_2\right)dr_2\le \int_{0}^{\infty}\exp\left(\nu_2 r_2^{k_2}\right)\exp\left(-\frac{r_2^{k_2}}{(\epsilon_0 r_{\mathcal{T}_{2}})^{k_2}}\Delta_2\right)dr_2
\end{equation}
where the last integral is finite provided that $r_{\mathcal{T}_2}$ is small enough.

The estimates (\ref{e1141b}) are deduced from the inequalities (\ref{e1152}), (\ref{e1203}) and (\ref{e1210}).

\end{proof}

\subsection{Main asymptotic result}

As before, for $0<\beta'<\beta$, we write $\mathbb{E}$ for the Banach space of bounded holomorphic functions defined on $\mathcal{T}_1\times\mathcal{T}_2\times H_{\beta'}$ equipped with the norm of the supremum, i.e. $\mathbb{E}:=\mathcal{O}_b(\mathcal{T}_1\times\mathcal{T}_2\times H_{\beta'})$.

\begin{theo}\label{teo2}
Under the previous assumptions, the solution (\ref{e776}) of (\ref{epral}) admits a splitting of the form
$$u_{d_1,d_2}(t_1,t_2,z,\epsilon)=b(t_1,t_2,z,\epsilon)+u_{d_1,d_2,1}(t_1,t_2,z,\epsilon)+u_{d_1,d_2,2}(t_1,t_2,z,\epsilon),$$
where $b(t_1,t_2,z,\epsilon)\in\mathbb{E}\{\epsilon\}$, and $u_{d_1,d_2,j}\in\mathcal{O}_b(\mathcal{T}_1\times\mathcal{T}_2\times H_{\beta'}\times \mathcal{E})$, for $j=1,2$. Moreover, there exist two formal power series in $\epsilon$ with coefficients in $\mathbb{E}$, say 
$$\hat{u}_{j}(t_1,t_2,z,\epsilon)=\sum_{k\ge0}H_{k}^{j}(t_1,t_2,z)\epsilon^k\in\mathbb{E}[[\epsilon]],$$ 
for $j=1,2$ which satisfy that $u_{d_1,d_2,j}$ admits $\hat{u}_{j}$ as its Gevrey asymptotic expansion of order $1/k_j$ with respect to $\epsilon$ on $\mathcal{E}$, which means that for all $\mathcal{W}\prec \mathcal{E}$, there exist $C,A>0$ with 
$$\sup_{(t_1,t_2,z)\in\mathcal{T}_1\times\mathcal{T}_2\times H_{\beta'}}\left|u_{d_1,d_2,j}(t_1,t_2,z,\epsilon)-\sum_{p=0}^{N-1}H_p^j(t_1,t_2,z)\epsilon^p \right|\le C A^{N}\Gamma\left(1+\frac{N}{k_j}\right)|\epsilon|^{N}\quad ,\epsilon\in \mathcal{W},$$
valid for all $N\ge1$.
\end{theo}
\begin{proof}

In view of the splitting in (\ref{e792}), the Gevrey expansions for $J_1$ determined in Proposition~\ref{propteo2}, together with the Gevrey bounds attained in Proposition~\ref{prop8} and Proposition~\ref{prop9}, we set:
\begin{multline*} 
b(t_1,t_2,z,\epsilon)=a(t_1,t_2,z,\epsilon)\hbox{ obtained in Proposition~\ref{propteo2}},\\
u_{d_1,d_2,1}(t_1,t_2,z,\epsilon)=J_{1,1,0}(t_1,t_2,z,\epsilon)+J_3(t_1,t_2,z,\epsilon),\\ 
u_{d_1,d_2,2}(t_1,t_2,z,\epsilon)=J_{1,2,0}(t_1,t_2,z,\epsilon)+J_2(t_1,t_2,z,\epsilon), 
\end{multline*}
and
$$\hat{u}_j(t_1,t_2,z,\epsilon)=\hat{J}_{1,j},\quad j=1,2.$$
The result of Theorem~\ref{teo2} is a straight consequence of Proposition~\ref{propteo2}, Proposition~\ref{prop8} and Proposition~\ref{prop9}.
\end{proof}

\section{Appendix}\label{secapendice}
\subsection{Inverse Fourier transform. Related function spaces}\label{secapendice1}

Let $\beta>0$ and $\mu>1$. For any $f\in\mathcal{C}(\R)$ such that there exists $C>0$ with
$$|f(m)|\le C \frac{1}{(1+|m|)^{\mu}}\exp(-\beta|m|),\qquad m\in\R,$$
one can define 
$$\mathcal{F}^{-1}(f)(x)=\frac{1}{\sqrt{2\pi}}\int_{-\infty}^{\infty}f(m)\exp(ixm)dm,\quad x\in\R.$$
It holds that $\mathcal{F}^{-1}(f)\in\mathcal{O}(H_{\beta}),$ where $H_{\beta}$ stands for the horizontal strip
$$H_{\beta}=\{z\in\C:|\hbox{Im}(z)|<\beta\}.$$
The set of functions satisfying the previous bounds determine the following Banach space.

\begin{defin}\label{defi937}
Let $\beta,\mu>0$. The set of continuous functions $h:\R\to\C$ such that
$$\left\|h(m)\right\|_{(\beta,\mu)}:=\sup_{m\in\R}(1+|m|)^{\mu}\exp(\beta|m|)|h(m)|$$ is a finite quantity determines a Banach space, denoted by $(E_{(\beta,\mu)},\left\|\cdot\right\|_{(\beta,\mu)})$.
\end{defin}

The proof of the following result can be found in detail in~\cite{lama}, Proposition 7.

\begin{prop}\label{prop1anexo}
Let $\beta>0$ and $\mu>1$. The following statements hold for every $h\in E_{(\beta,\mu)}$:
\begin{itemize}
\item Let $\varphi(m):= imh(m)$, for $m\in\R$. Then, $\varphi\in E_{(\beta,\mu-1)}$ and $\partial_z\mathcal{F}^{-1}(h)(z)=\mathcal{F}^{-1}(\varphi)(z)$, for all $z\in H_{\beta}$.
\item Let $g\in E_{(\beta,\mu)}$. Define the function $\psi$ by $\psi(m)=\frac{1}{(2\pi)^{1/2}}(h\star g)(m)$ for every $m\in\R$, where $\star$ denotes the convolution product of $h$ and $g$:
$$(h\star g)(m)=\int_{-\infty}^{\infty}h(m-m_1)g(m_1)dm_1,\quad m\in\R.$$ Then, it holds that $\psi\in E_{(\beta,\mu)}$, together with
$$\mathcal{F}^{-1}(h)(z)\mathcal{F}^{-1}(g)(z)=\mathcal{F}^{-1}(\psi)(z),$$
for every $z\in H_{\beta}$.
\end{itemize}
\end{prop}

\subsection{Two-level version of Ramis-Sibuya theorem}\label{secrs}

In this section, we recall the main facts about Gevrey asymptotic expansions concerning holomorphic functions.

\begin{defin}
Let $(\mathbb{E},\left\|\cdot\right\|_{\mathbb{E}})$ be a complex Banach space, and let $\mathcal{E}$ be a bounded sector with vertex at the origin. We also fix $k>0$. We say that $f(\epsilon)\in\mathcal{O}(\mathcal{E},\mathbb{E})$ admits $\hat{f}(\epsilon)=\sum_{p\ge0}a_p\epsilon^p\in\mathbb{E}[[\epsilon]]$ as its Gevrey asymptotic expansion of order $1/k$ (at the origin) if for all $T\prec \mathcal{E}$, there exist $C,A>0$ with 
$$\left\|f(\epsilon)-\sum_{p=0}^{N-1}a_p\epsilon^p \right\|_{\mathbb{E}}\le C A^{N}\Gamma\left(1+\frac{N}{k}\right)|\epsilon|^{N}\quad ,\epsilon\in T,$$
valid for all $N\ge1$.
\end{defin}

The classical Ramis-Sibuya Theorem can be found in detail in Theorem XI-2-3,~\cite{hssi} whereas the proof of the following generalization can be found in detail in Theorem (RS), page 63-65 in~\cite{lama1}.

\begin{theo}[RS]
Let $0<k_2<k_1$. We also fix a complex Banach space $(\mathbb{E},\left\|\cdot\right\|_{\mathbb{E}})$, and a good covering in $\C^{\star}$, $(\mathcal{E}_{p})_{0\le p\le \varsigma-1}$, for some integer $\varsigma\ge 2$ (see Definition~\ref{defi662} in Section~\ref{secpar}). Given $0\le p\le \varsigma-1$, we assume $G_p\in\mathcal{O}_b(\mathcal{E}_p,\mathbb{E})$, and define $\Delta_p(\epsilon)=G_{p+1}(\epsilon)-G_p(\epsilon)$, for $\epsilon\in \mathcal{E}_p\cap\mathcal{E}_{p+1}$, with the convention that $\mathcal{E}_{\varsigma}:=\mathcal{E}_0$ and $G_\varsigma:=G_0$. Assume moreover the existence of $I_1, I_2\subseteq\{0,\ldots,\varsigma-1\}$, such that $I_1,I_2\neq\emptyset$, and $I_1\cup I_2=\{0,\ldots,\varsigma-1\}$, with $I_1\cap I_2=\emptyset$, satisfying the following properties:
\begin{itemize}
\item for every $p\in I_1$, there exists $K_p,M_p>0$ such that 
$$\left\|\Delta_p(\epsilon)\right\|_{\mathbb{E}}\le K_p\exp\left(-\frac{M_p}{|\epsilon|^{k_1}}\right),\quad \epsilon\in \mathcal{E}_p\cap\mathcal{E}_{p+1},$$
\item for every $p\in I_2$, there exists $\tilde{K}_p,\tilde{M}_p>0$ such that 
$$\left\|\Delta_p(\epsilon)\right\|_{\mathbb{E}}\le \tilde{K}_p\exp\left(-\frac{\tilde{M}_p}{|\epsilon|^{k_2}}\right),\quad \epsilon\in \mathcal{E}_p\cap\mathcal{E}_{p+1}.$$
\end{itemize}
Then, there exist $a(\epsilon)\in\mathbb{E}\{\epsilon\}$, two formal power series $\hat{G}^1,\hat{G}^2\in\mathbb{E}[[\epsilon]]$, and for all $0\le p\le \varsigma-1$ two functions $G^1_p,G^2_p\in\mathcal{O}_b(\mathcal{E}_p,\mathbb{E})$ such that:
\begin{itemize}
\item[1)] For all $0\le p\le \varsigma-1$, the function $G_p$ admits the decomposition
$$G_p(\epsilon)=a_p(\epsilon)+G^1_p(\epsilon)+G^2_p(\epsilon),\quad \epsilon\in\mathcal{E}_p.$$
\item[2)] For $j=1,2$, and all $0\le p\le \varsigma-1$, the function $G^j_p$ admits $\hat{G}^j$ as its Gevrey asymptotic expansion of order $1/k_j$ on $\mathcal{E}_p$.
\end{itemize}

\end{theo}

\subsection{Auxiliary Banach spaces of analytic functions}\label{secauxban}

In this section, we state the definition of some auxiliary Banach spaces of functions which allow us to provide some important properties of analytic continuation of the solution to the auxiliary problem (\ref{e309}) in the Borel-Fourier space.

\subsubsection{First auxiliary Banach space}\label{sec452}

Let $\beta>0$ and $\mu>1$. Let us also fix $\rho_1,\rho_2>0$. 
 In the whole section, we fix positive integers $k_1,k_2$.

\begin{defin}\label{def963}
Let us consider the set of continuous maps $(\tau_1,\tau_2,m)\mapsto h(\tau_1,\tau_2,m)$ defined on $D(0,\rho_1)\times D(0,\rho_2)\times \R$, holomorphic with respect to its first two variables on $D(0,\rho_1)\times D(0,\rho_2)$ and such that there exists $C>0$ (which depends on $\beta,\mu,\rho_1,\rho_2$) with
$$|h(\tau_1,\tau_2,m)|\le C\frac{1}{(1+|m|)^{\mu}}e^{-\beta|m|}|\tau_1\tau_2|,$$
for every $(\tau_1,\tau_2,m)\in D(0,\rho_1)\times D(0,\rho_2)\times\R$. Such set is denoted by $B_{(\beta,\mu,\rho_1,\rho_2)}$. Given $h$ as before, we denote the minimum of such constant $C$ by $\left\| h(\tau_1,\tau_2,m)\right\|_{(\beta,\mu,\rho_1,\rho_2)}$. 
The pair $(B_{(\beta,\mu,\rho_1,\rho_2)},\left\| \cdot\right\|_{(\beta,\mu,\rho_1,\rho_2)})$ is a complex Banach space.
\end{defin}

We state some properties associated to the previous Banach space. The first one is a direct consequence of its definition.

\begin{prop}\label{prop477}
Let $(\tau_1,\tau_2,m)\mapsto b(\tau_1,\tau_2,m)$ be a continuous function defined on $D(0,\rho_1)\times D(0,\rho_2)\times \R$, holomorphic with respect to its first two variables on $D(0,\rho_1)\times D(0,\rho_2)$. Assume that 
$$C_b:=\sup_{(\tau_1,\tau_2,m)\in D(0,\rho_1)\times D(0,\rho_2)\times \R}|b(\tau_1,\tau_2,m)|$$
is finite. Then, for every $h\in B_{(\beta,\mu,\rho_1,\rho_2)}$, the function  $(\tau_1,\tau_2,m)\mapsto b(\tau_1,\tau_2,m)h(\tau_1,\tau_2,m)$ belongs to $B_{(\beta,\mu,\rho_1,\rho_2)}$ and it holds that
$$\left\|b(\tau_1,\tau_2,m)h(\tau_1,\tau_2,m)\right\|_{(\beta,\mu,\rho_1,\rho_2)}\le C_b 
\left\|h(\tau_1,\tau_2,m)\right\|_{(\beta,\mu,\rho_1,\rho_2)}.$$
\end{prop}

\begin{prop}\label{prop441}
Let $a(\tau_1,\tau_2,m)$ be a continuous function defined on $D(0,\rho_1)\times D(0,\rho_2)\times \R$, holomorphic with respect to its two first variables on $D(0,\rho_1)\times D(0,\rho_2)$. We assume this function satisfies there exists $\gamma_1\ge0$ with
$$|a(\tau_1,\tau_2,m)|\le \frac{C_1}{(1+|m|)^{\gamma_1}},\quad (\tau_1,\tau_2,m)\in D(0,\rho_1)\times D(0,\rho_2)\times \R.$$
We choose a mapping $m\mapsto h(m,\epsilon)$ which belongs to $E_{(\beta,\mu)}$, such that for every $m\in\R$, the function $D(0,\epsilon_0)\ni \epsilon\mapsto h(m,\epsilon)$ is holomorphic on $D(0,\epsilon_0)$ and there exists $K>0$ with
$$\sup_{\epsilon\in D(0,\epsilon_0)}\left\|h(m,\epsilon)\right\|_{(\beta,\mu)}\le K.$$
We consider a polynomial $P(X)\in\C[X]$. We assume that 
$$\gamma_1\ge\hbox{deg}(P),\qquad \mu> \hbox{deg}(P)+1.$$

Let us also fix $\sigma_j>-1$ for $j=1,\ldots,6$ with
$$k_1\sigma_1+\sigma_3+\sigma_5+\frac{1}{k_1}\ge 0,\quad k_2\sigma_2+\sigma_4+\sigma_6+\frac{1}{k_2}\ge0.$$

Then, for every $f\in B_{(\beta,\mu,\rho_1,\rho_2)}$ the function
\begin{multline*}\mathcal{B}_1(f):=a(\tau_1,\tau_2,m)\int_{-\infty}^{\infty}h(m-m_1,\epsilon)\tau_1^{\sigma_1k_1}\tau_2^{\sigma_2k_2}\int_{0}^{\tau_1^{k_1}}\int_{0}^{\tau_2^{k_2}}(\tau_1^{k_1}-s_1)^{\sigma_3}(\tau_2^{k_2}-s_2)^{\sigma_4}s_1^{\sigma_5}s_2^{\sigma_6}P(im_1)\\
\times f(s_1^{1/k_1},s_2^{1/k_2},m_1)ds_2ds_1dm_1
\end{multline*}
belongs to $B_{(\beta,\mu,\rho_1,\rho_2)}$. In addition to this, there exists $\tilde{C}_1>0$ with
$$\left\|\mathcal{B}_1(f)\right\|_{(\beta,\mu,\rho_1,\rho_2)}\le K\tilde{C}_1 \left\|f\right\|_{(\beta,\mu,\rho_1,\rho_2)}.$$
\end{prop}
\begin{proof}
We observe there exists $C_P>0$ such that $|P(im)|\le C_P(1+|m|)^{\hbox{deg}(P)}$, for all $m\in\R$. On the other hand, for every $(\tau_1,\tau_2,m)\in D(0,\rho_1)\times D(0,\rho_2)\times \R$ one has that
\begin{multline*}
\left|a(\tau_1,\tau_2,m)\int_{-\infty}^{\infty}h(m-m_1,\epsilon)\tau_1^{\sigma_1k_1}\tau_2^{\sigma_2k_2}\int_{0}^{\tau_1^{k_1}}\int_{0}^{\tau_2^{k_2}}(\tau_1^{k_1}-s_1)^{\sigma_3}(\tau_2^{k_2}-s_2)^{\sigma_4}s_1^{\sigma_5}s_2^{\sigma_6}P(im_1)\right.\\
\left.\times f(s_1^{1/k_1},s_2^{1/k_2},m_1)ds_2ds_1dm_1\right|\\
\le \frac{KC_1C_P}{(1+|m|)^{\gamma_1}}\int_{-\infty}^{\infty}\frac{1}{(1+|m-m_1|)^{\mu}}e^{-\beta|m-m_1|}|\tau_1|^{\sigma_1k_1}|\tau_2|^{\sigma_2k_2}\hfill\\
\times\int_{0}^{|\tau_1|^{k_1}}\int_{0}^{|\tau_2|^{k_2}}(|\tau_1|^{k_1}-r_1)^{\sigma_3}(|\tau_2|^{k_2}-r_2)^{\sigma_4}r_1^{\sigma_5}r_2^{\sigma_6}(1+|m_1|)^{\hbox{deg}(P)}\\
\times\left(|f(r_1^{1/k_1}e^{i\hbox{arg}(\tau_1)},r_2^{1/k_2}e^{i\hbox{arg}(\tau_2)},m_1)| (1+|m_1|)^{\mu}e^{\beta|m_1|}\frac{1}{r_1^{1/k_1}}\frac{1}{r_2^{1/k_2}}
\right)\\
\hfill\times \frac{1}{(1+|m_1|)^{\mu}}e^{-\beta|m_1|}r_1^{1/k_1}r_2^{1/k_2}dr_2dr_1dm_1.
\end{multline*}
Taking into account that $e^{-\beta|m-m_1|}e^{-\beta|m_1|}\le e^{-\beta|m|}$, the previous expression is upper bounded by
\begin{multline*}
KC_1C_P\left\|f\right\|_{(\beta,\mu,\rho_1,\rho_2)}e^{-\beta|m|}\frac{1}{(1+|m|)^{\gamma_1}}\int_{-\infty}^{\infty}\frac{1}{(1+|m-m_1|)^{\mu}}\frac{1}{(1+|m_1|)^{\mu-\hbox{deg}(P)}}dm_1 \\
\times|\tau_1|^{\sigma_1k_1}|\tau_2|^{\sigma_2k_2}\int_{0}^{|\tau_1|^{k_1}}\int_{0}^{|\tau_2|^{k_2}}(|\tau_1|^{k_1}-r_1)^{\sigma_3}(|\tau_2|^{k_2}-r_2)^{\sigma_4}r_1^{\sigma_5+1/k_1}r_2^{\sigma_6+1/k_2}
dr_2dr_1.
\end{multline*}
We observe that for $x>0$ and the change of variable $r=xs$ one has that
$$\int_0^{x}(x-r)^{\alpha_1}r^{\alpha_2}dr=x^{\alpha_1+\alpha_2+1}\int_0^{1}(1-s)^{\sigma_4}s^{\sigma_6+\frac{1}{k_2}}ds,$$
which allows us to conclude that
\begin{multline*}\left\|\mathcal{B}_1(f)\right\|_{(\beta,\mu,\rho_1,\rho_2)}\le KC_1C_P\left(\int_0^1(1-s)^{\sigma_3}s^{\sigma_5+1/k_1}ds\right)\left(\int_0^1(1-s)^{\sigma_4}s^{\sigma_6+1/k_2}ds\right)\\
\times\left(\sup_{m\in\R}\int_{-\infty}^{\infty}\frac{(1+|m|)^{\mu-\gamma_1}}{(1+|m-m_1|)^{\mu}(1+|m_1|)^{\mu-\hbox{deg}(P)}}dm_1\right)\rho_1^{\sigma_1k_1+\sigma_3+\sigma_5+\frac{1}{k_1}}\rho_2^{\sigma_2k_2+\sigma_4+\sigma_6+\frac{1}{k_2}}\\
\hfill\times \left\|f\right\|_{(\beta,\mu,\rho_1,\rho_2)}.
\end{multline*}

The proof concludes by taking into account the bounds on the integral on $m_1$, which are a direct consequence of Lemma 2.2,~\cite{cota2}, and the statements on the valid values of the parameters.
\end{proof}

\begin{prop}\label{prop487}
Let $a(\tau_1,\tau_2,m)$ be a continuous function defined on $D(0,\rho_1)\times D(0,\rho_2)\times \R$, holomorphic with respect to its two first variables on $D(0,\rho_1)\times D(0,\rho_2)$. We assume this function satisfies there exist $\gamma_1\ge0$ and $C_1>0$ with
$$|a(\tau_1,\tau_2,m)|\le \frac{C_1}{(1+|m|)^{\gamma_1}},\quad (\tau_1,\tau_2,m)\in D(0,\rho_1)\times D(0,\rho_2)\times \R.$$
Let $P_1(\epsilon,X),P_2(\epsilon,X)\in\mathcal{O}_b(D(0,\epsilon_0))[X]$ be polynomials with coefficients in the set of holomorphic and bounded functions on $D(0,\epsilon_0)$. We assume that
$$\gamma_1\ge\max\{\hbox{deg}(P_1),\hbox{deg}(P_2)\}.$$
Let us choose $\mu$ such that
$$\mu>\max\{\hbox{deg}(P_1),\hbox{deg}(P_2)\}+1.$$
For every $f,g\in B_{(\beta,\mu,\rho_1,\rho_2)}$, the function
\begin{multline*}
\mathcal{B}_2(f,g):=a(\tau_1,\tau_2,m)\tau_1^{k_1}\tau_2^{k_2}\int_{0}^{\tau_1^{k_1}}\int_{0}^{\tau_2^{k_2}}\int_{-\infty}^{\infty}P_1(\epsilon,i(m-m_1))f((\tau_1^{k_1}-s_1)^{1/k_1},(\tau_2^{k_2}-s_2)^{1/k_2},m-m_1)\\
\hfill\times P_2(\epsilon,im_1)g(s_1^{1/k_1},s_2^{1/k_2},m_1)\frac{1}{\tau_1^{k_1}-s_1}\frac{1}{s_1}\frac{1}{\tau_2^{k_2}-s_2}\frac{1}{s_2}dm_1ds_2ds_1
\end{multline*}
belongs to $B_{(\beta,\mu,\rho_1,\rho_2)}$. In addition to this, there exists $\tilde{C}_2>0$ such that
$$\left\|\mathcal{B}_2(f,g)\right\|_{(\beta,\mu,\rho_1,\rho_2)}\le \tilde{C}_2\left\|f\right\|_{(\beta,\mu,\rho_1,\rho_2)}\left\|g\right\|_{(\beta,\mu,\rho_1,\rho_2)}.$$
\end{prop}
\begin{proof}
Taking into account the definition of the Banach space $B_{(\beta,\mu,\rho_1,\rho_2)}$, we observe that for every $(\tau_1,\tau_2,m)\in D(0,\rho_1)\times D(0,\rho_2)\times \R$ one can follow an analogous line of arguments as in the previous result to arrive at
\begin{multline*}
\left|a(\tau_1,\tau_2,m)\tau_1^{k_1}\tau_2^{k_2}\int_{0}^{\tau_1^{k_1}}\int_{0}^{\tau_2^{k_2}}\int_{-\infty}^{\infty}P_1(\epsilon,i(m-m_1))f((\tau_1^{k_1}-s_1)^{1/k_1},(\tau_2^{k_2}-s_2)^{1/k_2},m-m_1)\right.\\
\left.\times P_2(\epsilon,im_1)g(s_1^{1/k_1},s_2^{1/k_2},m_1)\frac{1}{\tau_1^{k_1}-s_1}\frac{1}{s_1}\frac{1}{\tau_2^{k_2}-s_2}\frac{1}{s_2}dm_1ds_2ds_1\right|\\
\le \frac{C_1}{(1+|m|)^{\gamma_1}}|\tau_1|^{k_1}|\tau_2|^{k_2}\int_0^{|\tau_1|^{k_1}}\int_0^{|\tau_2|^{k_2}}\int_{-\infty}^{\infty}C_{P_1}(1+|m-m_1|)^{\hbox{deg}(P_1)}\frac{1}{(1+|m-m_1|)^{\mu}}e^{-\beta|m-m_1|}\\
\times|\tau_1^{k_1}-r_1e^{ik_1\hbox{arg}(\tau_1)}|^{1/k_1}|\tau_2^{k_2}-r_2e^{ik_2\hbox{arg}(\tau_2)}|^{1/k_2}C_{P_2}(1+|m_1|)^{\hbox{deg}(P_2)}\frac{1}{(1+|m_1|)^{\mu}}e^{-\beta|m_1|}r_1^{1/k_1}r_2^{1/k_2}\\
\times\frac{1}{|\tau_1^{k_1}-r_1e^{ik_1\hbox{arg}(\tau_1)}|}\frac{1}{r_1}\frac{1}{|\tau_2^{k_2}-r_2e^{ik_2\hbox{arg}(\tau_2)}|}\frac{1}{r_2}dm_1dr_2dr_1
\left\|f\right\|_{(\beta,\mu,\rho_1,\rho_2)}\left\|g\right\|_{(\beta,\mu,\rho_1,\rho_2)},
\end{multline*}
i.e.
\begin{multline*}
\left\|\mathcal{B}_2(f,g)\right\|_{(\beta,\mu.\rho_1,\rho_2)}\le C_1C_{P_1}C_{P_2}\\
\times\sup_{\tau_1\in D(0,\rho_1),\tau_2\in D(0,\rho_2)}\left[
|\tau_1|^{k_1-1}|\tau_2|^{k_2-1}\left(\int_0^{|\tau_1|^{k_1}}(|\tau_1|^{k_1}-r_1)^{1/k_1-1}r_1^{1/k_1-1}dr_1\right)\right.\\
\left.\times\left(\int_0^{|\tau_2|^{k_2}}(|\tau_2|^{k_2}-r_2)^{1/k_2-1}r_2^{1/k_2-1}dr_2\right)\right]\sup_{m\in\R}\left(\int_{-\infty}^{\infty}\frac{(1+|m|)^{\mu-\gamma_1}}{(1+|m-m_1|)^{\mu-\hbox{deg}(P_1)}(1+|m_1|)^{\mu-\hbox{deg}(P_2)}}dm_1\right)\\
\times\left\|f\right\|_{(\beta,\mu,\rho_1,\rho_2)}\left\|g\right\|_{(\beta,\mu,\rho_1,\rho_2)}.
\end{multline*}
An analogous change of variable as the one performed in the proof of the previous lemma yields
$$\int_{0}^{x}(x-r)^{\frac{1}{k}-1}r^{1/k-1}dr=x^{2/k-1}\frac{\Gamma(1/k)\Gamma(1/k)}{\Gamma(2/k)},$$
for every $k>0$. The assumptions made on $\gamma_1$ and $\hbox{deg}(P_j)$, for $j=1,2$ allow us to conclude with 
$$\tilde{C}_2=C_1C_{P_1}C_{P_2}\rho_1\rho_2\frac{(\Gamma(1/k_1)\Gamma(1/k_2))^2}{\Gamma(2/k_1)\Gamma(2/k_2)}\sup_{m\in\R}\left(\int_{-\infty}^{\infty}\frac{(1+|m|)^{\mu-\gamma_1}}{(1+|m-m_1|)^{\mu-\hbox{deg}(P_1)}(1+|m_1|)^{\mu-\hbox{deg}(P_2)}}dm_1\right),$$
which is a finite quantity owing to Lemma~2.2 from~\cite{cota2}. 
\end{proof}

\textbf{Remark:} Observe that $\tilde{C}_2$ approaches 0 when the quantities $\rho_1$ or $\rho_2$ are reduced.

\vspace{0.3cm}

\subsubsection{Second auxiliary Banach space}\label{sec453}

Let $\beta>0$ and $\mu>1$. We fix 
$\nu_1,\nu_2>0$. Let $S_{d_j}$ be an infinite sector of positive opening, with vertex at the origin and bisecting direction $d_{j}\in\R$, for $j=1,2$. As in the previous section, we fix positive integers $k_1,k_2$.

\begin{defin}
Let us denote by $E_{(\beta,\mu,\nu_1,\nu_2,S_{d_1},S_{d_2})}$ the set of all continuous maps $(\tau_1,\tau_2,m)\mapsto h(\tau_1,\tau_2,m)$ on $S_{d_1}\times S_{d_2}\times \R$, which are holomorphic on its two first variables on $S_{d_1}\times S_{d_2}$ and such that there exists $C>0$  (depending on $\beta,\mu, S_{d_1}, S_{d_2}$) such that
$$|h(\tau_1,\tau_2,m)|\le C\frac{1}{(1+|m|)^{\mu}}e^{-\beta|m|}\frac{\left|\tau_1\right|}{1+\left|\tau_1\right|^{2k_1}}\frac{\left|\tau_2\right|}{1+\left|\tau_2\right|^{2k_2}} \exp\left(\nu_1\left|\tau_1\right|^{k_1}+\nu_2\left|\tau_2\right|^{k_2}\right),
$$
for every $(\tau_1,\tau_2,m)\in S_{d_1}\times S_{d_2}\times \R$. Given such $h$, the minimum of the constant $C$ above is denoted by $\left\|h(\tau_1,\tau_2,m)\right\|_{(\beta,\mu,\nu_1,\nu_2,S_{d_1},S_{d_2})}$. 
The pair $(E_{(\beta,\mu,\nu_1,\nu_2,S_{d_1},S_{d_2})},\left\|\cdot\right\|_{(\beta,\mu,\nu_1,\nu_2,S_{d_1},S_{d_2})})$ turns out to be a complex Banach space.
\end{defin}

As in the previous section, one can state some properties associated to the application of certain operators on the previous Banach space. The following is a direct consequence of the definition.

\begin{prop}\label{prop477b}
Let $(\tau_1,\tau_2,m)\mapsto b(\tau_1,\tau_2,m)$ be a continuous function defined on $S_{d_1}\times S_{d_2}\times \R$, holomorphic with respect to its first two variables on $S_{d_1}\times S_{d_2}$. Assume that 
$$C_b:=\sup_{(\tau_1,\tau_2,m)\in S_{d_1}\times S_{d_2}\times \R}|b(\tau_1,\tau_2,m)|$$
is finite. Then, for every $h\in E_{(\beta,\mu,\nu_1,\nu_2,S_{d_1},S_{d_2})}$, the function  $(\tau_1,\tau_2,m)\mapsto b(\tau_1,\tau_2,m)h(\tau_1,\tau_2,m)$ belongs to $E_{(\beta,\mu,\nu_1,\nu_2,S_{d_1},S_{d_2})}$ and it holds that
$$\left\|b(\tau_1,\tau_2,m)h(\tau_1,\tau_2,m)\right\|_{(\beta,\mu,\nu_1,\nu_2,S_{d_1},S_{d_2})}\le C_b 
\left\|h(\tau_1,\tau_2,m)\right\|_{(\beta,\mu,\nu_1,\nu_2,S_{d_1},S_{d_2})}.$$
\end{prop}

\begin{prop}\label{prop441b}
Let $a(\tau_1,\tau_2,m)$ be a continuous function defined on $S_{d_1}\times S_{d_2}\times\R$, holomorphic with respect to its two first variables on $S_{d_1}\times S_{d_2}$. Assume there exist $\gamma_1,\delta_1,\delta_2\ge0$ and $C_1>0$ such that 
$$|a(\tau_1,\tau_2,m)|\le \frac{C_1}{(1+|m|)^{\gamma_1}(1+|\tau_1|^{\delta_1k_1}|\tau_2|^{\delta_2k_2})},\quad (\tau_1,\tau_2,m)\in S_{d_1}\times S_{d_2}\times\R.$$

Let $m\mapsto h(m,\epsilon)$ be a function in $E_{(\beta,\mu)}$, such that for every $m\in\R$, the function $D(0,\epsilon_0)\ni \epsilon\mapsto h(m,\epsilon)$ is holomorphic on $D(0,\epsilon_0)$ and there exists $K>0$ with
$$\sup_{\epsilon\in D(0,\epsilon_0)}\left\|h(m,\epsilon)\right\|_{(\beta,\mu)}\le K.$$
We consider a polynomial $P(X)\in\C[X]$. We assume that 
$$\gamma_1\ge\hbox{deg}(P),\qquad \mu> \hbox{deg}(P)+1.$$

Let us also fix $\sigma_j>-1$ for $j=1,\ldots,6$.
We assume that 
\begin{multline}\label{e774}
k_1(\sigma_1+\sigma_3+\sigma_5+1)=k_2(\sigma_2+\sigma_4+\sigma_6+1)\\
\delta_1k_1=\delta_2k_2\qquad \sigma_1+\sigma_3+\sigma_5+1\le \delta_1.
\end{multline}

Then, for every $f\in E_{(\beta,\mu,\nu_1,\nu_2,S_{d_1},S_{d_2})}$ the function
\begin{multline*}\mathcal{B}_1(f):=a(\tau_1,\tau_2,m)\int_{-\infty}^{\infty}h(m-m_1,\epsilon)\tau_1^{\sigma_1k_1}\tau_2^{\sigma_2k_2}\int_{0}^{\tau_1^{k_1}}\int_{0}^{\tau_2^{k_2}}(\tau_1^{k_1}-s_1)^{\sigma_3}(\tau_2^{k_2}-s_2)^{\sigma_4}s_1^{\sigma_5}s_2^{\sigma_6}P(im_1)\\
\times f(s_1^{1/k_1},s_2^{1/k_2},m_1)ds_2ds_1dm_1
\end{multline*}
belongs to $E_{(\beta,\mu,\nu_1,\nu_2,S_{d_1},S_{d_2})}$. Moreover, there exists $\tilde{C}_1>0$ with
$$\left\|\mathcal{B}_1(f)\right\|_{(\beta,\mu,\nu_1,\nu_2,S_{d_1},S_{d_2})}\le K\tilde{C}_1 \left\|f\right\|_{(\beta,\mu,\nu_1,\nu_2,S_{d_1},S_{d_2})}.$$
\end{prop}
\begin{proof}
We only sketch the parts of the proof which differ from that of Proposition~\ref{prop441}. The first part of the proof can be mimicked to check that for every $(\tau_1,\tau_2,m)\in S_{d_1}\times S_{d_2}\times\R$,
\begin{multline*}
|\mathcal{B}_1(f)|(1+|m|)^{\mu}e^{\beta|m|}\frac{1+|\tau_1|^{2k_1}}{|\tau_1|}\frac{1+|\tau_2|^{2k_2}}{|\tau_2|}\exp\left(-\nu_1|\tau_1|^{k_1}-\nu_2|\tau_2|^{k_2}\right)\\
\le \frac{KC_1}{(1+|m|)^{\gamma_1}(1+|\tau_1|^{\delta_1k_1}|\tau_2|^{\delta_2k_2})}\int_{-\infty}^{\infty}\frac{1}{(1+|m-m_1|)^{\mu}}e^{-\beta|m-m_1|}\\
\times|\tau_1|^{\sigma_1k_1}|\tau_2|^{\sigma_2k_2}\int_{0}^{|\tau_1|^{k_1}}\int_{0}^{|\tau_2|^{k_2}}(|\tau_1|^{k_1}-r_1)^{\sigma_3}(|\tau_2|^{k_2}-r_2)^{\sigma_4}r_1^{\sigma_5}r_2^{\sigma_6}C_P(1+|m_1|)^{\hbox{deg}(P)}\\
\times\left(|f(r_1^{1/k_1}e^{i\hbox{arg}(\tau_1)},r_2^{1/k_2}e^{i\hbox{arg}(\tau_2)},m_1)| (1+|m_1|)^{\mu}e^{\beta|m_1|}\frac{1+\left(r_1^{1/k_1}\right)^{2k_1}}{r_1^{1/k_1}}\frac{1+\left(r_2^{1/k_2}\right)^{2k_2}}{r_2^{1/k_2}}\right.\\
\hfill\left.\times\exp\left(-\nu_1\left(r_1^{1/k_1}\right)^{k_1}-\nu_2\left(r_2^{1/k_2}\right)^{k_2}\right)\right)\\
\times \frac{1}{(1+|m_1|)^{\mu}}e^{-\beta|m_1|}\frac{r_1^{1/k_1}}{1+r_1^2}\frac{r_2^{1/k_2}}{1+r_2^2}\exp\left(\nu_1r_1+\nu_2r_2\right)dr_2dr_1dm_1\\
\times (1+|m|)^{\mu}e^{\beta|m|}\frac{1+\left|\tau_1\right|^{2k_1}}{\left|\tau_1\right|}\frac{1+\left|\tau_2\right|^{2k_2}}{\left|\tau_2\right|}\exp\left(-\nu_1\left|\tau_1\right|^{k_1}-\nu_2\left|\tau_2\right|^{k_2}\right).
\end{multline*}
The previous expression is upper bounded by
\begin{multline*}
KC_1C_P\left\|f\right\|_{(\beta,\mu,\nu_1,\nu_2,S_{d_1},S_{d_2})}\left[(1+|m|)^{\mu-\gamma_1}\int_{-\infty}^{\infty}\frac{1}{(1+|m-m_1|)^{\mu}(1+|m_1|)^{\mu-\hbox{deg}(P)}}dm_1\right]\\
\times \frac{|\tau_1|^{\sigma_1k_1}|\tau_2|^{\sigma_2k_2}}{1+|\tau_1|^{\delta_1k_1}|\tau_2|^{\delta_2k_2}} \frac{1+\left|\tau_1\right|^{2k_1}}{\left|\tau_1\right|}\frac{1+\left|\tau_2\right|^{2k_2}}{\left|\tau_2\right|}\exp\left(-\nu_1\left|\tau_1\right|^{k_1}-\nu_2\left|\tau_2\right|^{k_2}\right)\\
\times\int_{0}^{|\tau_1|^{k_1}}\int_{0}^{|\tau_2|^{k_2}}(|\tau_1|^{k_1}-r_1)^{\sigma_3}(|\tau_2|^{k_2}-r_2)^{\sigma_4}r_1^{\sigma_5}r_2^{\sigma_6}\frac{r_1^{1/k_1}}{1+r_1^2}\frac{r_2^{1/k_2}}{1+r_2^2}\exp\left(\nu_1r_1+\nu_2r_2\right)dr_2dr_1.
\end{multline*}
After the change of variable $r_j=|\tau_j|^{k_j}s_j$, for $j=1,2$, it only rests to upper estimate the following expression:
\begin{multline}\label{e797}
\frac{|\tau_1|^{k_1(\sigma_1+\sigma_3+\sigma_5+1)}|\tau_2|^{k_2(\sigma_2+\sigma_4+\sigma_6+1)}}{1+|\tau_1|^{\delta_1k_1}|\tau_2|^{\delta_2k_2}}\\
\times\left(\int_{0}^{1}(1-s_1)^{\sigma_3}s_1^{\sigma_5+\frac{1}{k_1}}\frac{1+\left|\tau_1\right|^{2k_1}}{1+|\tau_1|^{2k_1}s_1^2}\exp\left(\nu_1|\tau_1|^{k_1}(s_1-1)\right)ds_1\right)\\
\times \left(\int_{0}^{1}(1-s_2)^{\sigma_4}s_2^{\sigma_6+\frac{1}{k_2}}\frac{1+\left|\tau_2\right|^{2k_2}}{1+|\tau_2|^{2k_2}s_2^2}\exp\left(\nu_2|\tau_2|^{k_2}(s_2-1)\right)ds_2\right).
\end{multline}

We consider the function
$$(0,\infty)\times [0,1]\ni(x,r)\mapsto \varphi(x,r)=\frac{1+x^2}{1+x^2r^2}\exp(\nu_1x(r-1)).$$

We observe that given $0<r_0<1$, then 
$$\varphi(x,r)\le (1+x^2)\exp(\nu_1x(r_0-1))\le M_0(r_0),$$
for every $x>0$ and $0\le r\le r_0$. Let us assume now that $r>r_0$. Then, one has that $\varphi(x,r)\le 2\exp(\nu_1 x(r-1))\le 2$ if $x<1$. Since for $x\ge1$ and $1\ge r> r_0$, one has 
$$\frac{1+x^2}{1+x^2r^2}\le \frac{1+x^2}{1+x^2r_0^2}=\frac{1+\frac{1}{x^2}}{r_0^2+\frac{1}{x^2}}\le \frac{2}{r_0^2},$$
we deduce
$$\varphi(x,r)\le \frac{2}{r_0^2}\exp(\nu_1x(r-1))\le\frac{2}{r_0^2}.$$

We can apply this result to the estimates in (\ref{e797}) twice to achieve that it only rests to upper estimate the expression
$$\frac{|\tau_1|^{k_1(\sigma_1+\sigma_3+\sigma_5+1)}|\tau_2|^{k_2(\sigma_2+\sigma_4+\sigma_6+1)}}{1+|\tau_1|^{\delta_1k_1}|\tau_2|^{\delta_2k_2}}\left(\int_{0}^{1}(1-s_1)^{\sigma_3}s_1^{\sigma_5+\frac{1}{k_1}}ds_1\right)\left(\int_{0}^{1}(1-s_2)^{\sigma_4}s_2^{\sigma_6+\frac{1}{k_2}}ds_2\right).$$
Regarding the definition of Beta function, see~\cite{ba2} Appendix B, the first integral of the previous expression equals $B(\sigma_3+1,\sigma_5+\frac{1}{k_1}+1)$ and the second integral equals $B(\sigma_4+1,\sigma_6+\frac{1}{k_2}+1)$.
Finally, we prove that  
$$\sup_{(\tau_1,\tau_2)\in S_{d_1}\times S_{d_2}}\frac{|\tau_1|^{k_1(\sigma_1+\sigma_3+\sigma_5+1)}|\tau_2|^{k_2(\sigma_2+\sigma_4+\sigma_6+1)}}{1+|\tau_1|^{\delta_1k_1}|\tau_2|^{\delta_2k_2}}$$
is upper bounded. Indeed, let us write $\alpha=\frac{\delta_1}{\sigma_1+\sigma_3+\sigma_5+1}$ and $\beta= \frac{\delta_2}{\sigma_2+\sigma_4+\sigma_6+1}$. From the hypotheses of the parameters involved, described in (\ref{e774}), it is sufficient to prove that the function $g$ is upper bounded for $x,y\ge0$, with
$$g(x,y)=\frac{xy}{1+x^{\alpha}y^{\beta}},$$
with $\alpha=\beta\ge 1$. It is straight to check that the previous function is upper bounded for $x,y\ge0$.

The result follows.
\end{proof}

The proof of the following result follows analogous lines as that of Proposition~\ref{prop487}. We omit the repetitive parts.

\begin{prop}\label{prop487b}
Let $a(\tau_1,\tau_2,m)$ be a continuous function defined on $S_{d_1}\times S_{d_2}\times \R$, holomorphic with respect to its two first variables on $S_{d_1}\times S_{d_2}$. We assume such function satisfies there exist $\gamma_1,\delta_1,\delta_2\ge0$ and $C_1>0$ with
$$|a(\tau_1,\tau_2,m)|\le \frac{C_1}{(1+|m|)^{\gamma_1}(1+|\tau_1|^{\delta_1k_1}|\tau_2|^{\delta_2k_2})},\quad (\tau_1,\tau_2,m)\in S_{d_1}\times S_{d_2}\times \R.$$
We also assume that 
$$\delta_1k_1=\delta_2k_2\ge1.$$
Let $P_1(\epsilon,X),P_2(\epsilon,X)\in\mathcal{O}_b(D(0,\epsilon_0))[X]$ be polynomials with coefficients being holomorphic and bounded functions on $D(0,\epsilon_0)$. We assume that
$$\gamma_1\ge\max\{\hbox{deg}(P_1),\hbox{deg}(P_2)\},$$
and choose $\mu$ such that
$$\mu>\max\{\hbox{deg}(P_1),\hbox{deg}(P_2)\}+1.$$
For every $f,g\in E_{(\beta,\mu,\nu_1,\nu_2,S_{d_1},S_{d_2})}$, the function
\begin{multline*}
\mathcal{B}_2(f,g):=a(\tau_1,\tau_2,m)\tau_1^{k_1}\tau_2^{k_2}\int_{0}^{\tau_1^{k_1}}\int_{0}^{\tau_2^{k_2}}\int_{-\infty}^{\infty}P_1(\epsilon,i(m-m_1))f((\tau_1^{k_1}-s_1)^{1/k_1},(\tau_2^{k_2}-s_2)^{1/k_2},m-m_1)\\
\hfill\times P_2(\epsilon,im_1)g(s_1^{1/k_1},s_2^{1/k_2},m_1)\frac{1}{\tau_1^{k_1}-s_1}\frac{1}{s_1}\frac{1}{\tau_2^{k_2}-s_2}\frac{1}{s_2}dm_1ds_2ds_1
\end{multline*}
belongs to $E_{(\beta,\mu,\nu_1,\nu_2,S_{d_1},S_{d_2})}$. In addition to this, there exists $\tilde{C}_2>0$ such that
$$\left\|\mathcal{B}_2(f,g)\right\|_{(\beta,\mu,\nu_1,\nu_2,S_{d_1},S_{d_2})}\le \tilde{C}_2 \left\|f\right\|_{(\beta,\mu,\nu_1,\nu_2,S_{d_1},S_{d_2})}\left\|g\right\|_{(\beta,\mu,\nu_1,\nu_2,S_{d_1},S_{d_2})}.$$
\end{prop}

\begin{proof}
We have (\ref{e706}) holds for some $C_{P_1},C_{P_2}>0$, and all $m\in\R$. For all $(\tau_1,\tau_2,m)\in S_{d_1}\times S_{d_2}\times\R$ one has
\begin{multline*}
|\mathcal{B}_2(f,g)|(1+|m|)^{\mu}e^{\beta|m|}\frac{1+|\tau_1|^{2k_1}}{|\tau_1|}\frac{1+|\tau_2|^{2k_2}}{|\tau_2|}\exp\left(-\nu_1|\tau_1|^{k_1}-\nu_2|\tau_2|^{k_2}\right)\\
\le\frac{C_1|\tau_1|^{k_1}|\tau_2|^{k_2}}{(1+|m|)^{\gamma_1}(1+|\tau_1|^{\delta_1k_1}|\tau_2|^{\delta_2k_2})}\\
\times\int_{0}^{|\tau_1|^{k_1}}\int_{0}^{|\tau_2|^{k_2}}\int_{-\infty}^{\infty}C_{P_1}(1+|m-m_1|)^{\hbox{deg}(P_1)}C_{P_2}(1+|m_1|)^{\hbox{deg}(P_2)}\\
\times\left\{|f((\tau_1^{k_1}-s_1e^{ik_1\hbox{arg}(\tau_1)})^{1/k_1},(\tau_2^{k_2}-s_2e^{ik_2\hbox{arg}(\tau_2)})^{1/k_2},m-m_1)|(1+|m-m_1|)^{\mu}e^{\beta|m-m_1|}\right.\\
\times\frac{1+\left(|\tau_1^{k_1}-s_1e^{ik_1\hbox{arg}(\tau_1)}|^{1/k_1}\right)^{2k_1}}{|\tau_1^{k_1}-s_1e^{ik_1\hbox{arg}(\tau_1)}|^{1/k_1}}\frac{1+\left(|\tau_2^{k_2}-s_2e^{ik_2\hbox{arg}(\tau_2)}|^{1/k_2}\right)^{2k_2}}{|\tau_2^{k_2}-s_2e^{ik_2\hbox{arg}(\tau_2)}|^{1/k_2}}\\
\left.\times
\exp\left(-\nu_1|\tau_1^{k_1}-s_1e^{ik_1\hbox{arg}(\tau_1)}|-\nu_2|\tau_2^{k_2}-s_2e^{ik_2\hbox{arg}(\tau_2)}|\right)\right\}\\
\times\frac{e^{-\beta|m-m_1|}}{(1+|m-m_1|)^{\mu}}
\frac{|\tau_1^{k_1}-s_1e^{ik_1\hbox{arg}(\tau_1)}|^{1/k_1}}{1+\left(|\tau_1^{k_1}-s_1e^{ik_1\hbox{arg}(\tau_1)}|^{1/k_1}\right)^{2k_1}}\frac{|\tau_2^{k_2}-s_2e^{ik_2\hbox{arg}(\tau_2)}|^{1/k_2}}{1+\left(|\tau_2^{k_2}-s_2e^{ik_2\hbox{arg}(\tau_2)}|^{1/k_2}\right)^{2k_2}}\\
\times \exp\left(\nu_1|\tau_1^{k_1}-s_1e^{ik_1\hbox{arg}(\tau_1)}|+\nu_2|\tau_2^{k_2}-s_2e^{ik_2\hbox{arg}(\tau_2)}|\right)\\
\times\left[|g((s_1e^{ik_1\hbox{arg}(\tau_1)})^{1/k_1},(s_2e^{ik_2\hbox{arg}(\tau_2)})^{1/k_2},m_1)(1+|m_1|)^{\mu}e^{\beta|m_1|}
\frac{1+(s_1e^{ik_1\hbox{arg}(\tau_1)})^{2}}{(s_1e^{ik_1\hbox{arg}(\tau_1)})^{1/k_1}}\right.\\
\left.\times \frac{1+(s_2e^{ik_2\hbox{arg}(\tau_2)})^{2}}{(s_2e^{ik_2\hbox{arg}(\tau_2)})^{1/k_2}}\exp\left(-\nu_1s_1e^{ik_1\hbox{arg}(\tau_1)}-\nu_2s_2e^{ik_2\hbox{arg}(\tau_2)}\right)|\right]\\
\times \frac{1}{(1+|m_1|)^{\mu}}e^{-\beta|m_1|}\frac{(s_1)^{1/k_1}}{1+(s_1)^{2}}\frac{(s_2)^{1/k_2}}{1+(s_2)^{2}}\exp\left(\nu_1s_1+\nu_2s_2\right)\frac{1}{|\tau_1^{k_1}|-s_1}\frac{1}{s_1}\\
\times\frac{1}{|\tau_2^{k_2}|-s_2}\frac{1}{s_2}dm_1ds_2ds_1 (1+|m|)^{\mu}e^{\beta|m|}\frac{1+\left|\tau_1\right|^{2k_1}}{\left|\tau_1\right|}\frac{1+\left|\tau_2\right|^{2k_2}}{\left|\tau_2\right|}\exp\left(-\nu_1\left|\tau_1\right|^{k_1}-\nu_2\left|\tau_2\right|^{k_2}\right)
\end{multline*}
which leads us to
\begin{multline*}
|\mathcal{B}_2(f,g)|(1+|m|)^{\mu}e^{\beta|m|}\frac{1+|\tau_1|^{2k_1}}{|\tau_1|}\frac{1+|\tau_2|^{2k_2}}{|\tau_2|}\exp\left(-\nu_1|\tau_1|^{k_1}-\nu_2|\tau_2|^{k_2}\right)\\
\le
C_1C_{P_1}C_{P_2}\left\|f\right\|_{(\beta,\mu,\nu_1,\nu_2,S_{d_1},S_{d_2})}\left\|g\right\|_{(\beta,\mu,\nu_1,\nu_2,S_{d_1},S_{d_2})}\\
\times \left[\int_{-\infty}^{\infty}\frac{(1+|m|)^{\mu-\gamma_1}dm_1}{(1+|m-m_1|)^{\mu-\hbox{deg}(P_1)}(1+|m_1|)^{\mu-\hbox{deg}(P_2)}}\right]\\
\times\frac{|\tau_1|^{k_1}|\tau_2|^{k_2}}{1+|\tau_1|^{\delta_1k_1}|\tau_2|^{\delta_2k_2}}\left[\int_{0}^{|\tau_1|^{k_1}}\frac{(|\tau_1^{k_1}|-s_1)^{1/k_1}}{1+\left((|\tau_1^{k_1}|-s_1)^{1/k_1}\right)^{2k_1}}\frac{s_1^{1/k_1}}{1+s_1^{2}}\frac{1}{|\tau_1^{k_1}|-s_1}\frac{1}{s_1}ds_1\right]\\
\times\left[\int_{0}^{|\tau_2|^{k_2}}\frac{(|\tau_2^{k_2}|-s_2)^{1/k_2}}{1+\left((|\tau_2^{k_2}|-s_2)^{1/k_2}\right)^{2k_2}}\frac{s_2^{1/k_2}}{1+s_2^{2}}\frac{1}{|\tau_2|^{k_2}-s_2}\frac{1}{s_2}ds_2\right]\frac{1+\left|\tau_1\right|^{2k_1}}{\left|\tau_1\right|}\frac{1+\left|\tau_2\right|^{2k_2}}{\left|\tau_2\right|}.
\end{multline*}
The change of variables $s_j=|\tau_j|^{k_j}r_j$, $j=1,2$ allow us to write

\begin{multline*}
|\mathcal{B}_2(f,g)|(1+|m|)^{\mu}e^{\beta|m|}\frac{1+|\tau_1|^{2k_1}}{|\tau_1|}\frac{1+|\tau_2|^{2k_2}}{|\tau_2|}\exp\left(-\nu_1|\tau_1|^{k_1}-\nu_2|\tau_2|^{k_2}\right)\\
\le
C_1C_{P_1}C_{P_2}\left\|f\right\|_{(\beta,\mu,\nu_1,\nu_2,S_{d_1},S_{d_2})}\left\|g\right\|_{(\beta,\mu,\nu_1,\nu_2,S_{d_1},S_{d_2})}\\
\times \left[\int_{-\infty}^{\infty}\frac{(1+|m|)^{\mu-\gamma_1}dm_1}{(1+|m-m_1|)^{\mu-\hbox{deg}(P_1)}(1+|m_1|)^{\mu-\hbox{deg}(P_2)}}\right]\\
\times\frac{|\tau_1||\tau_2|}{1+|\tau_1|^{\delta_1k_1}|\tau_2|^{\delta_2k_2}}\left[\int_{0}^{1}(1-r_1)^{1/k_1-1}r_1^{1/k_1-1}\frac{1+\left|\tau_1\right|^{2k_1}}{\left(1+r_1^2\left|\tau_1\right|^{2k_1}\right)\left(1+(1-r_1)^2\left|\tau_1\right|^{2k_1}\right)}dr_1\right]\\
\times \left[\int_{0}^{1}(1-r_2)^{1/k_2-1}r_2^{1/k_2-1}\frac{1+\left|\tau_2\right|^{2k_2}}{\left(1+r_2^2\left|\tau_2\right|^{2k_2}\right)\left(1+(1-r_2)^2\left|\tau_2\right|^{2k_2}\right)}dr_2\right]
\end{multline*}

Let $k\ge 1$. We consider the function
$$h(r,x)=\frac{1+x^{2k}}{(1+r^2x^{2k})(1+(1-r)^2x^{2k})},\quad (r,x)\in [0,1]\times \{x\in\R:x\ge0\}.$$
A partial fraction decomposition of $r\mapsto h(r,x)$ yields
$$h(r,x)=(1+x^{2k})\frac{1}{x^{2k}+4}\left[\frac{3-2r}{1+x^{2k}(1-r)^2}+\frac{2r+1}{1+x^{2kr^2}}\right]\le 3\frac{1+x^{2k}}{4+x^{2k}},$$
for all $0\le r\le 1$ and $x\ge0$, which implies that $h(r,x)$ is bounded on the domain $[0,1]\times [0, \infty)$.
This fact applied twice allows us to arrive at the existence of a positive constant $C_2$ such that

\begin{multline*}
|\mathcal{B}_2(f,g)|(1+|m|)^{\mu}e^{\beta|m|}\frac{1+|\tau_1|^{2k_1}}{|\tau_1|}\frac{1+|\tau_2|^{2k_2}}{|\tau_2|}\exp\left(-\nu_1|\tau_1|^{k_1}-\nu_2|\tau_2|^{k_2}\right)\\
\le
C_1C_{P_1}C_{P_2}C_2\frac{(\Gamma(1/k_1)\Gamma(1/k_2))^2}{\Gamma(2/k_1)\Gamma(2/k_2)}\left\|f\right\|_{(\beta,\mu,\nu_1,\nu_2,S_{d_1},S_{d_2})}\left\|g\right\|_{(\beta,\mu,\nu_1,\nu_2,S_{d_1},S_{d_2})}\\
\times \left[\int_{-\infty}^{\infty}\frac{(1+|m|)^{\mu-\gamma_1}dm_1}{(1+|m-m_1|)^{\mu-\hbox{deg}(P_1)}(1+|m_1|)^{\mu-\hbox{deg}(P_2)}}\right]\frac{|\tau_1||\tau_2|}{1+|\tau_1|^{\delta_1k_1}|\tau_2|^{\delta_2k_2}}.
\end{multline*}

We conclude the result by taking into account analogous estimates as in Proposition~\ref{prop441b} regarding both the first integral remaining, and the fact that
$$\sup_{x,y\ge0}\frac{xy}{1+x^{\alpha}y^{\beta}}<\infty,$$
for $\alpha=\beta\ge1$.

\end{proof}

\vspace{0.3cm}

\textbf{Remark:}
The constant $\tilde{C}_2$ approaches to 0 provided that the quantities $C_{P_1},C_{P_2}>0$ are taken small enough. 

\subsubsection{Third auxiliary Banach space}\label{sec924}

Let $\beta>0$ and $\mu>1$. We fix $\rho_1,\rho_2>0$, and $\nu_1,\nu_2>0$. Let $S_{d_j}$ be an infinite sector of some positive opening, with vertex at the origin, and bisecting direction $d_{j}\in\R$, for $j=1,2$. As in the previous sections, we fix positive integers $k_1,k_2$.

\begin{defin}
Let us denote by $F_{(\beta,\mu,\rho_1,\rho_2,S_{d_1},S_{d_2})}$ the set of all continuous maps $(\tau_1,\tau_2,m)\mapsto h(\tau_1,\tau_2,m)$ on $(S_{d_1}\cap D(0,\rho_1))\times( S_{d_2}\cap D(0,\rho_2))\times \R$, which are holomorphic on its two first variables on $(S_{d_1}\cap D(0,\rho_1))\times (S_{d_2}\cap D(0,\rho_2))$ and such that there exists $C>0$  (depending on $\beta,\mu, \rho_1,\rho_2, S_{d_1}, S_{d_2}$) such that
$$|h(\tau_1,\tau_2,m)|\le C\frac{1}{(1+|m|)^{\mu}}e^{-\beta|m|}\left|\tau_1\tau_2\right|,
$$
for every $(\tau_1,\tau_2,m)\in (S_{d_1}\cap D(0,\rho_1))\times (S_{d_2}\cap D(0,\rho_2))\times \R$. For such $h$, the minimum of such constant $C$ is denoted by $\left\|h(\tau_1,\tau_2,m)\right\|_{(\beta,\mu,\rho_1,\rho_2,S_{d_1},S_{d_2})}$. 
The pair $(F_{(\beta,\mu,\rho_1,\rho_2,S_{d_1},S_{d_2})},\left\|\cdot\right\|_{(\beta,\mu,\rho_1,\rho_2,S_{d_1},S_{d_2})})$ is a complex Banach space.
\end{defin}

Analogous results as those in the two previous sections regarding the action of certain operators acting on functions belonging to this Banach space can be stated. We omit their proof which follow analogous arguments as before.

\begin{prop}\label{prop477c}
Let $(\tau_1,\tau_2,m)\mapsto b(\tau_1,\tau_2,m)$ be a continuous function defined on $(S_{d_1}\cap D(0,\rho_1))\times( S_{d_2}\cap D(0,\rho_2))\times \R$, holomorphic with respect to its first two variables on $(S_{d_1}\cap D(0,\rho_1))\times(S_{d_2}\cap D(0,\rho_2))$. Assume that 
$$C_b:=\sup_{(\tau_1,\tau_2,m)\in (S_{d_1}\cap D(0,\rho_1))\times (S_{d_2}\cap D(0,\rho_2))\times \R}|b(\tau_1,\tau_2,m)|$$
is finite. Then, for every $h\in F_{(\beta,\mu,\rho_1,\rho_2,S_{d_1},S_{d_2})}$, the function  $(\tau_1,\tau_2,m)\mapsto b(\tau_1,\tau_2,m)h(\tau_1,\tau_2,m)$ belongs to $F_{(\beta,\mu,\rho_1,\rho_2,S_{d_1},S_{d_2})}$ and it holds that
$$\left\|b(\tau_1,\tau_2,m)h(\tau_1,\tau_2,m)\right\|_{(\beta,\mu,\rho_1,\rho_2,S_{d_1},S_{d_2})}\le C_b \left\|h(\tau_1,\tau_2,m)\right\|_{(\beta,\mu,\rho_1,\rho_2,S_{d_1},S_{d_2})}.$$
\end{prop}

\begin{prop}\label{prop441c}
Let $a(\tau_1,\tau_2,m)$ be a continuous function defined on $(S_{d_1}\cap D(0,\rho_1))\times (S_{d_2}\cap D(0,\rho_2))\times\R$, holomorphic with respect to its two first variables on $(S_{d_1}\cap D(0,\rho_1))\times (S_{d_2}\cap D(0,\rho_2))$. Assume there exist $\gamma_1\ge0$ and $C_1>0$ such that 
$$|a(\tau_1,\tau_2,m)|\le \frac{C_1}{(1+|m|)^{\gamma_1}},\quad (\tau_1,\tau_2,m)\in (S_{d_1}\cap D(0,\rho_1))\times (S_{d_2}\cap D(0,\rho_2))\times\R.$$

Let $m\mapsto h(m,\epsilon)$ be a function in $E_{(\beta,\mu)}$, such that for every $m\in\R$, the function $D(0,\epsilon_0)\ni \epsilon\mapsto h(m,\epsilon)$ is holomorphic on $D(0,\epsilon_0)$ and there exists $K>0$ with
$$\sup_{\epsilon\in D(0,\epsilon_0)}\left\|h(m,\epsilon)\right\|_{(\beta,\mu)}\le K.$$
We consider a polynomial $P(X)\in\C[X]$. We assume that 
$$\gamma_1\ge\hbox{deg}(P),\qquad \mu> \hbox{deg}(P)+1.$$

Let us also fix $\sigma_j>-1$ for $j=1,\ldots,6$, such that
$$k_1\sigma_1+\sigma_3+\sigma_5+\frac{1}{k_1}\ge0,\qquad k_2\sigma_2+\sigma_4+\sigma_6+\frac{1}{k_2}\ge0.$$

Then, for every $f\in F_{(\beta,\mu,\rho_1,\rho_2,S_{d_1},S_{d_2})}$ the function
\begin{multline*}\mathcal{B}_1(f):=a(\tau_1,\tau_2,m)\int_{-\infty}^{\infty}h(m-m_1,\epsilon)\tau_1^{\sigma_1k_1}\tau_2^{\sigma_2k_2}\int_{0}^{\tau_1^{k_1}}\int_{0}^{\tau_2^{k_2}}(\tau_1^{k_1}-s_1)^{\sigma_3}(\tau_2^{k_2}-s_2)^{\sigma_4}s_1^{\sigma_5}s_2^{\sigma_6}P(im_1)\\
\times f(s_1^{1/k_1},s_2^{1/k_2},m_1)ds_2ds_1dm_1
\end{multline*}
belongs to $F_{(\beta,\mu,\rho_1,\rho_2,S_{d_1},S_{d_2})}$. Moreover, there exists $\tilde{C}_1>0$ with
$$\left\|\mathcal{B}_1(f)\right\|_{(\beta,\mu,\rho_1,\rho_2,S_{d_1},S_{d_2})}\le K\tilde{C}_1 \left\|f\right\|_{(\beta,\mu,\rho_1,\rho_2,S_{d_1},S_{d_2})}.$$
\end{prop}
\begin{prop}\label{prop487c}
Let $a(\tau_1,\tau_2,m)$ be a continuous function defined on $(S_{d_1}\cap D(0,\rho_1))\times (S_{d_2}\cap D(0,\rho_2))\times \R$, holomorphic with respect to its two first variables on $(S_{d_1}\cap D(0,\rho_1))\times (S_{d_2}\cap D(0,\rho_2))$. We assume such function satisfies there exists $\gamma_1\ge0$ and $C_1>0$ with
$$|a(\tau_1,\tau_2,m)|\le \frac{C_1}{(1+|m|)^{\gamma_1}},\quad (\tau_1,\tau_2,m)\in (S_{d_1}\cap D(0,\rho_1))\times (S_{d_2}\cap D(0,\rho_2))\times \R.$$
Let $P_1(\epsilon,X),P_2(\epsilon,X)\in\mathcal{O}_b(D(0,\epsilon_0))[X]$ be polynomials with coefficients being holomorphic and bounded functions on $D(0,\epsilon_0)$. We assume that
$$\gamma_1\ge\max\{\hbox{deg}(P_1),\hbox{deg}(P_2)\},$$
and choose $\mu$ such that
$$\mu>\max\{\hbox{deg}(P_1),\hbox{deg}(P_2)\}+1.$$
For every $f,g\in F_{(\beta,\mu,\rho_1,\rho_2,S_{d_1},S_{d_2})}$, the function
\begin{multline*}
\mathcal{B}_2(f,g):=a(\tau_1,\tau_2,m)\tau_1^{k_1}\tau_2^{k_2}\int_{0}^{\tau_1^{k_1}}\int_{0}^{\tau_2^{k_2}}\int_{-\infty}^{\infty}P_1(\epsilon,i(m-m_1))f((\tau_1^{k_1}-s_1)^{1/k_1},(\tau_2^{k_2}-s_2)^{1/k_2},m-m_1)\\
\hfill\times P_2(\epsilon,im_1)g(s_1^{1/k_1},s_2^{1/k_2},m_1)\frac{1}{\tau_1^{k_1}-s_1}\frac{1}{s_1}\frac{1}{\tau_2^{k_2}-s_2}\frac{1}{s_2}dm_1ds_2ds_1
\end{multline*}
belongs to $F_{(\beta,\mu,\rho_1,\rho_2,S_{d_1},S_{d_2})}$. In addition to this, there exists $\tilde{C}_2>0$ such that
$$\left\|\mathcal{B}_2(f,g)\right\|_{(\beta,\mu,\nu_1,\nu_2,S_{d_1},S_{d_2})}\le \tilde{C}_2 \left\|f\right\|_{(\beta,\mu,\nu_1,\nu_2,S_{d_1},S_{d_2})}\left\|g\right\|_{(\beta,\mu,\nu_1,\nu_2,S_{d_1},S_{d_2})}.$$
\end{prop}

\vspace{0.3cm}

\textbf{Remark:} Observe the constant $\tilde{C}_2$ approaches 0 when the quantities $\rho_1$ or $\rho_2$ are reduced.

\vspace{0.4cm}

\textbf{Aknowledgements:} The second and third authors are partially supported by the project PID2022-139631NB-I00 of Ministerio de Ciencia e Innovaci\'on, Spain.


\begin{thebibliography}{99}
\bibitem{ba97} W. Balser, \emph{Multisummability of complete formal solutions for non-linear systems of meromorphic ordinary differential equations.} Complex Var. Theory Appl. 34 (1-2) (1997) 19--24. 
\bibitem{ba2} W. Balser, \emph{Formal power series and linear systems of meromorphic ordinary differential equations.} Universitext. Springer-Verlag, New York, 2000. xviii+299 pp.
\bibitem{ba4} W. Balser, \emph{Multisummability of formal power series solutions of partial differential equations with constant coefficients.} J. Differential Equations 201 (2004), no. 1, 63--74.
\bibitem{babrrasi} W. Balser, B. Braaksma, J.-P. Ramis, Y. Sibuya, \emph{Multisummability of formal power series solutions of linear ordinary differential equations.} Asymptot. Anal. 5(1) (1991) 27--45. 
\bibitem{br92} B. Braaksma, \emph{Multisummability of formal power series solutions of nonlinear meromorphic differential equations.} Ann. Inst. Fourier (Grenoble) 42(3) (1992) 517--540.
\bibitem{chenlastramalek} G. Chen, A. Lastra, S. Malek, \emph{Parametric Gevrey asymptotics in two complex time variables through truncated Laplace transforms.} Adv. Differ. Equ. 2020, 307 (2020). 
\bibitem{cota2} O. Costin, S. Tanveer, \emph{Short time existence and Borel summability in the Navier-Stokes equation in $\mathbb{R}^{3}$}, Comm. Partial Differential Equations 34 (2009), no. 7-9, 785--817.
\bibitem{hssi} P. Hsieh, Y. Sibuya, \emph{Basic theory of ordinary differential equations}. Universitext. Springer-Verlag, New York, 1999.
\bibitem{jikalasa} J. Jim\'enez-Garrido, S. Kamimoto, A. Lastra, J. Sanz, \emph{Multisummability in Carleman ultraholomorphic classes by means of nonzero proximate orders,} J. Math. Anal. Appl. 472, No. 1 (2019) 627--686. 
\bibitem{family1} A. Lastra, S. Malek, \emph{On parametric Gevrey asymptotics for some nonlinear initial value problems in symmetric complex time variables,} Asymptotic Anal. 118 (2020) No. 1--2, 49--79. 
\bibitem{family2} A. Lastra, S. Malek, \emph{On parametric Gevrey asymptotics for some initial value problems in two asymmetric complex time variables,} Results Math. 73 (2018), no. 4, Art. 155, 46 pp.
\bibitem{lama} A. Lastra, S. Malek, \emph{On parametric Gevrey asymptotics for some nonlinear initial value Cauchy problems.} J. Differential Equations 259 (2015), no. 10, 5220--5270.
\bibitem{lama1} A. Lastra, S. Malek, \emph{On parametric multisummable formal solutions to some nonlinear initial value Cauchy problems.} Adv. Difference Equ. 2015, 2015:200, 78 pp.
\bibitem{lamaboundarylayer} A. Lastra, S. Malek, \emph{Multiscale Gevrey asymptotics in boundary layer expansions for some initial value problem with merging turning points.} Adv. Differential Equations 24 no.1--2 (2019), 69--136. 
\bibitem{lamaboundarylayer2} A. Lastra, S. Malek, \emph{Boundary layer expansions for initial value problems with two complex time variables.} Adv. Difference Equations 2020, Paper No. 20 (2020).
\bibitem{lamisu4} A. Lastra, S. Michalik, M. Suwi\'nska, \emph{Multisummability of formal solutions ofr a family of generalized singularly perturbed moment differential equations.} Results Math. 78, 49 (2023). 
\bibitem{lo94} M. Loday-Richaud, \emph{Stokes phenomenon, multisummability and differential Galois groups.} Ann. Inst. Fourier (Grenoble) 44(3) (1994) 849--906. 
\bibitem{loday} M. Loday-Richaud, Divergent series, summability and resurgence. II. Simple and multiple summability. Lecture Notes in Mathematics, 2154. Springer, 2016. 
\bibitem{malek20} S. Malek, \emph{On a partial $q-$analog of a singularly perturbed problem with fuchsian and irregular time singularities}, Abstract and Applied Analysis, vol. 2020 (2020) Article ID 7985298. 
\bibitem{ma23} S. Malek, \emph{Double-scale expansions for a logarithmic type solution to a q-analog of a singular initial value problem}. Abstract and Applied Analysis, Volume 2023 (2023) Article ID 3025513. 
\bibitem{malra92} B. Malgrange, J.-P. Ramis, \emph{Fonctions multisommables.} Ann. Inst. Fourier (Grenoble) 42(1-2) (1992) 353--368.  
\bibitem{mi10} S. Michalik, \emph{On the multisummability of divergent solutions of linear partial differential equations with constant coefficients.} J. Differential Equations 249(3) (2010) 551--570.
\bibitem{michalik12} S. Michalik, \emph{Multisummability of formal solutions of inhomogeneous linear partial differential equations with constant coefficients}, J. Dyn. Control Syst. 18 (2012) 103--133.
\bibitem{ou} S. Ouchi, \emph{Multisummability of formal power series solutions of nonlinear partial differential equations in complex domains}. Asymptotic Anal. 47 (2006), no. 3-4, 187--225.
\bibitem{rasi94} J.-P. Ramis, Y. Sibuya, \emph{A new proof of multisummability of formal solutions of nonlinear meromorphic differential equations.} Ann. Inst. Fourier (Grenoble) 44(3) (1994) 811-848. 

\bibitem{ta} H. Tahara, \emph{Asymptotic existence theorem for formal solutions with singularities of nonlinear partial differential equations via multisummability.} J. Math. Soc. Japan 75 (2023) No. 3, 1055--1117.  
\bibitem{taya} H. Tahara, H. Yamazawa, \emph{Multisummability of formal solutions to the Cauchy problem for some linear partial differential equations}, J. Differential Equations, Volume 255, Issue 10, 15 November 2013, pages 3592--3637.
\bibitem{takei} Y. Takei, \emph{On the multisummability of WKB solutions of certain singularly perturbed linear ordinary differential equations}, Opuscula Math. 35 no. 5 (2015), 775--802.
\bibitem{ya} H. Yamazawa, \emph{On multisummability of formal solutions with logarithmic terms for some linear partial differential equations}, Funkc. Ekvacioj, 60, No. 3 (2017), 371--406.
\bibitem{yo} M. Yoshino, \emph{Parametric Borel summability of partial differential equations of irregular singular type}. Analytic, algebraic and geometric aspects of differential equations, 455–471, Trends Math., Birkh\"auser/Springer, Cham, 2017.
\end{thebibliography}
\end{document}